# Partial-State Deadzone-Adapted Disturbance Suppression Control For Strict-Feedback Systems

**Iasson Karafyllis[*], Alexandros Aslanidis[*] and Miroslav Krstic[**]**

[*]Dept. of Mathematics, National Technical University of Athens, Zografou Campus, 15780, Athens, Greece, email: iasonkar@central.ntua.gr; alex.aslanidis@hotmail.com

[**]Dept. of Mechanical and Aerospace Eng., University of California, San Diego, La Jolla, CA 92093-0411, U.S.A., email: krstic@ucsd.edu

## Abstract

In this paper we provide an extension of the partial-state Deadzone-Adapted Disturbance Suppression (DADS) adaptive control scheme in the case of strict-feedback systems. We study two cases: the cases of strict-feedback systems with and without integrators. The DADS controller achieves attenuation of the output to an assignable small level, despite the presence of (time-varying) disturbances, unknown (time-varying) parameters and unknown (time-varying) input coefficients; all of them with arbitrary and unknown bounds. The mechanism that makes DADS work is the same in all previously studied cases where DADS has been applied: the update law "pumps up" the nonlinear controller gain to a sufficient size, adapting the gain to the unknown uncertainties, and shuts off further pumping up of the gain once the output is in the residual set. The results are directly extended to the infinite-dimensional case. Two simple examples are provided: one finite-dimensional example and one infinite-dimensional example. The finite-dimensional example is used for a detailed comparison with existing results in the literature.

## 1. Introduction

Robustness in adaptive control is a topic that has attracted the attention of researchers in control theory. It is far more difficult to achieve robustness with respect to dynamic uncertainties compared to static uncertainties because dynamic uncertainties describe unmodeled dynamic phenomena, which typically have strong influence on the closed loop system performance due to the presence of feedback interconnections of unknown magnitude with the plant itself. The problem of achieving robustness with respect to dynamic uncertainties is equivalent to the problem of achieving robustness by means of a partial-state adaptive feedback law. The partial-state feedback design problem is important due to implementation issues (lack of measurements of some state components) but also due to the standard engineering practice: the feedback law is usually designed for a simplified system where a small or large part of its dynamics has been neglected.

Various approaches have been proposed in the literature for achieving robustness with respect to dynamic uncertainties in time-invariant nonlinear control systems for which no persistence of excitation condition is assumed. The nonlinear indirect adaptive laws in [16] produced by a backstepping procedure were modified by Praly and Jiang in [6] by using leakage (see [3, 4, 17]) together with dynamic (high) gains and nonlinear damping and by employing small-gain arguments (see [5]). The Praly-Jiang approach was later extended in [19, 20, 22, 23, 26]. A different approach, the universal $\lambda -$tracking approach, was used in [25] for the design of partial-state adaptive feedback laws. On the other hand, starting from [1], the prescribed performance adaptive control approach was utilized for the design of partial-state, time-varying, direct adaptive feedback laws. The prescribed performance adaptive control approach was later extended in [2, 24, 27, 28].

Recently, the Deadzone-Adapted Disturbance Suppression (DADS) adaptive control scheme (see [9, 10, 12]) provided a simple, direct, adaptive control scheme that combines nonlinear damping, single-gain adjustment (the dynamic feedback has only one state), and deadzone in the update law. The DADS controller achieves attenuation of the plant state to an assignable small level, despite the presence of (time-varying) disturbances and unknown (time-varying) parameters of arbitrary and unknown bounds. Moreover, the DADS controller prevents gain and state drift regardless of the size of the disturbance and unknown parameter. The DADS controller was recently applied to voltage regulation problems in [21] and was extended to the partial-state feedback case in [11, 13] under the assumption of the matching condition. It was shown in [11, 13] that DADS can work even if the neglected dynamics are infinite-dimensional, e.g., even if Partial Differential Equations (PDEs) are neglected. It was also shown in [11, 13] that DADS can allow an "escape" of the requirements imposed by using small-gain arguments.

In the present paper, we provide the extension of the partial-state DADS adaptive control scheme in the case of strict-feedback systems. We study two cases: the cases of strict-feedback systems with and without integrators. The DADS controller achieves attenuation of the output to an assignable small level, despite the presence of (time-varying) disturbances, unknown (time-varying) parameters and unknown (time-varying) input coefficients; all of them with arbitrary and unknown bounds. The mechanism that makes DADS work is the same in all previously studied cases where DADS was applied (see [9, 10, 12]): the update law "pumps up" the nonlinear controller gain to a sufficient size, adapting the gain to the unknown uncertainties, and shuts off further pumping up of the gain once the output is in the residual set. The results are directly extended to the infinite-dimensional case.

The structure of the paper is as follows. In Section 2 we provide all the results (Theorem 1 and Theorem 2) for the finite-dimensional case. In Section 3 we provide all results (Theorem 3 and Theorem 4) for the infinite-dimensional case. Section 4 is devoted to the presentation of two simple examples (one finite-dimensional and one infinite-dimensional). The finite-dimensional example is used for a detailed comparison with the controller proposed in [6]. The proofs of all results are provided in Section 5. Finally, the concluding remarks of the present work are provided in Section 6.

**Notation.** Throughout this paper, we adopt the following notation.

* $\mathbb{R}_+ := [0,+\infty)$. For a vector $x \in \mathbb{R}^n$, $|x|$ denotes its Euclidean norm and $x'$ denotes its transpose. We use the notation $x^+$ for the positive part of the real number $x \in \mathbb{R}$, i.e., $x^+ = \max(x,0)$.

* Let $D \subseteq \mathbb{R}^n$ be an open set and let $S \subseteq \mathbb{R}^n$ be a set that satisfies $D \subseteq S \subseteq cl(D)$, where $cl(D)$ is the closure of $D$. By $C^0(S\,;\Omega)$, we denote the class of continuous functions on $S$, which take values in $\Omega \subseteq \mathbb{R}^m$. By $C^k(S\,;\Omega)$, where $k \geq 1$ is an integer, we denote the class of functions on $S \subseteq \mathbb{R}^n$, which take values in $\Omega \subseteq \mathbb{R}^m$ and have continuous derivatives of order $k$. In other words, the functions of class $C^k(S\,;\Omega)$ are the functions which have continuous derivatives of order $k$ in $D = \mathrm{int}(S)$ that can be continued continuously to all points in $\partial D \cap S$. When $\Omega = \mathbb{R}$ then we write $C^0(S)$ or $C^k(S)$. A function $f \in C^\infty(S\,;\Omega) = \bigcap_{k=0}^{\infty} C^k(S\,;\Omega)$ is called a smooth function.

* By $L^\infty\left(\mathbb{R}_+;\Omega\right)$, where $\Omega \subseteq \mathbb{R}^m$ we denote the class of essentially bounded, Lebesgue measurable functions $d:\mathbb{R}_+ \to \Omega$. For $d \in L^\infty\left(\mathbb{R}_+;\Omega\right)$ we define $\|d\|_\infty = \sup_{t\geq 0}\left(|d(t)|\right)$, where $\sup_{t\geq 0}\left(|d(t)|\right)$ is the essential supremum. When $\Omega = \mathbb{R}$ then we write $L^\infty\left(\mathbb{R}_+\right)$.

* By $K$ we denote the class of increasing continuous functions $a:\mathbb{R}_+ \to \mathbb{R}_+$ with $a(0)=0$. By $K_\infty$ we denote the class of increasing continuous functions $a:\mathbb{R}_+ \to \mathbb{R}_+$ with $a(0)=0$ and $\lim_{s\to+\infty}\left(a(s)\right) = +\infty$. By $KL$ we denote the set of all continuous functions $\sigma:\mathbb{R}_+ \times \mathbb{R}_+ \to \mathbb{R}_+$ with the properties: (i) for each $t \geq 0$ the mapping $\sigma(\cdot,t)$ is of class $K$; (ii) for each $s \geq 0$, the mapping $\sigma(s,\cdot)$ is non-increasing with $\lim_{t\to+\infty}\left(\sigma(s,t)\right) = 0$.

* Let $S \subseteq \mathbb{R}^n$ be a non-empty set with $0 \in S$. We say that a function $V: S \to \mathbb{R}_+$ is positive definite if $V(x) > 0$ for all $x \in S$ with $x \neq 0$ and $V(0)=0$. We say that a continuous function $V: S \to \mathbb{R}_+$ is radially unbounded if the following property holds: "for every $M > 0$ the set $\left\{\, x \in S : V(x) \leq M \,\right\}$ is compact". For $V \in C^1(S;\mathbb{R}_+)$ we define $\nabla V(x) = \left(\frac{\partial V}{\partial x_1}(x),\ldots,\frac{\partial V}{\partial x_n}(x)\right)$.

* Let $X$ be a given normed linear space with norm $\|\cdot\|_X$, let $B \subseteq X$ be a subset of $X$ and let $I \subseteq \mathbb{R}$ be a non-empty interval. By $C^0\left(I;B\right)$ we denote the class of continuous mappings $f: I \to B$, i.e., the class of mappings $f: I \to B$ for which the following property holds: for every $t_0 \in I$ and for every $\varepsilon > 0$ there exists $\delta > 0$ such that $\|f(t) - f(t_0)\|_X < \varepsilon$ for all $t \in I \cap (t_0 - \delta, t_0 + \delta)$. By $C_c^0\left(I;X\right)$ we denote the class of continuous mappings $f: I \to X$ with compact support. By $C_c^1\left(I\right)$, where $I \subseteq \mathbb{R}$ is an open interval, we denote the class of continuously differentiable functions $f: I \to \mathbb{R}$ with compact support. By $C^1\left(I;B\right)$ we denote the class of continuous mappings $f \in C^0\left(I;B\right)$ for which the following property holds: there exists $g \in C^0\left(I;X\right)$ such for every $t_0 \in I$ and for every $\varepsilon > 0$ there exists $\delta > 0$ with $\left\|\frac{f(t)-f(t_0)}{t-t_0} - g(t_0)\right\|_X < \varepsilon$ for all $t \in I \cap \left((t_0-\delta,t_0) \cup (t_0,t_0+\delta)\right)$. In this case we write $f_t = g$.

* Let $X$ be a given Banach space with norm $\|\cdot\|_X$, $D \subseteq X$ be a subset of $X$ and let $I \subseteq \mathbb{R}$ be a given interval with non-empty interior. A function $f: I \to D$ is measurable if there exists a set $E \subset I$ of measure zero and a sequence $\left\{ f_n \in C_c^0\left(I; X\right): n \geq 1 \right\}$ such that $\lim_{n \to +\infty} \left( f_n(t) \right) = f(t)$ for all $t \in I \setminus E$. For $p \in [1, +\infty)$, $L^p(I; D)$ is the set of equivalence classes of measurable functions $f: I \to D$ with $\|f\|_p := \left( \int_I \|f(t)\|_X^p \, dt \right)^{1/p} < +\infty$. $L^\infty(I; D)$ is the set of equivalence classes of measurable functions $f: I \to D$ with $\|f\|_\infty := \operatorname{ess\,sup}_{t \in I} \left( \|f(t)\|_X \right) < +\infty$. When $D = \mathbb{R}$ then we write $L^p(I)$ and $L^\infty(I)$ instead of $L^p(I; \mathbb{R})$ and $L^\infty(I; \mathbb{R})$.

* Let $a < b$ be given constants and $k \geq 1$ be an integer. For $p \in [1, +\infty]$, $W^{k,p}(a,b)$ denotes the Sobolev space of functions in $L^p(a,b)$ with weak derivatives up to order $k$, all in $L^p(a,b)$. We set $H^k(a,b) = W^{k,2}(a,b)$. The closure of $C_c^1\left((a,b)\right)$ in $H^1(a,b)$ is denoted by $H_0^1(a,b)$.

* Let $I \subseteq \mathbb{R}$ be an interval and let $a < b$ be given constants. Let the function $u: I \times (a,b) \to \mathbb{R}$ be given. We use the notation $u[t]$ to denote the profile of $u$ at certain $t \in I$, i.e., $(u[t])(x) = u(t,x)$ for all $x \in (a,b)$. We also use the notation $u[t]$ when $u: I \to X$ and $X$ is a normed linear space and the notation $u(t)$ when $u: I \to \mathbb{R}^n$, i.e., when $X = \mathbb{R}^n$.

## 2. Results for the finite-dimensional case

In this section we provide results for the case of strict feedback systems of the form

$$\dot{\xi}_j = h_j(w, \xi_1, ..., \xi_j, \theta, d) + b_j g_j(w, \xi_1, ..., \xi_j, \theta)\xi_{j+1}$$
$$\text{for } j = 1, ..., m-1 \text{ and}$$
$$\dot{\xi}_m = h_m(w, \xi, \theta, d) + b_m g_m(w, \xi, \theta) u \tag{2.1}$$

$$\dot{w} = f(\xi_1, w, \delta, \theta) \tag{2.2}$$

where $(\xi, w) \in \mathbb{R}^m \times \mathbb{R}^l$ is the state with $\xi = (\xi_1, ..., \xi_m)$, $u \in \mathbb{R}$ is the control input, $d \in \mathbb{R}^q$, $\delta \in \mathbb{R}^\alpha$, $\theta \in \Theta \subseteq \mathbb{R}^p$, $b_j \in (0, +\infty) \cap B_j$ for $j = 1, ..., m$ are disturbances and $\Theta \subseteq \mathbb{R}^p$, $B_j \subseteq \mathbb{R}$ for $j = 1, ..., m$ are non-empty, closed sets, $f: \mathbb{R} \times \mathbb{R}^l \times \mathbb{R}^\alpha \times \Theta \to \mathbb{R}^l$, $h_j: \mathbb{R}^l \times \mathbb{R}^j \times \Theta \times \mathbb{R}^q \to \mathbb{R}$, $g_j: \mathbb{R}^l \times \mathbb{R}^j \times \Theta \to \mathbb{R}$ for $j = 1, ..., m$ are smooth mappings with $f(0,0,0,\theta) = 0$ and $h_j(0,0,\theta,0) = 0$ for $j = 1, ..., m$ and for all $\theta \in \Theta$.

The case (2.1), (2.2) is the case of strict feedback systems with no integrators. Our aim is to design a partial-state robust adaptive feedback law that depends only on $\xi \in \mathbb{R}^m$ that regulates the output $y = \xi_1$, i.e., we assume that the state component $w \in \mathbb{R}^l$ is not measured. Equivalently, $w \in \mathbb{R}^l$ is a dynamic uncertainty that affects the strict feedback system (2.1).

The case (2.1), (2.2) is to be distinguished from the case of strict feedback systems of the form

$$\dot{\zeta} = A\zeta + \bar{b}\xi_1 \tag{2.3}$$

$$\dot{\xi}_j = h_j(w,\zeta,\xi_1,...,\xi_j,\theta,d) + b_j g_j(w,\zeta,\xi_1,...,\xi_j,\theta)\xi_{j+1}$$
$$\text{for } j=1,...,m-1 \text{ and}$$
$$\dot{\xi}_m = h_m(w,\zeta,\xi,\theta,d) + b_m g_m(w,\zeta,\xi,\theta)u \tag{2.4}$$

$$\dot{w} = f(\zeta,\xi_1,w,\delta,\theta) \tag{2.5}$$

where $(\zeta,\xi,w)\in\mathbb{R}^n\times\mathbb{R}^m\times\mathbb{R}^l$ is the state with $\xi=(\xi_1,...,\xi_m)$, $u\in\mathbb{R}$ is the control input, $d\in\mathbb{R}^q$, $\delta\in\mathbb{R}^\alpha$, $\theta\in\Theta\subseteq\mathbb{R}^p$, $b_j\in(0,+\infty)\cap B_j$ for $j=1,...,m$ are disturbances and $\Theta\subseteq\mathbb{R}^p$, $B_j\subseteq\mathbb{R}$ for $j=1,...,m$ are non-empty, closed sets, $A\in\mathbb{R}^{n\times n}$ is a real matrix, $\bar{b}\in\mathbb{R}^n$, $f:\mathbb{R}^n\times\mathbb{R}\times\mathbb{R}^l\times\mathbb{R}^\alpha\times\Theta\to\mathbb{R}^l$, $h_j:\mathbb{R}^l\times\mathbb{R}^n\times\mathbb{R}^j\times\Theta\times\mathbb{R}^q\to\mathbb{R}$, $g_j:\mathbb{R}^l\times\mathbb{R}^n\times\mathbb{R}^j\times\Theta\to\mathbb{R}$ for $j=1,...,m$ are smooth mappings with $f(0,0,0,0,\theta)=0$ and $h_j(0,0,0,\theta,0)=0$ for $j=1,...,m$ and for all $\theta\in\Theta$.

The case (2.3), (2.4), (2.5) is the case of strict feedback systems with integrators. Our aim for the case (2.3), (2.4), (2.5) is to design a partial-state robust adaptive feedback law that depends only on $(\zeta,\xi)\in\mathbb{R}^n\times\mathbb{R}^m$ that regulates the output $y=(\zeta,\xi_1)\in\mathbb{R}^{n+1}$, i.e., we assume that the state component $w\in\mathbb{R}^l$ is not measured. Equivalently, $w\in\mathbb{R}^l$ is a dynamic uncertainty that affects the strict feedback system with integrators (2.3), (2.4).

2.1. Results

For the case with no integrators (2.1), (2.2) we employ the following assumption.

**Assumption (A):** *There exist smooth functions $\Phi,R:\mathbb{R}^l\to\mathbb{R}_+$ being positive definite and radially unbounded, smooth functions $\eta_j,\mu_j:\mathbb{R}^j\to(0,+\infty)$ for $j=1,...,m$, $\psi:\mathbb{R}\to(0,+\infty)$, a continuous function $\phi:\mathbb{R}\times\mathbb{R}^l\to(0,+\infty)$, a constant $\Lambda\ge 0$ and non-decreasing functions $\gamma,\bar{\gamma}\in C^0(\mathbb{R}_+;\mathbb{R}_+)$ such that the following inequalities hold:*

$$\eta_j(\xi_1,...,\xi_j)\le g_j(w,\xi_1,...,\xi_j,\theta),$$
$$\text{for all } (w,\xi_1,...,\xi_j)\in\mathbb{R}^l\times\mathbb{R}^j,\ \theta\in\Theta,\ j=1,...,m, \tag{2.6}$$

$$g_j(w,\xi_1,...,\xi_j,\theta)\le\mu_j(\xi_1,...,\xi_j)(1+|\theta|)\left(1+\sqrt{R(w)+\Lambda}\right),$$
$$\text{for all } (w,\xi_1,...,\xi_j)\in\mathbb{R}^l\times\mathbb{R}^j,\ \theta\in\Theta,\ b_j\in(0,+\infty)\cap B_j,\ j=1,...,m-1, \tag{2.7}$$

$$\left|h_j(w,\xi_1,...,\xi_j,\theta,d)\right|\le(1+|\theta|)\mu_j(\xi_1,...,\xi_j)\left(|(\xi_1,...,\xi_j)|+\sqrt{R(w)+\Lambda}\right)$$
$$+(1+|\theta|)\mu_j(\xi_1,...,\xi_j)|d|\left(1+\sqrt{R(w)+\Lambda}\right)$$
$$\text{for all } (w,\xi_1,...,\xi_j)\in\mathbb{R}^l\times\mathbb{R}^j,\ \theta\in\Theta,\ d\in\mathbb{R}^q,\ j=1,...,m, \tag{2.8}$$

$$\nabla\Phi(w)\, f(\xi_1, w, \delta, \theta) \le -R(w) + \gamma\left(|\delta|\right) + \bar{\gamma}\left(|\theta|\right) + \xi_1^2 \psi\left(\xi_1\right) + |\theta||\xi_1||\phi(\xi_1, w)|\,,$$
$$\text{for all } (w, \xi_1) \in \mathbb{R}^l \times \mathbb{R}\,,\ \theta \in \Theta\,,\ \delta \in \mathbb{R}^\alpha\,, \tag{2.9}$$

$$|\phi(\xi_1, w)| \le \mu_1(\xi_1)\left(\sqrt{R(w) + \Lambda} + |\xi_1|\right),$$
$$\text{for all } (w, \xi_1) \in \mathbb{R}^l \times \mathbb{R}\,, \tag{2.10}$$

**Remark: (a)** It follows from (2.7) that the following inequalities are valid:

$$\left|b_j g_j(w, \xi_1, ..., \xi_j, \theta)\right| \le \mu_j\left(\xi_1, ..., \xi_j\right)\left(1 + |\theta|\right)\left(1 + b_j\right)\left(1 + \sqrt{R(w) + \Lambda}\right),$$
$$\text{for all } (w, \xi_1, ..., \xi_j) \in \mathbb{R}^l \times \mathbb{R}^j\,,\ \theta \in \Theta\,,\ b_j \in (0, +\infty) \cap B_j\,,\ j = 1, ..., m-1\,, \tag{2.11}$$

Inequalities (2.11) allow us to use a backstepping lemma (see Lemma 1 below) that deals with the problem of partial-state feedback design.

**(b)** The constant $\Lambda \ge 0$ plays an important role in achieving exact regulation (possible when $\Lambda = 0$; see the statements of the main results below). However, it should be noted that $\Lambda$ has to be positive when $\limsup_{w \to 0}\left(\frac{|h_j(w, 0, ..., 0, 0, 0)|}{\sqrt{R(w)}}\right) = +\infty$ for some $j = 1, ..., m$ or $\limsup_{w \to 0}\left(\frac{|\phi(0, w)|}{\sqrt{R(w)}}\right) = +\infty$ (e.g., when $|h_j(w, 0, ..., 0, 0, 0)| = K_1 |w|$ and $R(w) = K_2 |w|^4$ for some constants $K_1, K_2 > 0$).

**(c)** The functions $\gamma, \bar{\gamma} \in C^0\left(\mathbb{R}_+; \mathbb{R}_+\right)$ appearing in (2.9) are not needed to be known for the control design. On the other hand, the functions $\psi : \mathbb{R} \to (0, +\infty)$, $R : \mathbb{R}^l \to \mathbb{R}_+$, $\phi : \mathbb{R} \times \mathbb{R}^l \to (0, +\infty)$ appearing on the right-hand side of (2.9) have to be known. The functions $\gamma, \bar{\gamma} \in C^0\left(\mathbb{R}_+; \mathbb{R}_+\right)$ appearing in (2.9) play an important role in achieving exact regulation (possible when $\gamma(0) = 0$ and $\bar{\gamma}(s) \equiv 0$; see the statements of the main results below).

**(d)** The linear growth with respect to $|\theta|$ in (2.7), (2.8) is consistent with previous works on DADS (see [9, 10, 11, 12, 13]) where the vector of the unknown parameters $\theta$ appears in an affine way in the differential equations. The growth with respect to $w$ through the term $\sqrt{R(w) + \Lambda}$ in (2.7), (2.8), (2.10) is closely related to the dissipation term $-R(w)$ in (2.9) and indicates the nature of interconnections between the subsystems (2.1) and (2.2) that can be handled/compensated by a partial-state controller that depends only on $\xi$.

The nature of inequality (2.9) is explained in detail below.

For the case with integrators (2.3), (2.4), (2.5) we employ a different assumption.

**<u>Assumption (B):</u>** *There exist smooth functions* $\Phi, R:\mathbb{R}^l \to \mathbb{R}_+$ *being positive definite and radially unbounded, smooth functions* $\eta_j, \mu_j : \mathbb{R}^n \times \mathbb{R}^j \to (0,+\infty)$ *for* $j=1,...,m$, $Q:\mathbb{R}^n \times \mathbb{R} \to \mathbb{R}_+$ *with* $Q(0,0)=0$, *a constant* $\Lambda \geq 0$ *and non-decreasing functions* $\gamma, \bar{\gamma} \in C^0(\mathbb{R}_+;\mathbb{R}_+)$ *such that the following inequalities hold:*

$$\eta_j(\zeta,\xi_1,...,\xi_j) \leq g_j(w,\zeta,\xi_1,...,\xi_j,\theta),$$
$$\text{for all } (w,\zeta,\xi_1,...,\xi_j) \in \mathbb{R}^l \times \mathbb{R}^n \times \mathbb{R}^j,\ \theta \in \Theta,\ j=1,...,m, \tag{2.12}$$

$$g_j(w,\zeta,\xi_1,...,\xi_j,\theta) \leq \mu_j(\zeta,\xi_1,...,\xi_j)(1+|\theta|)\left(1+\sqrt{R(w)+\Lambda}\right),$$
$$\text{for all } (w,\zeta,\xi_1,...,\xi_j) \in \mathbb{R}^l \times \mathbb{R}^n \times \mathbb{R}^j,\ \theta \in \Theta,\ j=1,...,m-1, \tag{2.13}$$

$$\left|h_j(w,\zeta,\xi_1,...,\xi_j,\theta,d)\right| \leq (1+|\theta|)\mu_j(\zeta,\xi_1,...,\xi_j)\left(\left|(\zeta,\xi_1,...,\xi_j)\right|+\sqrt{R(w)+\Lambda}\right)$$
$$+(1+|\theta|)\mu_j(\zeta,\xi_1,...,\xi_j)|d|\left(1+\sqrt{R(w)+\Lambda}\right)$$
$$\text{for all } (w,\zeta,\xi_1,...,\xi_j) \in \mathbb{R}^l \times \mathbb{R}^n \times \mathbb{R}^j,\ \theta \in \Theta,\ d \in \mathbb{R}^q,\ j=1,...,m, \tag{2.14}$$

$$\nabla\Phi(w) f(\zeta,\xi_1,w,\delta,\theta) \leq -R(w)+\gamma(|\delta|)+\bar{\gamma}(|\theta|)+Q(\zeta,\xi_1),$$
$$\text{for all } (w,\zeta,\xi_1) \in \mathbb{R}^l \times \mathbb{R}^n \times \mathbb{R},\ \theta \in \Theta,\ \delta \in \mathbb{R}^\alpha, \tag{2.15}$$

**Remarks: (a)** It should be noticed that inequalities (2.13) guarantee the following inequalities

$$\left|b_j g_j(w,\zeta,\xi_1,...,\xi_j,\theta)\right| \leq (1+|\theta|)(1+b_j)\mu_j(\zeta,\xi_1,...,\xi_j)\left(1+\sqrt{R(w)+\Lambda}\right),$$
$$\text{for all } (w,\zeta,\xi_1,...,\xi_j) \in \mathbb{R}^l \times \mathbb{R}^n \times \mathbb{R}^j,\ \theta \in \Theta,\ b_j \in (0,+\infty) \cap B_j,\ j=1,...,m-1, \tag{2.16}$$

Again, inequalities (2.16) allow us to use a backstepping lemma (see Lemma 1 below) that deals with the problem of partial-state feedback design.

**(b)** The difference between Assumption (A) and Assumption (B) is twofold: (i) first, in Assumption (B) all inequalities involve the state component $\zeta \in \mathbb{R}^n$ which is absent in Assumption (A), and (ii) secondly, there is a difference between inequality (2.15) and inequality (2.9). The difference is explained below.

**(c)** Again, the constant $\Lambda \geq 0$ plays an important role in achieving exact regulation (possible when $\Lambda = 0$; see the statements of the main results below). However, it should be noted that $\Lambda$ has to be positive when $\limsup_{w \to 0}\left(\frac{\left|h_j(w,0,0,...,0,0,0)\right|}{\sqrt{R(w)}}\right) = +\infty$ for some $j=1,...,m$.

**(d)** The functions $\gamma, \bar{\gamma} \in C^0(\mathbb{R}_+;\mathbb{R}_+)$ appearing in (2.15) are not needed to be known for the control design. On the other hand, the functions $Q:\mathbb{R}^n \times \mathbb{R} \to \mathbb{R}_+$, $R:\mathbb{R}^l \to \mathbb{R}_+$ appearing on the right-hand side of (2.15) have to be known. The functions $\gamma, \bar{\gamma} \in C^0(\mathbb{R}_+;\mathbb{R}_+)$ appearing in (2.15) play an important role in achieving exact regulation (possible when $\gamma(0)=0$ and $\bar{\gamma}(s) \equiv 0$; see the statements of the main results below).

**(e)** The linear growth with respect to $|\theta|$ in (2.13), (2.14) is consistent with previous works on DADS (see [9, 10, 11, 12, 13]) where the vector of the unknown parameters $\theta$ appears in an affine way in the differential equations. The growth with respect to $w$ through the term $\sqrt{R(w)+\Lambda}$ in (2.13), (2.14) is closely related to the dissipation term $-R(w)$ in (2.15) and indicates the nature of interconnections between the subsystems (2.3), (2.4) and (2.5) that can be handled/compensated by a partial-state controller that depends only on $(\zeta,\xi)$.

The nature of inequality (2.15) is explained in detail below.

Our result for the case with no integrators is stated next. It is the analogue of Theorem 6.2 on page 92 in [10].

**Theorem 1 (Partial-State DADS controller for systems in parametric strict feedback form with no integrators):** *Consider the control system (2.1), (2.2) and suppose that Assumption (A) holds. Then for every constants* $\kappa,\Gamma,\varepsilon>0$ *there exist functions* $k,V\in C^{\infty}\left(\mathbb{R}^m\times\mathbb{R}\right)$ *with* $V(0,z)=k(0,z)=0$ *for* $z\in\mathbb{R}$ *and* $V(\xi,z)>0$ *for* $\xi\neq 0$ *and a function* $\Xi\in C^0\left(\mathbb{R}_+^7\times(0,+\infty)\right)$ *such that for every* $(w_0,\xi_0,z_0)\in\mathbb{R}^l\times\mathbb{R}^m\times\mathbb{R}$, $d\in L^{\infty}\left(\mathbb{R}_+;\mathbb{R}^q\right)$, $\theta\in L^{\infty}\left(\mathbb{R}_+;\Theta\right)$, $\delta\in L^{\infty}\left(\mathbb{R}_+;\mathbb{R}^{\alpha}\right)$, $b_j\in L^{\infty}\left(\mathbb{R}_+;(0,+\infty)\cap B_j\right)$ *with* $\inf_{t\geq 0}\left(b_j(t)\right)>0$ *for* $j=1,...,m$ *the unique solution of (2.1), (2.2) with*

$$u=k(\xi,z) \tag{2.17}$$

$$\dot{z}=\Gamma\exp(-z)\left(V(\xi,z)-\varepsilon\right)^{+} \tag{2.18}$$

*and* $\left(w(0),\xi(0),z(0)\right)=\left(w_0,\xi_0,z_0\right)$ *exists and satisfies the following estimates for all* $t\geq 0$:

$$\limsup_{t\to+\infty}\left(|\xi_1(t)|\right)\leq\sqrt{2\varepsilon}\ \text{ and }\ \limsup_{t\to+\infty}\left(V(\xi(t),z(t))\right)\leq\varepsilon\,, \tag{2.19}$$

$$|w(t)|+|\xi(t)|\leq\Xi\left(|\xi_0|,|w_0|,|z_0|,\|d\|_{\infty},\|\delta\|_{\infty},\|\theta\|_{\infty},\|b\|_{\infty},\min_{i=1,...,m}\left(\inf_{s\geq 0}\left(b_i(s)\right)\right)\right), \tag{2.20}$$

$$z_0\leq z(t)\leq\lim_{s\to+\infty}\left(z(s)\right)\leq\Xi\left(|\xi_0|,|w_0|,|z_0|,\|d\|_{\infty},\|\delta\|_{\infty},\|\theta\|_{\infty},\|b\|_{\infty},\min_{i=1,...,m}\left(\inf_{s\geq 0}\left(b_i(s)\right)\right)\right). \tag{2.21}$$

*Finally, if* $\gamma(0)=\Lambda=0$, $\bar{\gamma}(s)\equiv 0$, $\lim_{t\to+\infty}\left(d(t)\right)=0$, $\lim_{t\to+\infty}\left(\delta(t)\right)=0$, $\liminf_{t\to+\infty}\left(b_i(t)\right)\geq 1/\kappa$ *for* $i=1,...,m$ *then one of the following holds: either* $\lim_{t\to+\infty}\left(|\xi(t)|\right)=\lim_{t\to+\infty}\left(|w(t)|\right)=0$ *or* $\kappa+\exp\left(\lim_{s\to+\infty}\left(z(s)\right)\right)<\limsup_{t\to+\infty}\left(|\theta(t)|+|b(t)|+|\theta(t)||b(t)|\right)$.

**Remark: (a)** The proof of Theorem 1 provides many stability estimates that are not included in the statement of Theorem 1. For example, when $\gamma(0)=\Lambda=0$, $\bar{\gamma}(s)\equiv 0$, $d(t)\equiv 0$, $\delta(t)\equiv 0$, $b_i$, $i=1,...,m$ and $\theta$ are constant with $b_i \geq \frac{1}{\kappa}$ for $i=1,...,m$ and $|\theta|+|b|+|\theta||b| \leq \kappa$, then the closed-loop system (2.1), (2.2) with (2.17), (2.18) is UGAOS with output $y=(\xi,w)$.

**(b)** Notice that Theorem 1 does not exclude the possibility of regulation of the full plant state $(\xi,w)$, when $\gamma(0)=\Lambda=0$, $\bar{\gamma}(s)\equiv 0$. Notice that Assumption (A) with $\gamma(0)=\Lambda=0$, $\bar{\gamma}(s)\equiv 0$ implies that $(\xi,w)=(0,0)$ is an equilibrium point for system (2.1), (2.2).

**(c)** The main idea of the proof of Theorem 1 is to show the two following differential inequalities for constants $c,a,\kappa>0$ with $\left(2^{m-1}+1\right)a<\kappa$ and appropriate non-negative functions $Q(\xi_1,z)$, $V(\xi,z)$:

$$\dot{V} \leq -2^{1-m}cV(\xi,z)-Q(\xi_1,z)+\frac{2^{m-1}a}{\rho}\left(\min_{i=1,...,m}(b_i)\left(\left(\frac{1}{\min_{i=1,...,m}(b_i)}-\rho\right)^+\right)^2+|d|^2+|d|^4+R(w)+\Lambda\right)$$
$$+\frac{2^{m-1}a}{\rho}\left(\left((|\theta|+|b|+|b||\theta|-\rho)^+\right)^2+\left((|\theta|+|b|+|b||\theta|-\rho)^+\right)^4+\left((|\theta|+|b|+|b||\theta|-\rho)^+\right)^8\right) \tag{2.22}$$

$$\dot{\Phi} \leq -\left(1-\frac{a}{\kappa}\right)R(w)+\gamma(|\delta|)+\bar{\gamma}(|\theta|)+Q(\xi_1,z)+\frac{a}{\rho}\left(\Lambda+\left((|\theta|-\rho)^+\right)^2+\left((|\theta|-\rho)^+\right)^4\right) \tag{2.23}$$

where $\rho=\kappa+\exp(z)$. Then by forming the function $U(w,\xi,z)=\Phi(w)+V(\xi,z)$, we can establish (2.20). By exploiting Lemma 4.1 on page 63 in [10] and the differential inequality (2.22) we can establish (2.21). Finally, the asymptotic estimates (2.19) are obtained by means of Barbălat's Lemma (see [8, 14]) and by establishing the inequality $V(\xi,z)\geq \xi_1^2/2$. Notice that the differential inequalities (2.22), (2.23) can be used for the derivation of many interesting estimates for the solutions of the closed-loop system (2.1), (2.2) with (2.17), (2.18).

Our result for the case with integrators is stated next. It is the analogue of Theorem 6.1 on page 91 in [10].

**Theorem 2 (Partial-State DADS controller for systems in parametric strict feedback form with integrators):** *Consider the control system (2.3), (2.4), (2.5) and suppose that Assumption (B) holds. Moreover, suppose that the pair of matrices* $(A,\bar{b})$ *is a stabilizable pair. Then for every constants* $\kappa,\Gamma,\varepsilon>0$ *there exist functions* $k,V\in C^\infty\left(\mathbb{R}^n\times\mathbb{R}^m\times\mathbb{R}\right)$ *with* $V(0,0,z)=k(0,0,z)=0$ *for* $z\in\mathbb{R}$ *and* $V(\zeta,\xi,z)>0$ *for* $(\zeta,\xi)\neq(0,0)$ *and a function* $\Xi\in C^0\left(\mathbb{R}_+^8\times(0,+\infty)\right)$ *such that for every* $(w_0,\zeta_0,\xi_0,z_0)\in\mathbb{R}^l\times\mathbb{R}^n\times\mathbb{R}^m\times\mathbb{R}$, $d\in L^\infty\left(\mathbb{R}_+;\mathbb{R}^q\right)$, $\theta\in L^\infty\left(\mathbb{R}_+;\Theta\right)$, $\delta\in L^\infty\left(\mathbb{R}_+;\mathbb{R}^\alpha\right)$, $b_j\in L^\infty\left(\mathbb{R}_+;(0,+\infty)\cap B_j\right)$ *with* $\inf_{t\geq 0}\left(b_j(t)\right)>0$ *for* $j=1,...,m$ *the unique solution of (2.3), (2.4), (2.5) with*

$$u=k(\zeta,\xi,z) \tag{2.24}$$

$$\dot{z} = \Gamma \exp(-z)\left(V(\zeta,\xi,z)-\varepsilon\right)^{+} \tag{2.25}$$

*and* $\left(w(0),\zeta(0),\xi(0),z(0)\right)=\left(w_0,\zeta_0,\xi_0,z_0\right)$ *exists and satisfies the following estimates for all* $t \geq 0$ *:*

$$\limsup_{t\to+\infty}\left(\left|\left(\zeta(t),\xi_1(t)\right)\right|\right) \leq \sqrt{\varepsilon} \text{ and } \limsup_{t\to+\infty}\left(V(\zeta(t),\xi(t),z(t))\right) \leq \varepsilon\,, \tag{2.26}$$

$$|w(t)|+|\zeta(t)|+|\xi(t)| \leq \Xi\left(|\zeta_0|,|\xi_0|,|w_0|,|z_0|,\|d\|_\infty,\|\delta\|_\infty,\|\theta\|_\infty,\|b\|_\infty,\min_{i=1,\dots,m}\left(\inf_{s\geq 0}\left(b_i(s)\right)\right)\right), \tag{2.27}$$

$$z_0 \leq z(t) \leq \lim_{s\to+\infty}\left(z(s)\right) \leq \Xi\left(|\zeta_0|,|\xi_0|,|w_0|,|z_0|,\|d\|_\infty,\|\delta\|_\infty,\|\theta\|_\infty,\|b\|_\infty,\min_{i=1,\dots,m}\left(\inf_{s\geq 0}\left(b_i(s)\right)\right)\right). \tag{2.28}$$

*Finally, if* $\gamma(0)=\Lambda=0$, $\bar{\gamma}(s)\equiv 0$, $\lim_{t\to+\infty}\left(d(t)\right)=0$, $\lim_{t\to+\infty}\left(\delta(t)\right)=0$, $\liminf_{t\to+\infty}\left(b_i(t)\right)\geq 1/\kappa$ *for* $i=1,\dots,m$ *then one of the following holds: either* $\lim_{t\to+\infty}\left(|\zeta(t)|\right)=\lim_{t\to+\infty}\left(|\xi(t)|\right)=\lim_{t\to+\infty}\left(|w(t)|\right)=0$ *or* $\kappa+\exp\left(\lim_{s\to+\infty}\left(z(s)\right)\right)<\limsup_{t\to+\infty}\left(|\theta(t)|+|b(t)|+|\theta(t)||b(t)|\right)$.

**Remark: (a)** Again, the proof of Theorem 2 provides many stability estimates that are not included in the statement of Theorem 2. For example, when $\gamma(0)=\Lambda=0$, $\bar{\gamma}(s)\equiv 0$, $d(t)\equiv 0$, $\delta(t)\equiv 0$, $b_i$, $i=1,\dots,m$ and $\theta$ are constant with $b_i \geq \frac{1}{\kappa}$ for $i=1,\dots,m$ and $|\theta|+|b|+|\theta||b|\leq\kappa$, then the closed-loop system (2.3), (2.4), (2.5) with (2.24), (2.25) is UGAOS with output $y=\left(\zeta,\xi,w\right)$.

**(b)** Again, it should be noticed that Theorem 2 does not exclude the possibility of regulation of the full plant state $\left(\zeta,\xi,w\right)$, when $\gamma(0)=\Lambda=0$, $\bar{\gamma}(s)\equiv 0$. Notice that Assumption (B) with $\gamma(0)=\Lambda=0$, $\bar{\gamma}(s)\equiv 0$ implies that $\left(\zeta,\xi,w\right)=\left(0,0,0\right)$ is an equilibrium point for system (2.3), (2.4), (2.5).

**(c)** Again, the main idea of the proof of Theorem 2 is to show the following differential inequality for constants $c,a,\kappa>0$ with $2^{m-1}a<\kappa$ and appropriate non-negative function $V(\zeta,\xi,z)$:

$$\begin{aligned}\dot{V} \leq & -2^{1-m}cV(\xi,z)-Q(\zeta,\xi_1)+\frac{2^{m-1}a}{\rho}\left(\min_{i=1,\dots,m}\left(b_i\right)\left(\left(\frac{1}{\min_{i=1,\dots,m}\left(b_i\right)}-\rho\right)^{+}\right)^{2}+|d|^2+|d|^4+R(w)+\Lambda\right)\\ & +\frac{2^{m-1}a}{\rho}\left(\left(\left(|\theta|+|b|+|b||\theta|-\rho\right)^{+}\right)^{2}+\left(\left(|\theta|+|b|+|b||\theta|-\rho\right)^{+}\right)^{4}+\left(\left(|\theta|+|b|+|b||\theta|-\rho\right)^{+}\right)^{8}\right)\end{aligned} \tag{2.29}$$

where $\rho=\kappa+\exp\left(z\right)$ and $Q(\zeta,\xi_1)$ is the function involved in (2.15). Then by forming the function $U\left(w,\zeta,\xi,z\right)=\Phi\left(w\right)+V\left(\zeta,\xi,z\right)$ and using the differential inequality (2.15), we can establish (2.27). By exploiting Lemma 4.1 on page 63 in [10] and the differential inequality (2.29) we can establish (2.28). Finally, the asymptotic estimates (2.26) are obtained by means of Barbălat's Lemma (see [8, 14]) and by establishing the inequality $V(\zeta,\xi,z)\geq\xi_1^2+|\zeta|^2$. Notice that the differential inequalities (2.15), (2.29) can be used for the derivation of many interesting estimates for the solutions of the closed-loop system (2.3), (2.4), (2.5) with (2.24), (2.25).

2.2. The nature of inequalities (2.9), (2.15) and the difference between them

We are now in a position to explain the nature of inequalities (2.9), (2.15) and the difference between them. Let us call the function $Q:\mathbb{R}^n\times\mathbb{R}\to\mathbb{R}_+$ appearing on the right hand-side of (2.15) the "debt function" (the "debt" that the unmeasured subsystem incurs and has to be "paid"/compensated by the controller design). Clearly, the debt function cannot depend on unknown parameters and unknown disturbances because in such a case the controller would not be able to handle the effect of this term. However, it should be noticed that there is no growth restriction on the debt function: this is achieved by utilizing a special construction of the non-negative function $V(\zeta,\xi,z)$ in the proof of Theorem 2 (see Part I of the proof of Theorem 2).

In order to explain inequalities (2.9), (2.15) let us next focus on the case $\bar{\gamma}(s)\equiv 0$ (the case where exact regulation is possible).

We first note that inequality (2.15) with $\bar{\gamma}(s)\equiv 0$ holds in the linear case where

$$f(\zeta,\xi_1,w,\delta,\theta)=\bar{A}w+\bar{f}\left(\zeta,\xi_1\right)+\bar{g}\left(\delta\right) \tag{2.30}$$

with $\bar{A}\in\mathbb{R}^{l\times l}$ being a Hurwitz matrix and $\bar{f}:\mathbb{R}^n\times\mathbb{R}\to\mathbb{R}^l$, $\bar{g}:\mathbb{R}^\alpha\to\mathbb{R}^l$ being smooth maps with $\bar{f}(0,0)=0$, $\bar{g}(0)=0$. Indeed, in this case there exists a positive definite, symmetric matrix $\bar{P}\in\mathbb{R}^{l\times l}$ such that $\bar{P}\bar{A}+\bar{A}'\bar{P}=-I_l$ and consequently by defining $\Phi(w)=w'\bar{P}w$ we get

$$\begin{aligned}\nabla\Phi(w)f(\zeta,\xi_1,w,\delta,\theta)&=-|w|^2+2w'\bar{P}\bar{f}\left(\zeta,\xi_1\right)+2w'\bar{P}\bar{g}\left(\delta\right)\\&\leq-\frac{1}{2}|w|^2+4\left|\bar{P}\right|^2\left|\bar{f}\left(\zeta,\xi_1\right)\right|^2+4\left|\bar{P}\right|^2\left|\bar{g}\left(\delta\right)\right|^2\end{aligned} \tag{2.31}$$

Therefore, inequality (2.15) with $\bar{\gamma}(s)\equiv 0$ holds $R(w)=\frac{1}{2}|w|^2$, $Q(\zeta,\xi_1)=4\left|\bar{P}\right|^2\left|\bar{f}\left(\zeta,\xi_1\right)\right|^2$ and $\gamma\in C^0\left(\mathbb{R}_+;\mathbb{R}_+\right)$ being any non-decreasing function that satisfies $\gamma\left(|\delta|\right)\geq 4\left|\bar{P}\right|^2\left|\bar{g}\left(\delta\right)\right|^2$ (here we can also achieve $\gamma(0)=0$). In this case the term $\sqrt{R(w)+\Lambda}$ appearing in (2.13), (2.14) becomes proportional to $|w|$ when $\Lambda=0$. Thus, in this case inequalities (2.13), (2.14) essentially require a linear growth with respect to $|w|$.

On the other hand, inequality (2.9) with $\bar{\gamma}(s)\equiv 0$ holds in the (semi)linear case even if $\theta$ is present in $f$ in an affine way, i.e., if

$$f(\xi_1,w,\delta,\theta)=\bar{A}w+\xi_1\tilde{f}\left(\xi_1\right)+\bar{g}\left(\delta\right)+\xi_1\tilde{A}\left(w,\xi_1\right)\theta \tag{2.32}$$

with $\bar{A} \in \mathbb{R}^{l \times l}$ being a Hurwitz matrix and $\tilde{f} : \mathbb{R} \to \mathbb{R}^l$, $\bar{g} : \mathbb{R}^\alpha \to \mathbb{R}^l$ being smooth maps with $\bar{g}(0) = 0$ and $\tilde{A}(w, \xi_1) \in \mathbb{R}^{l \times p}$ is a matrix with smooth elements. Indeed, by defining $\Phi(w) = w'\bar{P}w$ we get

$$\begin{aligned} \nabla\Phi(w) f(\xi_1, w, \delta, \theta) &= -|w|^2 + 2\xi_1 w'\bar{P}\tilde{f}(\xi_1) + 2w'\bar{P}\bar{g}(\delta) + 2\xi_1 w'\bar{P}\tilde{A}(w, \xi_1)\theta \\ &\leq -\frac{1}{2}|w|^2 + 4\xi_1^2 |\bar{P}|^2 |\tilde{f}(\xi_1)|^2 + 4|\bar{P}|^2 |\bar{g}(\delta)|^2 + 2|\theta||\xi_1||\tilde{A}'(w, \xi_1)\bar{P}w| \end{aligned} \tag{2.33}$$

Therefore, inequality (2.9) with $\bar{\gamma}(s) \equiv 0$ holds $R(w) = \frac{1}{2}|w|^2$, $\psi(\xi_1) = 4|\bar{P}|^2 |\tilde{f}(\xi_1)|^2$, $\phi(\xi_1, w) = 2|\tilde{A}'(w, \xi_1)\bar{P}w|$ and $\gamma \in C^0(\mathbb{R}_+; \mathbb{R}_+)$ being any non-decreasing function that satisfies $\gamma(|\delta|) \geq 4|\bar{P}|^2 |\bar{g}(\delta)|^2$ (here again we can also achieve $\gamma(0) = 0$).

The proof of Theorem 1 shows that inequalities (2.9), (2.10) guarantee that inequality (2.23) holds for all $(\xi, w) \in \mathbb{R}^m \times \mathbb{R}^l$, $\delta \in \mathbb{R}^\alpha$, $\theta \in \Theta$, $z \in \mathbb{R}$ with

$$Q(\xi_1, z) := \xi_1^2 \left( \psi(\xi_1) + \frac{\rho^3}{16a^3} \xi_1^2 \mu_1^4(\xi_1) + \frac{\rho^3}{2a} \mu_1^2(\xi_1) + \frac{\rho}{4a} \xi_1^2 \mu_1^2(\xi_1) + \rho\mu_1(\xi_1) \right) \tag{2.34}$$

where $\rho := \kappa + \exp(z)$ and $a = r\kappa$. Notice that the debt function in the case of Theorem 1 is dependent on the dynamic gain $z$; this is not allowed in Theorem 2 where the debt function depends on $(\zeta, \xi_1)$. The reason that explains this difference has to do with the backstepping procedure: while in Theorem 1 the state component $\xi_1$ is to be led to zero, in Theorem 2 the state component $\xi_1$ tracks the manifold $\xi_1 = \tilde{k}'\zeta$, where $\xi_1 = \tilde{k}'\zeta$ is a linear stabilizing feedback law for (2.3). While in Theorem 1 we can lead $\xi_1$ to zero sufficiently fast (in order to overcome the effect of the dynamic gain $z$), this is not possible in Theorem 2 (where $(\zeta, \xi_1)$ is led to zero with a specific exponential convergence rate determined by the feedback law $\xi_1 = \tilde{k}'\zeta$).

## 2.3. Formulae for the feedback law

Although Theorem 1 is an existential result, its proof is constructive by means of a backstepping lemma (Lemma 1 below). Therefore, we can provide an algorithm for the DADS feedback law (2.17), (2.18) whose existence is guaranteed by Theorem 1 under assumption (A).

Define for $i = 2, ..., m$:

$$M_i(\xi_1, ..., \xi_i) := \left(1 + |(\xi_1, ..., \xi_i)|^2\right) \sum_{j=1}^{i} \mu_j(\xi_1, ..., \xi_j) \tag{2.35}$$

**Algorithm:** Let arbitrary $\kappa, \Gamma, \varepsilon, c > 0$, $\beta \in (0,1)$ be given. Define $a = \dfrac{\beta\kappa}{1+2^{m-1}}$.

Step 1: Set

$$k_1(\xi_1, z) := -\rho P_1(\xi_1, z)\frac{4a + P_1(\xi_1, z)\xi_1^2}{4a\eta_1(\xi_1)}\xi_1$$
$$V_1(\xi_1, z) := \frac{1}{2}\xi_1^2 \tag{2.36}$$
$$\sigma_1(\xi_1, z) := \frac{1}{2}$$

where $\rho := \kappa + \exp(z)$ and $P_1(\xi_1, z)$ is any smooth function that satisfies the following inequality for all $(\xi_1, z) \in \mathbb{R} \times \mathbb{R}$:

$$\begin{aligned} P_1(\xi_1, z) \ge \frac{c}{2} + \psi(\xi_1) + \frac{\rho}{2a}\left(3(1+\rho)^2 + \rho^2\right)\mu_1^2(\xi_1) + \mu_1(\xi_1)(1+2\rho) \\ + \frac{\rho}{2a}\left(\left(\frac{5}{4} + (1+\rho)^4\right)\frac{\rho^2}{a^2}\mu_1^2(\xi_1) + 1\right)\mu_1^2(\xi_1)\xi_1^2 + \frac{\rho^7}{16a^7}\mu_1^8(\xi_1)\xi_1^6 \end{aligned} \tag{2.37}$$

Step $j=2,\ldots,m$: Set

$$V_j(\xi_1,\ldots,\xi_j,z) := V_{j-1}(\xi_1,\ldots,\xi_{j-1},z) + \frac{1}{2}\left(\xi_j - k_{j-1}(\xi_1,\ldots,\xi_{j-1},z)\right)^2 \tag{2.38}$$

Find a smooth function $S_{j-1}(\xi_1,\ldots,\xi_{j-1},z)$ that satisfies

$$\left|\frac{\partial V_{j-1}}{\partial \xi}(\xi_1,\ldots,\xi_{j-1},z)\right| + \left|k_{j-1}(\xi_1,\ldots,\xi_{j-1},z)\right| \le S_{j-1}(\xi_1,\ldots,\xi_{j-1},z)\left|(\xi_1,\ldots,\xi_{j-1})\right| \tag{2.39}$$

and set

$$\sigma_j(\xi_1,\ldots,\xi_j,z) := \left(1 + 2S_{j-1}^2(\xi_1,\ldots,\xi_{j-1},z)\right)\sigma_{j-1}(\xi_1,\ldots,\xi_{j-1},z) + 4 \tag{2.40}$$

Let $P_j(\xi_1,\ldots,\xi_j,z)$ be any smooth function that satisfies

$$\begin{aligned} P_j(\xi_1,\ldots,\xi_j,z) &\ge \frac{9c}{20} + G_{4,j}(\xi_1,\ldots,\xi_j,z) + G_{1,j}(\xi_1,\ldots,\xi_j,z)M_j^2(\xi_1,\ldots,\xi_j) \\ &+ G_{2,j}(\xi_1,\ldots,\xi_j,z)M_j^4(\xi_1,\ldots,\xi_j)\left(\xi_j - k_{j-1}(\xi_1,\ldots,\xi_{j-1},z)\right)^2 \\ &+ G_{3,j}(\xi_1,\ldots,\xi_j,z)M_j^8(\xi_1,\ldots,\xi_j)\left(\xi_j - k_{j-1}(\xi_1,\ldots,\xi_{j-1},z)\right)^6 \end{aligned} \tag{2.41}$$

where

$$
\begin{aligned}
G_{1,j}\left(\xi_1,...,\xi_j,z\right) &:= \frac{\rho}{a}\left(4\left(1+\rho\right)^2+\left|\left(\xi_1,...,\xi_j\right)\right|^2\right)\left(1+\left|\frac{\partial k_{j-1}}{\partial \xi}\left(\xi_1,...,\xi_{j-1},z\right)\right|^2\right) \\
&+\frac{\rho}{a}\left(1+3\left(1+\rho\right)^2\right)\left(\left|\frac{\partial V_{j-1}}{\partial \xi}\left(\xi_1,...,\xi_{j-1},z\right)\right|^2+\left|\frac{\partial k_{j-1}}{\partial \xi}\left(\xi_1,...,\xi_{j-1},z\right)\right|^2\left|\xi_j\right|^2\right) \\
&+\frac{5\left(1+\rho\right)^2\sigma_j\left(\xi_1,...,\xi_j,z\right)}{8c}\left(S_{j-1}^2\left(\xi_1,...,\xi_{j-1},z\right)+2\right)\left(3+\left|\frac{\partial k_{j-1}}{\partial \xi}\left(\xi_1,...,\xi_{j-1},z\right)\right|^2\right) \\
&+\frac{5\left(1+\rho\right)^2\sigma_j\left(\xi_1,...,\xi_j,z\right)}{8c}\left|\frac{\partial k_{j-1}}{\partial \xi}\left(\xi_1,...,\xi_{j-1},z\right)\right|^2\left(3+\left|\frac{\partial k_{j-1}}{\partial \xi}\left(\xi_1,...,\xi_{j-1},z\right)\right|^2\right)
\end{aligned} \tag{2.42}
$$

$$
\begin{aligned}
G_{2,j}\left(\xi_1,...,\xi_j,z\right) &:= \frac{\rho^3}{a^3}\left(9\left(1+\rho\right)^4+16\right)\left(1+\left|\frac{\partial k_{j-1}}{\partial \xi}\left(\xi_1,...,\xi_{j-1},z\right)\right|^4\right) \\
&+\frac{9\rho^3}{a^3}\left(\left|\frac{\partial V_{j-1}}{\partial \xi}\left(\xi_1,...,\xi_{j-1},z\right)\right|^4+\left|\frac{\partial k_{j-1}}{\partial \xi}\left(\xi_1,...,\xi_{j-1},z\right)\right|^4\left|\xi_j\right|^4\right)
\end{aligned} \tag{2.43}
$$

$$
G_{3,j}\left(\xi_1,...,\xi_j,z\right) := \frac{81\rho^7}{2a^7}\left(1+\left|\frac{\partial k_{j-1}}{\partial \xi}\left(\xi_1,...,\xi_{j-1},z\right)\right|^8\right) \tag{2.44}
$$

$$
\begin{aligned}
G_{4,j}\left(\xi_1,...,\xi_j,z\right) &:= \frac{\Gamma\exp(-z)}{4}\left(1+\left|\frac{\partial V_{j-1}}{\partial z}\left(\xi_1,...,\xi_{j-1},z\right)\right|^2\right) \\
&+\frac{\Gamma\exp(-z)}{4}\left(1+\left|\frac{\partial k_{j-1}}{\partial z}\left(\xi_1,...,\xi_{j-1},z\right)\right|^2\right)\left(1+\left(\xi_j-k_{j-1}\left(\xi_1,...,\xi_{j-1},z\right)\right)^2\right) \\
&+\frac{5\Gamma^2\exp(-2z)}{2c}\left|\frac{\partial k_{j-1}}{\partial z}\left(\xi_1,...,\xi_{j-1},z\right)\right|^2 V_{j-1}\left(\xi_1,...,\xi_{j-1},z\right)
\end{aligned} \tag{2.45}
$$

Set

$$
\begin{aligned}
&k_j\left(\xi_1,...,\xi_j,z\right) \\
&:= -\rho P_j\left(\xi_1,...,\xi_j,z\right)\frac{4a+P_j\left(\xi_1,...,\xi_j,z\right)\left(\xi_j-k_{j-1}\left(\xi_1,...,\xi_{j-1},z\right)\right)^2}{4a\eta_j\left(\xi_1,...,\xi_j\right)}\left(\xi_j-k_{j-1}\left(\xi_1,...,\xi_{j-1},z\right)\right)
\end{aligned} \tag{2.46}
$$

The above algorithm can be applied when assumption (A) holds but gives an "aggressive" feedback law which can become smoother when additional assumptions are valid. For example, when we know that $b_j(t)\geq 1$ for all $t\geq 0$ and for some $j=2,...,m$ (i.e., if $B_j\subseteq[1,+\infty)$) then (2.46) can be replaced by the simpler formula

$$
k_j\left(\xi_1,...,\xi_j,z\right) := -\frac{P_j\left(\xi_1,...,\xi_j,z\right)}{\eta_j\left(\xi_1,...,\xi_j\right)}\left(\xi_j-k_{j-1}\left(\xi_1,...,\xi_{j-1},z\right)\right) \tag{2.47}
$$

Similarly, if $b_1(t) \geq 1$ for all $t \geq 0$ then we can use

$$k_1(\xi_1, z) := -\frac{P_1(\xi_1, z)}{\eta_1(\xi_1)}\xi_1 \tag{2.48}$$

instead of (2.36).

Unfortunately, we cannot give a similar algorithm for the DADS feedback law (2.24), (2.25) whose existence is guaranteed by Theorem 2 under assumption (B) because additional technical steps are used in the proof of Theorem 2.

## 3. Results for the infinite-dimensional case

We provide results for the case of strict feedback systems of the form

$$\begin{gathered} \dot{\xi}_j = h_j(w, \xi_1, ..., \xi_j, \theta, d) + b_j g_j(w, \xi_1, ..., \xi_j, \theta)\xi_{j+1} \\ \text{for } j = 1, ..., m-1 \text{ and} \\ \dot{\xi}_m = h_m(w, \xi, \theta, d) + b_m g_m(w, \xi, \theta)u \end{gathered} \tag{3.1}$$

where $(\xi, w) \in \mathbb{R}^m \times X$ is the state with $\xi = (\xi_1, ..., \xi_m)$, $X$ is a Banach space with norm $\|w\|$, $u \in \mathbb{R}$ is the control input, $d \in \mathbb{R}^q$, $\theta \in \Theta \subseteq \mho$ where $\mho$ is a Banach space with norm $\|\theta\|_{\mho}$, $b_j \in (0, +\infty) \cap B_j$ for $j = 1, ..., m$ are disturbances, $B_j \subseteq \mathbb{R}$ for $j = 1, ..., m$ are non-empty, closed sets, $h_j : X \times \mathbb{R}^j \times \Theta \times \mathbb{R}^q \to \mathbb{R}$, $g_j : X \times \mathbb{R}^j \times \Theta \to \mathbb{R}$ for $j = 1, ..., m$ are continuous mappings with $h_j(0, 0, \theta, 0) = 0$ for $j = 1, ..., m$ and for all $\theta \in \Theta$.

Again, the case (3.1) is the case of strict feedback systems with no integrators. Our aim is to design a partial-state robust adaptive feedback law that depends only on $\xi \in \mathbb{R}^m$ that regulates the output $y = \xi_1$, i.e., we assume that the state component $w \in X$ is not measured. Equivalently, $w \in X$ is a dynamic uncertainty that affects the strict feedback system (2.1). However here, we do not provide an explicit description of the dynamics of the infinite-dimensional state component $w \in X$ .

The case (3.1) is to be distinguished from the case of strict feedback systems of the form

$$\dot{\zeta} = A\zeta + \bar{b}\xi_1 \tag{3.2}$$

$$\begin{gathered} \dot{\xi}_j = h_j(w, \zeta, \xi_1, ..., \xi_j, \theta, d) + b_j g_j(w, \zeta, \xi_1, ..., \xi_j, \theta)\xi_{j+1} \\ \text{for } j = 1, ..., m-1 \text{ and} \\ \dot{\xi}_m = h_m(w, \zeta, \xi, \theta, d) + b_m g_m(w, \zeta, \xi, \theta)u \end{gathered} \tag{3.3}$$

where $(\zeta, \xi, w) \in \mathbb{R}^n \times \mathbb{R}^m \times X$ is the state with $\xi = (\xi_1, ..., \xi_m)$, $u \in \mathbb{R}$ is the control input, $X$ is a Banach space with norm $\|w\|$, $d \in \mathbb{R}^q$, $\theta \in \Theta \subseteq \mho$ where $\mho$ is a Banach space with norm $\|\theta\|_{\mho}$, $b_j \in (0, +\infty) \cap B_j$ for $j = 1, ..., m$ are disturbances, $B_j \subseteq \mathbb{R}$ for $j = 1, ..., m$ are non-empty, closed

sets, $A\in\mathbb{R}^{n\times n}$ is a real matrix, $\bar{b}\in\mathbb{R}^{n}$, $h_j : X\times\mathbb{R}^n\times\mathbb{R}^j\times\Theta\times\mathbb{R}^q\to\mathbb{R}$, $g_j : X\times\mathbb{R}^n\times\mathbb{R}^j\times\Theta\to\mathbb{R}$ for $j=1,...,m$ are continuous mappings with $h_j(0,0,0,\theta,0)=0$ for $j=1,...,m$ and for all $\theta\in\Theta$.

Again, the case (3.2), (3.3) is the case of strict feedback systems with integrators. Our aim for the case (3.2), (3.3) is to design a partial-state robust adaptive feedback law that depends only on $(\zeta,\xi)\in\mathbb{R}^n\times\mathbb{R}^m$ that regulates the output $y=(\zeta,\xi_1)\in\mathbb{R}^{n+1}$, i.e., we assume that the state component $w\in X$ is not measured. Equivalently, $w\in X$ is a dynamic uncertainty that affects the strict feedback system with integrators (3.2), (3.3).

For the case with no integrators (3.1) we employ the following assumption.

**Assumption (A'):** *There exist a Banach space $S\subseteq X$, continuous functionals $\Phi,R : X\to\mathbb{R}_+$, $\phi:\mathbb{R}\times S\to(0,+\infty)$, functions $a_1,a_2,a_3,a_4\in K_\infty$, smooth functions $\eta_j,\mu_j:\mathbb{R}^j\to(0,+\infty)$ for $j=1,...,m$, $\psi:\mathbb{R}\to(0,+\infty)$, a constant $\Lambda\geq 0$ and non-decreasing functions $\gamma,\bar{\gamma}\in C^0(\mathbb{R}_+;\mathbb{R}_+)$ such that inequality (2.10) and the following inequalities hold:*

$$\eta_j(\xi_1,...,\xi_j)\leq g_j(w,\xi_1,...,\xi_j,\theta),$$
$$\text{for all } (w,\xi_1,...,\xi_j)\in S\times\mathbb{R}^j,\ \theta\in\Theta,\ j=1,...,m, \tag{3.4}$$

$$g_j(w,\xi_1,...,\xi_j,\theta)\leq\mu_j(\xi_1,...,\xi_j)\left(1+\|\theta\|_{\mho}\right)\left(1+\sqrt{R(w)+\Lambda}\right),$$
$$\text{for all } (w,\xi_1,...,\xi_j)\in S\times\mathbb{R}^j,\ \theta\in\Theta,\ b_j\in(0,+\infty)\cap B_j,\ j=1,...,m-1, \tag{3.5}$$

$$\left|h_j(w,\xi_1,...,\xi_j,\theta,d)\right|\leq\left(1+\|\theta\|_{\mho}\right)\mu_j(\xi_1,...,\xi_j)\left(\left|(\xi_1,...,\xi_j)\right|+\sqrt{R(w)+\Lambda}\right)$$
$$+\left(1+\|\theta\|_{\mho}\right)\mu_j(\xi_1,...,\xi_j)|d|\left(1+\sqrt{R(w)+\Lambda}\right)$$
$$\text{for all } (w,\xi_1,...,\xi_j)\in S\times\mathbb{R}^j,\ \theta\in\Theta,\ d\in\mathbb{R}^q,\ j=1,...,m, \tag{3.6}$$

$$\dot{\Phi}(w,\xi_1,\theta,\delta)\leq -R(w)+\gamma\left(\|\delta\|_{\mathcal{E}}\right)+\bar{\gamma}\left(\|\theta\|_{\mho}\right)+\xi_1^2\psi(\xi_1)+\|\theta\|_{\mho}|\xi_1||\phi(\xi_1,w)|,$$
$$\text{for all } (w,\xi_1)\in S\times\mathbb{R},\ \theta\in\Theta,\ \delta\in\mathcal{E}, \tag{3.7}$$

$$a_1(\Phi(w))\leq R(w)\leq a_2(\Phi(w)),\ \text{for all } w\in X \tag{3.8}$$

$$a_3(\|w\|)\leq\Phi(w)\leq a_4(\|w\|),\ \text{for all } w\in X \tag{3.9}$$

*where $\mathcal{E}$ is a Banach space with norm $\|\delta\|_{\mathcal{E}}$ and $\dot{\Phi}:S\times\mathbb{R}\times\Theta\times\mathcal{E}\to\mathbb{R}$ is a functional that is equal to the upper right Dini derivative of $\Phi:X\to\mathbb{R}_+$ for every solution for which $w[t]\in S$, i.e., $\dot{\Phi}(w[t],\xi_1(t),\theta[t],\delta[t])=\limsup_{h\to 0^+}\left(h^{-1}\left(\Phi(w[t+h])-\Phi(w[t])\right)\right)$ for all $t\geq 0$ for which $w[t]\in S$ and $\limsup_{h\to 0^+}\left(h^{-1}\left(\Phi(w[t+h])-\Phi(w[t])\right)\right)$ exists.*

For the case with integrators (3.2), (3.3) we need to employ a different assumption.

**Assumption (B'):** *There exist a Banach space* $S \subseteq X$, *continuous functionals* $\Phi, R : X \to \mathbb{R}_+$, *functions* $a_1, a_2, a_3, a_4 \in K_\infty$, *smooth functions* $\eta_j, \mu_j : \mathbb{R}^n \times \mathbb{R}^j \to (0,+\infty)$ *for* $j = 1,...,m$, $Q : \mathbb{R}^n \times \mathbb{R} \to \mathbb{R}_+$ *with* $Q(0,0) = 0$, *a constant* $\Lambda \ge 0$ *and non-decreasing functions* $\gamma, \bar{\gamma} \in C^0(\mathbb{R}_+;\mathbb{R}_+)$ *such that inequalities (3.8), (3.9) and the following inequalities hold:*

$$\eta_j(\zeta,\xi_1,...,\xi_j) \le g_j(w,\zeta,\xi_1,...,\xi_j,\theta), \quad \text{for all } (w,\zeta,\xi_1,...,\xi_j) \in S \times \mathbb{R}^n \times \mathbb{R}^j,\ \theta \in \Theta,\ j = 1,...,m, \tag{3.10}$$

$$g_j(w,\zeta,\xi_1,...,\xi_j,\theta) \le \mu_j(\zeta,\xi_1,...,\xi_j)\left(1+\|\theta\|_\mho\right)\left(1+\sqrt{R(w)+\Lambda}\right), \quad \text{for all } (w,\zeta,\xi_1,...,\xi_j) \in S \times \mathbb{R}^n \times \mathbb{R}^j,\ \theta \in \Theta,\ j = 1,...,m-1, \tag{3.11}$$

$$\begin{aligned} \left|h_j(w,\zeta,\xi_1,...,\xi_j,\theta,d)\right| \le& \left(1+\|\theta\|_\mho\right)\mu_j(\zeta,\xi_1,...,\xi_j)\left(\left|(\zeta,\xi_1,...,\xi_j)\right| + \sqrt{R(w)+\Lambda}\right) \\ &+ \left(1+\|\theta\|_\mho\right)\mu_j(\zeta,\xi_1,...,\xi_j)|d|\left(1+\sqrt{R(w)+\Lambda}\right) \\ &\text{for all } (w,\zeta,\xi_1,...,\xi_j) \in S \times \mathbb{R}^n \times \mathbb{R}^j,\ \theta \in \Theta,\ d \in \mathbb{R}^q,\ j = 1,...,m, \end{aligned} \tag{3.12}$$

$$\dot{\Phi}(w,\zeta,\xi_1,\theta,\delta) \le -R(w) + \gamma\left(\|\delta\|_{\mathcal{E}}\right) + \bar{\gamma}\left(\|\theta\|_\mho\right) + Q(\zeta,\xi_1), \quad \text{for all } (w,\zeta,\xi_1) \in S \times \mathbb{R}^n \times \mathbb{R},\ \theta \in \Theta,\ \delta \in \mathcal{E}, \tag{3.13}$$

*where* $\mathcal{E}$ *is a Banach space with norm* $\|\delta\|_{\mathcal{E}}$ *and* $\dot{\Phi} : S \times \mathbb{R}^n \times \mathbb{R} \times \Theta \times \mathcal{E} \to \mathbb{R}$ *is a functional that is equal to the upper right Dini derivative of* $\Phi : X \to \mathbb{R}_+$ *for every solution for which* $w[t] \in S$, *i.e.,* $\dot{\Phi}\left(w[t],\zeta(t),\xi_1(t),\theta[t],\delta[t]\right) = \limsup_{h\to 0^+}\left(h^{-1}\left(\Phi\left(w[t+h]\right) - \Phi\left(w[t]\right)\right)\right)$ *for all* $t \ge 0$ *for which* $w[t] \in S$ *and* $\limsup_{h\to 0^+}\left(h^{-1}\left(\Phi\left(w[t+h]\right) - \Phi\left(w[t]\right)\right)\right)$ *exists.*

**Remark:** A difference of the finite-dimensional case with the infinite-dimensional case is that Assumption (A') and Assumption (B') require the existence of functions $a_1, a_2, a_3, a_4 \in K_\infty$ so that inequalities (3.8), (3.9) hold. This requirement does not appear in Assumption (A) and Assumption (B) for the finite-dimensional case. Instead, in the finite-dimensional case it is assumed that the functions $\Phi, R : \mathbb{R}^l \to \mathbb{R}_+$ are positive definite and radially unbounded.

Our result for the case with no integrators is stated next. Its proof is exactly the same with the proof of Theorem 1 and is omitted.

**Theorem 3 (Partial-State DADS controller for infinite-dimensional systems in parametric strict feedback form with no integrators-The analogue of Theorem 1 in the infinite-dimensional case):** *Consider the control system (3.1) and suppose that Assumption (A') holds. Then for every constants* $\kappa, \Gamma, \varepsilon > 0$ *there exist functions* $k, V \in C^\infty(\mathbb{R}^m \times \mathbb{R})$ *with* $V(0,z) = k(0,z) = 0$ *for* $z \in \mathbb{R}$ *and* $V(\xi,z) > 0$ *for* $\xi \ne 0$ *and a function* $\Xi \in C^0\left(\mathbb{R}_+^7 \times (0,+\infty)\right)$ *such that every solution*

$$w \in C^0\left(\mathbb{R}_+;X\right) \text{ with } \mathbb{R}_+ \ni t \to \Phi\left(w[t]\right) \text{ being locally absolutely continuous}$$
$$\text{and } w \in L^1\left((0,T);S\right) \text{ for each } T>0$$
$$z \in C^1\left(\mathbb{R}_+\right), \xi \in C^0\left(\mathbb{R}_+;\mathbb{R}^m\right) \text{ being locally absolutely continuous}$$

*of the closed-loop system (3.1) with (2.17), (2.18) that corresponds to some* $d \in L^\infty\left(\mathbb{R}_+;\mathbb{R}^q\right)$, $\theta \in L^\infty\left(\mathbb{R}_+;\Theta\right)$, $\delta \in L^\infty\left(\mathbb{R}_+;\mathcal{E}\right)$, $b_j \in L^\infty\left(\mathbb{R}_+;(0,+\infty)\cap B_j\right)$ *with* $\inf_{t\geq 0}\left(b_j(t)\right)>0$ *for* $j=1,...,m$ *satisfies the following estimates for all* $t\geq 0$*:*

$$\limsup_{t\to+\infty}\left(\left|\xi_1(t)\right|\right) \leq \sqrt{2\varepsilon} \text{ and } \limsup_{t\to+\infty}\left(V(\xi(t),z(t))\right) \leq \varepsilon\,, \tag{3.14}$$

$$\left\|w[t]\right\| + \left|\xi(t)\right| \leq \Xi\left(\left|\xi(0)\right|,\left\|w[0]\right\|,\left|z(0)\right|,\left\|d\right\|_\infty,\left\|\delta\right\|_\infty,\left\|\theta\right\|_\infty,\left\|b\right\|_\infty,\min_{i=1,...,m}\left(\inf_{s\geq 0}\left(b_i(s)\right)\right)\right), \tag{3.15}$$

$$z(0) \leq z(t) \leq \lim_{s\to+\infty}\left(z(s)\right) \leq \Xi\left(\left|\xi(0)\right|,\left\|w[0]\right\|,\left|z(0)\right|,\left\|d\right\|_\infty,\left\|\delta\right\|_\infty,\left\|\theta\right\|_\infty,\left\|b\right\|_\infty,\min_{i=1,...,m}\left(\inf_{s\geq 0}\left(b_i(s)\right)\right)\right). \tag{3.16}$$

*Finally, if* $\gamma(0)=\Lambda=0$, $\bar{\gamma}(s)\equiv 0$, $\lim_{t\to+\infty}\left(d(t)\right)=0$, $\lim_{t\to+\infty}\left(\delta(t)\right)=0$, $\liminf_{t\to+\infty}\left(b_i(t)\right)\geq \frac{1}{\kappa}$ *for* $i=1,...,m$ *then one of the following holds: either* $\lim_{t\to+\infty}\left(\left|\xi(t)\right|\right)=\lim_{t\to+\infty}\left(\left\|w[t]\right\|\right)=0$ *or* $\kappa+\exp\left(\lim_{s\to+\infty}\left(z(s)\right)\right) < \limsup_{t\to+\infty}\left(\left\|\theta[t]\right\|_\mho + \left|b(t)\right| + \left\|\theta[t]\right\|_\mho\left|b(t)\right|\right)$.

Our result for the case with integrators is stated next. Its proof is exactly the same with the proof of Theorem 2 and is omitted.

**Theorem 4 (Partial-State DADS controller for infinite-dimensional systems in parametric strict feedback form with integrators-The analogue of Theorem 2 in the infinite-dimensional case):** *Consider the control system (3.2), (3.3) and suppose that Assumption (B') holds. Moreover, suppose that the pair of matrices* $\left(A,\bar{b}\right)$ *is a stabilizable pair. Then for every constants* $\kappa,\Gamma,\varepsilon>0$ *there exist functions* $k,V \in C^\infty\left(\mathbb{R}^n\times\mathbb{R}^m\times\mathbb{R}\right)$ *with* $V(0,0,z)=k(0,0,z)=0$ *for* $z\in\mathbb{R}$ *and* $V(\zeta,\xi,z)>0$ *for* $(\zeta,\xi)\neq(0,0)$ *and a function* $\Xi \in C^0\left(\mathbb{R}_+^8\times(0,+\infty)\right)$ *such that every solution*

$$w \in C^0\left(\mathbb{R}_+;X\right) \text{ with } \mathbb{R}_+ \ni t \to \Phi\left(w[t]\right) \text{ being locally absolutely continuous}$$
$$\text{and } w \in L^1\left((0,T);S\right) \text{ for each } T>0$$
$$z \in C^1\left(\mathbb{R}_+\right), (\zeta,\xi) \in C^0\left(\mathbb{R}_+;\mathbb{R}^n\times\mathbb{R}^m\right) \text{ being locally absolutely continuous}$$

*of the closed-loop system (3.2), (3.3) with (2.24), (2.25) that corresponds to some* $d \in L^\infty\left(\mathbb{R}_+;\mathbb{R}^q\right)$, $\theta \in L^\infty\left(\mathbb{R}_+;\Theta\right)$, $\delta \in L^\infty\left(\mathbb{R}_+;\mathcal{E}\right)$, $b_j \in L^\infty\left(\mathbb{R}_+;(0,+\infty)\cap B_j\right)$ *with* $\inf_{t\geq 0}\left(b_j(t)\right)>0$ *for* $j=1,...,m$ *satisfies the following estimates for all* $t\geq 0$*:*

$$\limsup_{t\to+\infty}\left(\left|\left(\zeta(t),\xi_1(t)\right)\right|\right) \leq \sqrt{\varepsilon} \text{ and } \limsup_{t\to+\infty}\left(V(\zeta(t),\xi(t),z(t))\right) \leq \varepsilon\,, \tag{3.17}$$

$$\|w[t]\|+|\zeta(t)|+|\xi(t)|\le \Xi\left(|\zeta(0)|,|\xi(0)|,\|w[0]\|,|z(0)|,\|d\|_\infty,\|\delta\|_\infty,\|\theta\|_\infty,\|b\|_\infty,\min_{i=1,...,m}\left(\inf_{s\ge 0}\left(b_i(s)\right)\right)\right),\quad (3.18)$$

$$z(0)\le z(t)\le \lim_{s\to+\infty}\left(z(s)\right)\le \Xi\left(|\zeta(0)|,|\xi(0)|,\|w[0]\|,|z(0)|,\|d\|_\infty,\|\delta\|_\infty,\|\theta\|_\infty,\|b\|_\infty,\min_{i=1,...,m}\left(\inf_{s\ge 0}\left(b_i(s)\right)\right)\right).\quad (3.19)$$

*Finally, if* $\gamma(0)=\Lambda=0$, $\bar{\gamma}(s)\equiv 0$, $\lim_{t\to+\infty}\left(d(t)\right)=0$, $\lim_{t\to+\infty}\left(\delta(t)\right)=0$, $\liminf_{t\to+\infty}\left(b_i(t)\right)\ge 1/\kappa$ *for* $i=1,...,m$ *then one of the following holds: either* $\lim_{t\to+\infty}\left(|\zeta(t)|\right)=\lim_{t\to+\infty}\left(|\xi(t)|\right)=\lim_{t\to+\infty}\left(|w(t)|\right)=0$ *or* $\kappa+\exp\left(\lim_{s\to+\infty}\left(z(s)\right)\right)<\limsup_{t\to+\infty}\left(|\theta(t)|+|b(t)|+|\theta(t)||b(t)|\right)$.

## 4. Examples

**Example 1:** Consider the following system

$$\begin{aligned}\dot{w}&=-w+\xi_1^2+\delta\\ \dot{\xi}_1&=2w+\theta\xi_1^2+\xi_2+d\\ \dot{\xi}_2&=u\end{aligned}\quad (4.1)$$

where $w\in\mathbb{R}$ is an unmeasured state, $\xi=\left(\xi_1,\xi_2\right)\in\mathbb{R}^2$ is the measured state, $u\in\mathbb{R}$ is the control input, $\left(\delta,d\right)\in\mathbb{R}^2$ is the time-varying disturbance and $\theta\in\mathbb{R}$ is an unknown, possibly time-varying parameter. System (4.1) with constant $\theta$ was studied in the paper [6], where the authors designed a robust partial-state adaptive controller that can regulate practically the output

$$y=\xi_1\quad (4.2)$$

and avoid the phenomenon of parameter drift. The controller in [6] is described by the equations

$$\begin{aligned}u&=-2\zeta-\xi_1+A\left(\xi_2+\hat{\theta}\xi_1^2\right)-\xi_1^2\Gamma\left(\xi_1^3-A\xi_1^2\zeta-\sigma\hat{\theta}\right)\\ &-p\left(\tanh\left(2\delta r\xi_1\right)+\frac{2\delta r\xi_1}{\cosh^2\left(2\delta r\xi_1\right)}\right)\left(-r+2.5\xi_1^4+0.625\right)\\ &-B\left(-\lambda\sigma\left(p-1\right)+\lambda B\xi_1+\lambda\zeta C\right)-pC\end{aligned}\quad (4.3)$$

$$\begin{aligned}\dot{r}&=-r+2.5\xi_1^4+0.625\\ \frac{d\hat{\theta}}{dt}&=\Gamma\left(\xi_1^3-A\xi_1^2\zeta-\sigma\hat{\theta}\right)\\ \dot{p}&=-\lambda\sigma\left(p-1\right)+\lambda B\xi_1+\lambda\zeta C\end{aligned}\quad (4.4)$$

where $r, p \geq 0$ are non-negative controller states, $\hat{\theta} \in \mathbb{R}$ is the parameter estimate, $\lambda, \sigma, \Gamma > 0$ are parameters of the controller and

$$\begin{aligned} &\zeta = \xi_2 + 2\xi_1 + \hat{\theta}\xi_1^2 + 0.25\xi_1^3 + Bp \\ &A = -2 - 2\hat{\theta}\xi_1 - 0.75\xi_1^2 - 2\bar{\delta}\, p\left( \frac{1}{\cosh^2\left(2\bar{\delta}\xi_1\right)} + \frac{4}{\cosh^2\left(4\bar{\delta}\xi_1\right)} + \frac{r^2}{\cosh^2\left(2\bar{\delta} r \xi_1\right)} \right) \\ &B = \tanh\left(2\bar{\delta}\xi_1\right) + 2\tanh\left(4\bar{\delta}\xi_1\right) + r\tanh\left(2\bar{\delta} r\xi_1\right) \\ &C = A\tanh\left(2\bar{\delta} A\zeta\right) + r\tanh\left(2\bar{\delta} r\zeta\right) + 2\tanh\left(4\bar{\delta}\zeta\right) \end{aligned} \tag{4.5}$$

with $\bar{\delta} > 0$ being an additional parameter of the controller. It is shown in [6] that for every given bounds on $\delta, d, \theta$ and for every $\varepsilon > 0$ there exist controller parameters $\lambda, \sigma, \Gamma, \bar{\delta} > 0$ such that for every initial condition the solution of the closed-loop system (4.1), (4.2), (4.3), (4.4) is bounded and satisfies the asymptotic estimate

$$\limsup_{t \to +\infty}\left(|y(t)|\right) \leq \sqrt{2\varepsilon} \tag{4.6}$$

The performance guaranteed by the controller (4.3), (4.4) for the case

$$\delta(t) \equiv 0.5,\ d(t) = 0.6\sin(2t) \tag{4.7}$$

was compared with the performance guaranteed by adaptive controllers designed using methodologies in the papers [15, 18].

System (4.1) is a system of the form (2.1), (2.2) and a DADS partial-state controller can be designed even if $\theta$ is time-varying. More, specifically, let $a, c, \Gamma, \kappa, \varepsilon, Q > 0$, $r \in (0,4)$ be the parameters of the controller. We use the functions

$$\begin{aligned} &\Phi(w) = \frac{8a}{\kappa(4-r)} w^2 \\ &V(\xi, z) = \frac{1}{2}\xi_1^2 + \frac{1}{2}\zeta^2 \end{aligned} \tag{4.8}$$

where

$$\begin{aligned} &\zeta = \xi_2 + \left(2c + \frac{5\rho}{4a} + \frac{\rho}{2}\right)\xi_1 + \frac{2Q+\rho}{2}\xi_1^3 + \frac{\rho}{4a}\xi_1^5 \\ &\rho = \kappa + \exp(z) \end{aligned} \tag{4.9}$$

in order to design a DADS controller of the form

$$\begin{aligned} &u = k(\xi, z) \\ &\dot{z} = \Gamma\exp(-z)\left(V(\xi, z) - \varepsilon\right)^{+} \end{aligned} \tag{4.10}$$

that guarantees the following differential inequalities:

$$\dot{\Phi} \le -\left(2-\frac{r}{2}-\frac{8a}{\kappa Q(4-r)}\right)\Phi(w)+Q\xi_1^4+\frac{16a}{\kappa r(4-r)}\delta^2$$
$$\dot{V} \le -2cV(\xi,z)-Q\xi_1^4+\frac{2a}{\rho}\left(d^2+\left((|\theta|-\rho)^+\right)^2+\frac{\kappa(4-r)}{8a}\Phi(w)\right) \tag{4.11}$$

Here instead of applying the methodology of the proof of Theorem 1, we can give a simple locally Lipschitz partial-state DADS feedback law that guarantees inequalities (4.11):

$$u=-\xi_1-A(\xi,z)\xi_2-\Gamma B(\xi,z)\xi_1\left(V(\xi,z)-\varepsilon\right)^+-\left(c+L(\xi,z)\right)\zeta$$
$$\dot{z}=\Gamma\exp(-z)\left(V(\xi,z)-\varepsilon\right)^+ \tag{4.12}$$

where

$$A(\xi,z)=2c+\frac{5\rho}{4a}+\frac{\rho}{2}+3\frac{2Q+\rho}{2}\xi_1^2+\frac{5\rho}{4a}\xi_1^4$$
$$B(\xi,z)=\frac{5+2a}{4a}+\frac{1}{2}\xi_1^2+\frac{1}{4a}\xi_1^4$$
$$L(\xi,z)=\frac{\rho}{4a}\left(5+\xi_1^4+\frac{a\rho}{c}\xi_1^2\right)A^2(\xi,z)$$

Controller (4.12) is a cancellation-type controller while the general recursion in Lemma 1 and Theorem 1 dominates rather than cancels. The distinction between the present example and the general case and the departure from the general design procedure is made for pedagogical reasons but also because formula (4.12) is simpler and avoid conservativeness.

Exploiting the function

$$U(w,\xi,z)=V(\xi,z)+\Phi(w) \tag{4.13}$$

and assuming that

$$Q>\frac{32}{(4-r)^2}\frac{a}{\kappa}$$

we obtain from (4.11), (4.13) and (4.10) the differential inequality

$$\dot{U} \le -\mu U(w,\xi,z)+\frac{16a}{\kappa r(4-r)}\delta^2+\frac{2a}{\kappa}\left(d^2+\left((|\theta|-\rho)^+\right)^2\right) \tag{4.14}$$

with $\mu=\min\left(\frac{Q\kappa(4-r)^2-32a}{4\kappa Q(4-r)},2c\right)>0$. Inequalities (4.11), (4.14) guarantee (see the proof of Theorem 1) that for every initial condition and for every $d,\delta,\theta\in L^\infty(\mathbb{R}_+)$ the solution of the closed-loop system (4.1), (4.2), (4.12) is bounded and satisfies the asymptotic estimate (4.6).

The difference between the controller (4.3), (4.4) and the DADS controller (4.12) is clear: for every given $\varepsilon>0$ the controller (4.3), (4.4) achieves the control objective (4.6) for given $d, \delta, \theta \in L^{\infty}(\mathbb{R}_+)$ if $\theta$ is constant and we know a constant $R > 0$ for which $\|d\|_{\infty} + \|\delta\|_{\infty} + |\theta| \leq R$. Indeed, if such a constant $R > 0$ is known then we can tune the controller (4.3), (4.4) so that the control objective (4.6) is achieved. On the other hand, the DADS controller (4.12) achieves the control objective (4.6) always even if $\theta$ is time-varying and without any knowledge of the magnitude of $\|d\|_{\infty} + \|\delta\|_{\infty} + \|\theta\|_{\infty}$. The only requirement for DADS is that the inequality $Q > \frac{32}{(4-r)^2} \frac{a}{\kappa}$ must hold for some $r \in (0,4)$ and this inequality involves parameters of the controller which are known.

As remarked above, the DADS controller (4.12) achieves many more things than the p-UBIBS property (2.20), (2.21) (notice that for this example the input coefficients are not affected by an external disturbance and the $b$'s are all identically 1) and the zero p-OAG property (4.6). The differential inequality (4.14) implies that every solution of the closed-loop system (4.1), (4.2), (4.12) satisfies the following estimate for all $t \geq 0$:

$$U\left(w(t), \xi(t), z(t)\right) \leq \exp(-\mu t) U\left(w(t), \xi(t), z(t)\right) + \frac{2a}{\mu\kappa}\left( \|d\|_{\infty}^2 + \left( \left( \|\theta\|_{\infty} - \kappa \right)^{+} \right)^2 + \frac{8}{r(4-r)} \|\delta\|_{\infty}^2 \right) \tag{4.15}$$

Estimate (4.15) shows that when $\theta$ is constant then the closed-loop system (4.1), (4.2), (4.12) satisfies the p-IOS property from the inputs $\delta, d$ with output $y = (\xi_1, w)$. As noted in [10], estimate (4.15) in conjunction with definitions (4.8), (4.13) show the first phase of convergence where the output $y = (\xi_1, w)$ converges uniformly (and here exponentially) to the set $\left\{ (\xi_1, w) \in \mathbb{R}^2 : \xi_1^2 + w^2 \leq \frac{4a}{\mu\kappa}\left( \|d\|_{\infty}^2 + \left( \left( \|\theta\|_{\infty} - \kappa \right)^{+} \right)^2 + \frac{8}{r(4-r)} \|\delta\|_{\infty}^2 \right) \right\}$. Then the second (non-uniform) phase of convergence continues in order to achieve the zero p-OAG property (4.6) and there is possibly a third (again non-uniform) phase of convergence where the whole plant state is led to zero, i.e., $\lim_{t \to +\infty}\left(|\xi(t)|\right) = \lim_{t \to +\infty}\left(|w(t)|\right) = 0$. Indeed, estimate (4.15) guarantees that the third phase of convergence will certainly take place if $\lim_{t \to +\infty}\left(\delta(t)\right) = \lim_{t \to +\infty}\left(d(t)\right) = 0$ and $\kappa \geq \limsup_{t \to +\infty}\left(|\theta(t)|\right)$. This type of convergence is a characteristic feature of DADS, as noted in [10].

We next evaluate the two controllers with the parameters used in [6] for (4.3), (4.4), i.e., the parameters $\bar{\delta} = 0.2785$, $\lambda = \sigma = \Gamma = 1$. For the DADS controller (4.12) we use the parameters

$$a = \kappa = 0.08,\ c = 0.04,\ Q = 2.01,\ \Gamma = 2$$

We also use the same initial condition that was used in [6] for (4.3), (4.4), i.e.,

$$w(0) = \xi_1(0) = \xi_2(0) = r(0) = 1,\ \hat{\theta}(0) = 0.5,\ p(0) = 0$$

For the DADS controller (4.12) we use the initial condition

$$w(0) = \xi_1(0) = \xi_2(0) = 1,\ z(0) = -2\ln(10)$$

For both controllers we select $\varepsilon = 1/32$ so that $\sqrt{2\varepsilon} = 0.25$. Moreover, $\theta$ is constant with value $\theta = 0.1$ (the same value was used in [6] and [18]).

For the disturbance given by (4.7) both controllers achieve the control objective (4.6) with $\sqrt{2\varepsilon} = 0.25$. This is shown in Fig. 1.

For $t \geq 1$ the control inputs are very similar. This is shown in Fig.2. However, for $t \in [0,1]$ the DADS controller is more aggressive: we have $\max_{t\in[0,1]}\left(|u(t)|\right) = 1651.94$ for the DADS controller (4.12) while we have $\max_{t\in[0,1]}\left(|u(t)|\right) = 152.57$ for the controller (4.3), (4.4).

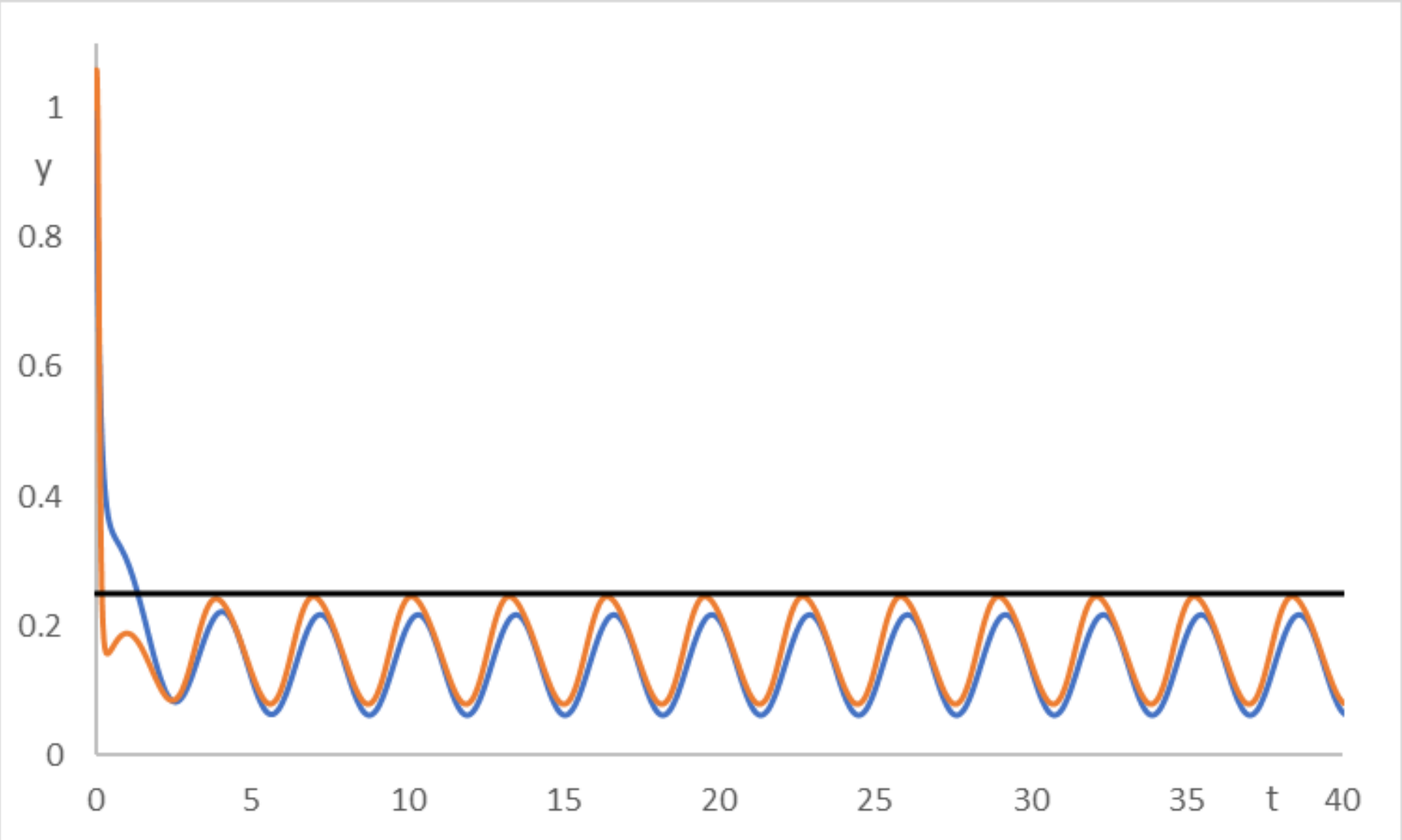

**Fig. 1:** The output $y(t)$ for the controller (4.3), (4.4) (red line) and the DADS controller (4.12) (blue line). Disturbance given by (4.7). The black line indicates $\sqrt{2\varepsilon} = 0.25$.

When the disturbance changes, namely when

$$\delta(t) \equiv 2,\ d(t) = 0.6\sin(2t) \tag{4.16}$$

then the controller (4.3), (4.4) fails to achieve the control objective (4.6) with $\sqrt{2\varepsilon} = 0.25$ while the DADS controller (4.12) achieves the control objective. This is shown in Fig. 3. It is clear that in the case (4.16) the controller (4.3), (4.4) must be tuned in a different way in order to achieve the control objective (4.6) with $\sqrt{2\varepsilon} = 0.25$. $\triangleleft$

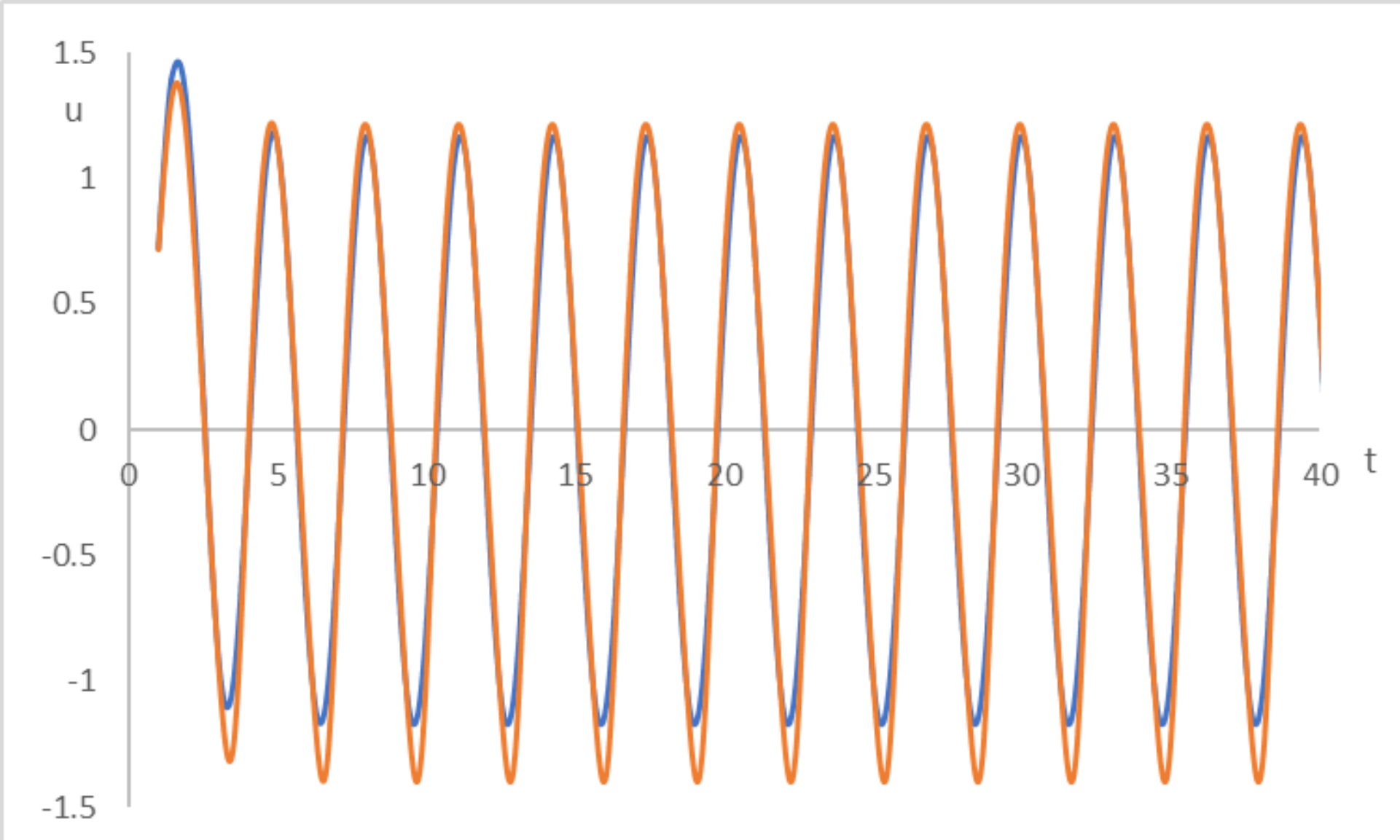

**Fig. 2:** The control input $u(t)$ for the controller (4.3), (4.4) (red line) and the DADS controller (4.12) (blue line). Disturbance given by (4.7).

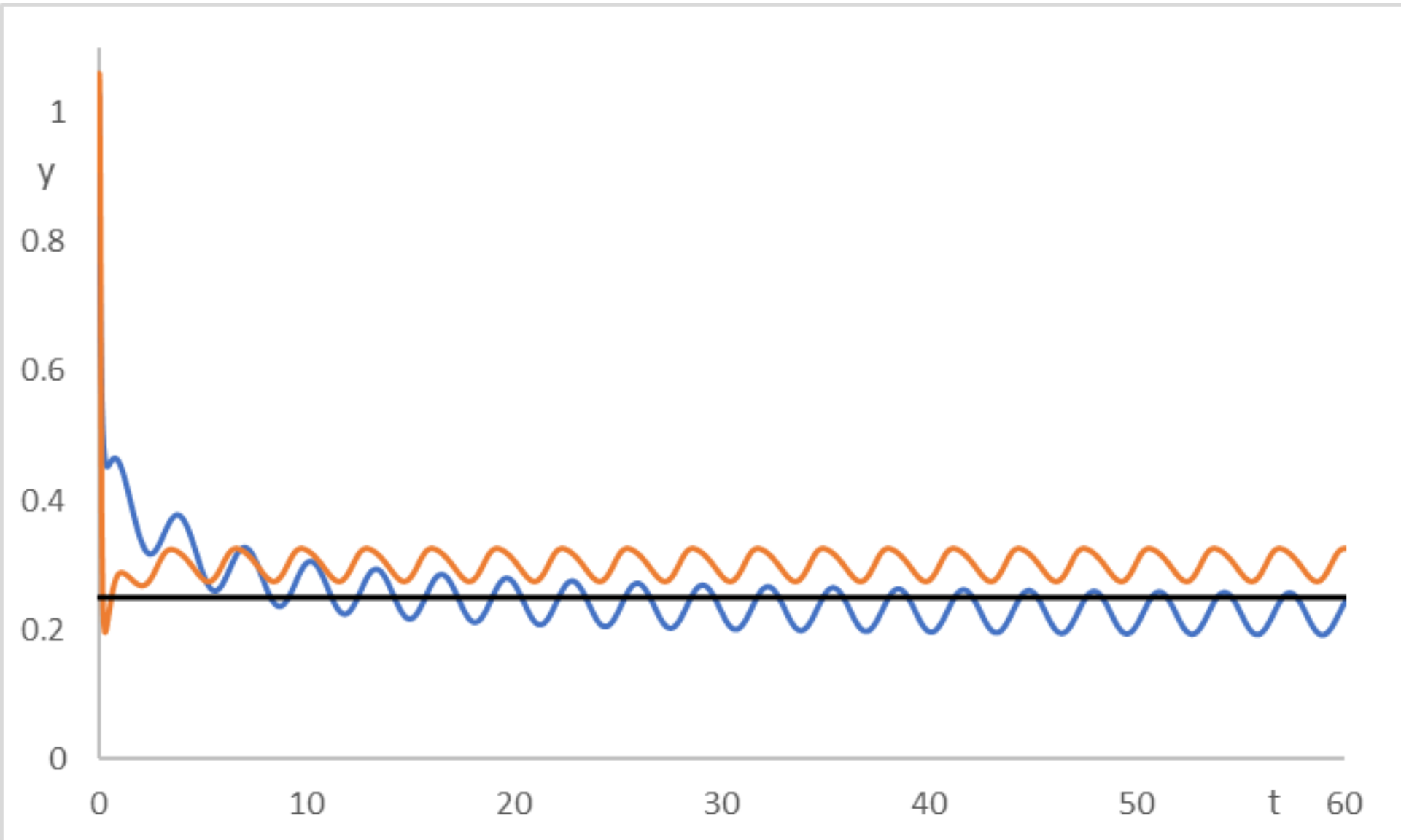

**Fig. 3:** The outpput $y(t)$ for the controller (4.3), (4.4) (red line) and the DADS controller (4.12) (blue line). Disturbance given by (4.16). The black line indicates $\sqrt{2\varepsilon} = 0.25$.

Our next example shows that partial-state DADS is not limited to finite-dimensional systems. We are dealing with a parabolic PDE but other kinds of PDEs (e.g., first-order or second-order hyperbolic PDEs) can be handled equally well (just like in [13]).

**Example 2:** Consider the following system

$$\begin{aligned}
&w_t(t,x) = \bar{p}w_{xx}(t,x) + \theta_1(t,x)\zeta(t) \\
&w(t,0) = w(t,1) = 0 \\
&\dot{\zeta}(t) = \xi_1(t) \\
&\dot{\xi}_1(t) = \langle \theta_2[t], w[t] \rangle + b(t)u(t) + d(t)
\end{aligned} \tag{4.17}$$

where $x \in [0,1]$ is the spatial variable, $\zeta(t), \xi_1(t) \in \mathbb{R}$, $w[t] \in H^2(0,1) \cap H_0^1(0,1)$ are the state components, $u(t) \in \mathbb{R}$ is the control input, $d(t) \in \mathbb{R}$ is a disturbance, $\theta_1[t], \theta_2[t] \in L^2(0,1)$, $b(t) > 0$ are unknown parameters, $\bar{p} > 0$ is an unknown constant parameter and $\langle \theta_2[t], w[t] \rangle$ denotes the inner product in $L^2(0,1)$, i.e., $\langle \theta_2[t], w[t] \rangle = \int_0^1 (\theta_2[t])(x)(w[t])(x)dx = \int_0^1 \theta_2(t,x)w(t,x)dx$.

In what follows we do not study existence/uniqueness issues for system (4.17) and we focus only on control issues. We seek a robust adaptive controller that depends only on $(\zeta, \xi_1) \in \mathbb{R}^2$ and guarantees robust practical regulation of the output $y = (\zeta, \xi_1) \in \mathbb{R}^2$ of system (4.17).

Theorem 4 can be applied to system (4.17). However, we intend to provide an explicit feedback for system (4.17) that guarantees (3.17), (3.18) and (3.19) for appropriate $\Xi \in C^0\left(\mathbb{R}_+^8 \times (0,+\infty)\right)$.

In what follows $\|\cdot\|$ denotes the norm of $L^2(0,1)$. We start with the functional

$$\Phi(w) = \|w\|^2 / (2\bar{p}) \tag{4.18}$$

which is defined for all $w \in X = L^2(0,1)$ and satisfies for every $\beta \in (0,\pi^2), C \geq 0$, $w \in S = H^2(0,1) \cap H_0^1(0,1)$ (using Wirtinger's inequality):

$$\begin{aligned}
\dot{\Phi} &= \langle w, w_{xx} \rangle + \frac{1}{\bar{p}} \langle \theta_1, w \rangle \zeta \leq -\|w_x\|^2 + \frac{\|\theta_1\|}{\bar{p}} \|w\| |\zeta| \\
&\leq -\pi^2 \|w\|^2 + (\pi^2 - \beta)\|w\|^2 + \frac{\|\theta_1\|^2}{4(\pi^2 - \beta)\bar{p}^2} \zeta^2 \\
&\leq -\beta \|w\|^2 + \left( \frac{\|\theta_1\|^2}{4(\pi^2 - \beta)\bar{p}^2} - C \right)^+ \zeta^2 + C\zeta^2 \\
&\leq -\beta \|w\|^2 + \frac{1}{4}\left( \left( \frac{\|\theta_1\|^2}{4(\pi^2 - \beta)\bar{p}^2} - C \right)^+ \right)^2 + \zeta^4 + C\zeta^2
\end{aligned} \tag{4.19}$$

Consequently, (3.13) holds with $R(w) = \beta \|w\|^2$, $Q(\zeta, \xi_1) = \zeta^4 + C\zeta^2$, $\bar{\gamma}(s) = \frac{1}{4}\left( \left( \frac{\|\theta_1\|^2}{4(\pi^2 - \beta)\bar{p}^2} - C \right)^+ \right)^2$ and $\gamma(s) \equiv 0$. Let $\varepsilon, c, \Gamma, G, H > 0$ be arbitrary constants with

$$G \geq \max\left(\frac{2C}{c}, \frac{2\left(\sqrt{c^2+1}+c\right)}{\sqrt{c^2+1}-c}\right) \tag{4.20}$$

$$H \geq \frac{4}{3c} \quad , \quad \kappa > a$$

We claim that the DADS controller

$$u = -P(\xi_1, \zeta)(\xi_1 + 2c\zeta) \tag{4.21}$$

with (2.25) and

$$V(\xi_1, \zeta) := G\bar{V}(\xi_1, \zeta) + H\bar{V}^2(\xi_1, \zeta) \tag{4.22}$$

$$\bar{V}(\xi_1, \zeta) := \frac{1}{2}\zeta^2 + \frac{1}{2}(\xi_1 + 2c\zeta)^2 \tag{4.23}$$

$$P(\xi_1, \zeta) \geq \frac{\rho S(\xi_1, \zeta, z)}{G + 2H\bar{V}(\xi_1, \zeta)}\left(\frac{S(\xi_1, \zeta, z)}{4a}(\xi_1 + 2c\zeta)^2 + 1\right) \tag{4.24}$$

$$\begin{aligned} S(\xi_1, \zeta, z) \geq{}& \left(5G + 9H\bar{V}(\xi_1, \zeta)\right)\frac{c}{2} + \frac{\left(1-4c^2\right)^2}{4c}\left(G + 2H\bar{V}(\xi_1, \zeta)\right) \\ &+ \left(\frac{2\rho^2}{\beta} + 1\right)\frac{\rho}{4a}\left(G + 2H\bar{V}(\xi_1, \zeta)\right)^2 + \frac{\rho^3}{16a^3\beta^2}\left(G + 2H\bar{V}(\xi_1, \zeta)\right)^4(\xi_1 + 2c\zeta)^2 \end{aligned} \tag{4.25}$$

where $\rho = \kappa + \exp(z)$ achieves all the conclusions of Theorem 4 with

$$n = q = m = 1,\ \Theta = \mho = \left(L^2(0,1)\right)^2 \text{ with norm } \|(\theta_1, \theta_2)\|_{\mho} = \sqrt{\|\theta_1\|^2 + \|\theta_2\|^2},$$
$$B_1 = \mathbb{R},\ X = L^2(0,1),\ S = H^2(0,1) \cap H_0^1(0,1)$$

while $\delta \in L^\infty(\mathbb{R}_+; \mathcal{E})$ is irrelevant for this example.

The constructions of the functions $V, \bar{V}$ in (4.22), (4.23) is made by following the methodology in the proof of Theorem 2. Indeed, definition (4.23) implies the estimate

$$\bar{V}(\xi_1, \zeta) \geq \frac{\sqrt{c^2+1}-c}{2\left(\sqrt{c^2+1}+c\right)}\left(\xi_1^2 + \zeta^2\right)$$

which in conjunction with (4.20) and definition (4.22) guarantees the estimate

$$V(\xi_1, \zeta) \geq \xi_1^2 + \zeta^2 \tag{4.26}$$

Moreover, by combining the inequality

$$\left(1-4c^2\right)\zeta(\xi_1 + 2c\zeta) \leq c\zeta^2 + \frac{\left(1-4c^2\right)^2}{4c}(\xi_1 + 2c\zeta)^2$$

with (4.17), (4.21), (4.22), (4.23) we obtain:

$$
\begin{aligned}
\dot{V} &= -cV(\xi_1,\zeta) - \zeta^4 - C\zeta^2 + \left( G\frac{c}{2} + H\frac{c}{2}\bar{V}(\xi_1,\zeta) + \zeta^2 + C \right)\zeta^2 \\
&-2c\left(G + 2H\bar{V}(\xi_1,\zeta)\right)\zeta^2 + \left(5G + 9H\bar{V}(\xi_1,\zeta)\right)\frac{c}{2}(\xi_1 + 2c\zeta)^2 \\
&+\left(G + 2H\bar{V}(\xi_1,\zeta)\right)(\xi_1 + 2c\zeta)\left(\langle\theta_2, w\rangle + bu + d + (1-4c^2)\zeta\right) \\
&\le -cV(\xi_1,\zeta) - \zeta^4 - C\zeta^2 + \left( \zeta^2 + C - G\frac{c}{2} - H\frac{3c}{2}\bar{V}(\xi_1,\zeta) \right)\zeta^2 \\
&+\left(5G + 9H\bar{V}(\xi_1,\zeta)\right)\frac{c}{2}(\xi_1 + 2c\zeta)^2 + \frac{(1-4c^2)^2}{4c}\left(G + 2H\bar{V}(\xi_1,\zeta)\right)(\xi_1 + 2c\zeta)^2 \\
&+\left(G + 2H\bar{V}(\xi_1,\zeta)\right)|\xi_1 + 2c\zeta|\|\theta_2\|\|w\| + \left(G + 2H\bar{V}(\xi_1,\zeta)\right)|\xi_1 + 2c\zeta||d| \\
&-b\left(G + 2H\bar{V}(\xi_1,\zeta)\right)P(\xi_1,\zeta)(\xi_1 + 2c\zeta)^2
\end{aligned} \tag{4.27}
$$

Exploiting (4.20) and definition (4.23) we get from (4.27):

$$
\begin{aligned}
\dot{V} &\le -cV(\xi_1,\zeta) - \zeta^4 - C\zeta^2 - b\left(G + 2H\bar{V}(\xi_1,\zeta)\right)P(\xi_1,\zeta)(\xi_1 + 2c\zeta)^2 \\
&+\left(5G + 9H\bar{V}(\xi_1,\zeta)\right)\frac{c}{2}(\xi_1 + 2c\zeta)^2 + \frac{(1-4c^2)^2}{4c}\left(G + 2H\bar{V}(\xi_1,\zeta)\right)(\xi_1 + 2c\zeta)^2 \\
&+\left(G + 2H\bar{V}(\xi_1,\zeta)\right)|\xi_1 + 2c\zeta|\|\theta_2\|\|w\| + \left(G + 2H\bar{V}(\xi_1,\zeta)\right)|\xi_1 + 2c\zeta||d|
\end{aligned} \tag{4.28}
$$

Working as in the proof of Theorem 2 and utilizing the inequalities

$$
\left(G + 2H\bar{V}(\xi_1,\zeta)\right)|\xi_1 + 2c\zeta||d| \le \frac{a}{\rho}|d|^2 + \frac{\rho}{4a}\left(G + 2H\bar{V}(\xi_1,\zeta)\right)^2(\xi_1 + 2c\zeta)^2
$$

$$
\left(G + 2H\bar{V}(\xi_1,\zeta)\right)|\xi_1 + 2c\zeta|\|\theta_2\|\|w\| \le \frac{a}{\rho}\beta\|w\|^2 + \frac{\rho}{4a\beta}\left(G + 2H\bar{V}(\xi_1,\zeta)\right)^2\|\theta_2\|^2(\xi_1 + 2c\zeta)^2
$$

$$
\|\theta_2\|^2 \le 2\rho^2 + 2\left((\|\theta_2\| - \rho)^+\right)^2
$$

$$
\begin{aligned}
&\left(G + 2H\bar{V}(\xi_1,\zeta)\right)|\xi_1 + 2c\zeta|\|\theta_2\|\|w\| \le \frac{a}{\rho}\beta\|w\|^2 + \frac{\rho^3}{2a\beta}\left(G + 2H\bar{V}(\xi_1,\zeta)\right)^2(\xi_1 + 2c\zeta)^2 \\
&+\frac{\rho}{2a\beta}\left(G + 2H\bar{V}(\xi_1,\zeta)\right)^2\left((\|\theta_2\| - \rho)^+\right)^2(\xi_1 + 2c\zeta)^2 \\
&\le \frac{a}{\rho}\left(\beta\|w\|^2 + \left((\|\theta_2\| - \rho)^+\right)^4\right) + \frac{\rho^3}{2a\beta}\left(G + 2H\bar{V}(\xi_1,\zeta)\right)^2(\xi_1 + 2c\zeta)^2 \\
&+\frac{\rho^3}{16a^3\beta^2}\left(G + 2H\bar{V}(\xi_1,\zeta)\right)^4(\xi_1 + 2c\zeta)^4
\end{aligned}
$$

we obtain from (4.25) and (4.28) the estimate

$$\begin{aligned}\dot{V} \le & -cV(\xi_1,\zeta)-\zeta^4-C\zeta^2-b\left(G+2H\bar{V}(\xi_1,\zeta)\right)P(\xi_1,\zeta)(\xi_1+2c\zeta)^2 \\ & +S(\xi_1,\zeta,z)(\xi_1+2c\zeta)^2+\frac{a}{\rho}\left(\beta\|w\|^2+\left(\left(\|\theta_2\|-\rho\right)^+\right)^4+|d|^2\right)\end{aligned} \tag{4.29}$$

The inequalities

$$1=b\left(\frac{1}{b}-\rho\right)+b\rho\le b\left(\frac{1}{b}-\rho\right)^+ +b\rho$$

$$\left(\frac{1}{b}-\rho\right)^+ S(\xi_1,\zeta,z)(\xi_1+2c\zeta)^2\le\frac{a}{\rho}\left(\left(\frac{1}{b}-\rho\right)^+\right)^2+\frac{\rho}{4a}S^2(\xi_1,\zeta,z)(\xi_1+2c\zeta)^4$$

and (4.24) guarantee the following differential inequality:

$$\dot{V}\le -cV(\xi,\zeta)-\zeta^4-C\zeta^2+\frac{a}{\rho}\left(\beta\|w\|^2+\left(\left(\|\theta\|-\rho\right)^+\right)^4+|d|^2+b\left(\left(\frac{1}{b}-\rho\right)^+\right)^2\right) \tag{4.30}$$

The conclusions of Theorem 4 are consequences of the differential inequalities (4.19), (4.30) and inequalities (4.20), (4.26). ⊲

## 5. Proofs of Main Results

The following backstepping lemma deals with the parameterized control system

$$\begin{aligned}\dot{\xi}&=f(w,\xi,b,\theta,d)+g(w,\xi,b,\theta)\chi \\ \dot{\chi}&=h(w,\xi,\chi,b,\theta,d)+\tilde{g}(w,\xi,\chi,\theta)\tilde{b}u\end{aligned} \tag{5.1}$$

where $(\xi,\chi)\in\mathbb{R}^n\times\mathbb{R}$ is the state, $w\in\mathbb{R}^l$ is the vector of the parameters, $u\in\mathbb{R}$ is the control input and $\theta\in\Theta\subseteq\mathbb{R}^p$, $\Theta\subseteq\mathbb{R}^p$ is a non-empty closed set, $d\in\mathbb{R}^m$, $b=(b_1,...,b_q)\in(0,+\infty)^q\cap B$, $B\subseteq\mathbb{R}^q$ is a non-empty closed set, $\tilde{b}\in(0,+\infty)$ are unknown disturbances, $f:\mathbb{R}^l\times\mathbb{R}^n\times\left((0,+\infty)^q\cap B\right)\times\Theta\times\mathbb{R}^m\to\mathbb{R}^n$, $h:\mathbb{R}^l\times\mathbb{R}^n\times\mathbb{R}\times\left((0,+\infty)^q\cap B\right)\times\Theta\times\mathbb{R}^m\to\mathbb{R}$, $g:\mathbb{R}^l\times\mathbb{R}^n\times\left((0,+\infty)^q\cap B\right)\times\mathbb{R}^p\to\mathbb{R}^n$, $\tilde{g}:\mathbb{R}^l\times\mathbb{R}^n\times\mathbb{R}\times\Theta\to\mathbb{R}$ are smooth mappings with $f(0,0,b,\theta,0)=0$, $h(0,0,0,b,\theta,0)=0$ for all $\theta\in\Theta$, $b\in(0,+\infty)^q\cap B$.

**Lemma 1:** *Consider system (5.1) and assume that there exist constants* $a,c,\Gamma,\varepsilon,\kappa>0$, $\Lambda\ge 0$, *smooth functions* $\sigma,k,V\in C^\infty(\mathbb{R}^n\times\mathbb{R})$, $Q\in C^\infty(\mathbb{R}^n\times\mathbb{R};\mathbb{R}_+)$, $R\in C^\infty(\mathbb{R}^l;\mathbb{R}_+)$, $\mu,\eta\in C^\infty\left(\mathbb{R}^n\times\mathbb{R};(0,+\infty)\right)$ *with* $V(0,z)=k(0,z)=0$ *for* $z\in\mathbb{R}$ *and* $V(\xi,z)>0$ *for* $\xi\ne 0$, $z\in\mathbb{R}$, *such that the following inequalities hold for all* $\theta \in \Theta$, $d \in \mathbb{R}^m$, $b \in (0,+\infty)^q \cap B$, $(w,\xi,\chi,z) \in \mathbb{R}^l \times \mathbb{R}^n \times \mathbb{R} \times \mathbb{R}$:

$$|\xi|^2 \leq \sigma(\xi,z)V(\xi,z) \tag{5.2}$$

$$\begin{aligned}&\frac{\partial V}{\partial \xi}(\xi,z)\big(f(w,\xi,b,\theta,d)+g(w,\xi,b,\theta)k(\xi,z)\big)+\frac{\partial V}{\partial z}(\xi,z)\Gamma\exp(-z)\big(V(\xi,z)-\varepsilon\big)^+\\&\leq -cV(\xi,z)-Q(\xi,z)+\frac{a}{\rho}\left(\min_{i=1,\dots,q}(b_i)\left(\left(\frac{1}{\min_{i=1,\dots,q}(b_i)}-\rho\right)^+\right)^2+|d|^2+|d|^4+R(w)+\Lambda\right)\\&+\frac{a}{\rho}\left(\left(\left(|\theta|+|b|+|b||\theta|-\rho\right)^+\right)^2+\left(\left(|\theta|+|b|+|b||\theta|-\rho\right)^+\right)^4+\left(\left(|\theta|+|b|+|b||\theta|-\rho\right)^+\right)^8\right)\end{aligned} \tag{5.3}$$

$$\tilde{g}(w,\xi,\chi,\theta) \geq \eta(\xi,\chi) \tag{5.4}$$

$$\begin{aligned}&|f(w,\xi,b,\theta,d)|+|h(w,\xi,\chi,b,\theta,d)|\\&\leq (1+|\theta|)(1+|b|)\mu(\xi,\chi)\left(1+\sqrt{R(w)+\Lambda}\right)|d|\\&+(1+|\theta|)(1+|b|)\mu(\xi,\chi)\left(|(\xi,\chi)|+\sqrt{R(w)+\Lambda}\right)\end{aligned} \tag{5.5}$$

$$|g(w,\xi,b,\theta)| \leq (1+|\theta|)(1+|b|)\mu(\xi,\chi)\left(1+\sqrt{R(w)+\Lambda}\right) \tag{5.6}$$

*where* $\rho := \kappa + \exp(z)$. *Define the following function for all* $(\xi,\chi,z) \in \mathbb{R}^n \times \mathbb{R} \times \mathbb{R}$:

$$\bar{V}(\xi,\chi,z) := V(\xi,z) + \frac{1}{2}\big(\chi - k(\xi,z)\big)^2 \tag{5.7}$$

*Then there exist functions* $\bar{\sigma} \in C^\infty(\mathbb{R}^n \times \mathbb{R} \times \mathbb{R})$, $M \in C^\infty(\mathbb{R}^n \times \mathbb{R} \times \mathbb{R};(0,+\infty))$ *such that, with*

$$\bar{k}(\xi,\chi,z) := -\frac{M(\xi,\chi,z)}{\eta(\xi,\chi)}\big(\chi - k(\xi,z)\big), \text{ for all } (\xi,\chi,z) \in \mathbb{R}^n \times \mathbb{R} \times \mathbb{R} \tag{5.8}$$

*the following inequalities hold for all* $\theta \in \Theta$, $d \in \mathbb{R}^m$, $b \in (0,+\infty)^q \cap B$, $\tilde{b} \in (0,+\infty)$, $(w,\xi,\chi,z) \in \mathbb{R}^l \times \mathbb{R}^n \times \mathbb{R} \times \mathbb{R}$:

$$|\xi|^2 + |\chi|^2 \leq \bar{\sigma}(\xi,\chi,z)\bar{V}(\xi,\chi,z) \tag{5.9}$$

$$\begin{aligned}
&\frac{\partial \bar{V}}{\partial \xi}(\xi,\chi,z)\left(f(w,\xi,b,\theta,d)+g(w,\xi,b,\theta)\chi\right)\\
&+\frac{\partial \bar{V}}{\partial z}(\xi,\chi,z)\Gamma\exp(-z)\left(\bar{V}(\xi,\chi,z)-\varepsilon\right)^{+}\\
&+\frac{\partial \bar{V}}{\partial \chi}(\xi,\chi,z)\left(h(w,\xi,\chi,b,\theta,d)+\tilde{g}(w,\xi,\chi,\theta)\tilde{b}\,\bar{k}(\xi,\chi,z)\right)\\
&\leq -\frac{c}{2}\bar{V}(\xi,\chi,z)-Q(\xi,z)+\frac{2a}{\rho}\left(\min\left(\min_{i=1,\dots,q}(b_i),\tilde{b}\right)\left(\left(\frac{1}{\min\left(\min_{i=1,\dots,q}(b_i),\tilde{b}\right)}-\rho\right)^{+}\right)^{2}\right)\\
&+\frac{2a}{\rho}\left(\left(\left(|\theta|+|b|+|b||\theta|-\rho\right)^{+}\right)^{2}+\left(\left(|\theta|+|b|+|b||\theta|-\rho\right)^{+}\right)^{4}+\left(\left(|\theta|+|b|+|b||\theta|-\rho\right)^{+}\right)^{8}\right)\\
&+\frac{2a}{\rho}\left(|d|^{2}+|d|^{4}+R(w)+\Lambda\right)
\end{aligned} \tag{5.10}$$

*where* $\rho:=\kappa+\exp(z)$. *The functions* $\bar{V},\bar{k}\in C^{\infty}\left(\mathbb{R}^n\times\mathbb{R}\times\mathbb{R}\right)$ *defined by (5.7) and (5.8) satisfy* $\bar{V}(0,0,z)=\bar{k}(0,0,z)=0$ *for* $z\in\mathbb{R}$ *and* $\bar{V}(\xi,\chi,z)>0$ *for* $(\xi,\chi)\neq 0$, $z\in\mathbb{R}$. *Moreover, if the set* $\left\{(\xi,z)\in\mathbb{R}^n\times\left[-\tilde{M},\tilde{M}\right]:V(\xi,z)\leq\tilde{M}\right\}$ *is bounded for each* $\tilde{M}>0$ *then the set* $\left\{(\xi,\chi,z)\in\mathbb{R}^n\times\mathbb{R}\times\left[-\tilde{M},\tilde{M}\right]:\bar{V}(\xi,\chi,z)\leq\tilde{M}\right\}$ *is bounded for each* $\tilde{M}>0$.

**Proof of Lemma 1:** Without loss of generality, we may assume that $\sigma(\xi,z)>0$ for all $\xi\in\mathbb{R}^n$, $z\in\mathbb{R}$.

We start by using the fact that for every smooth function $\psi:\mathbb{R}^n\times\mathbb{R}\to\mathbb{R}^m$ with $\psi(0,z)=0$ for all $z\in\mathbb{R}$ there exists a smooth positive function $L_\psi:\mathbb{R}^n\times\mathbb{R}\to(0,+\infty)$ such that the inequality $|\psi(\xi,z)|\leq L_\psi(\xi,z)|\xi|$ holds for all $(\xi,z)\in\mathbb{R}^n\times\mathbb{R}$. Since $V(0,z)=0$ for all $z\in\mathbb{R}$ and $V(\xi,z)>0$ for all $\xi\neq 0$, $z\in\mathbb{R}$, it follows that for each fixed $z\in\mathbb{R}$, the mapping $\xi\mapsto V(\xi,z)$ has a global minimum at $\xi=0$. Consequently, $\frac{\partial V}{\partial \xi}(0,z)=0$ for all $z\in\mathbb{R}$. Moreover, since $k(0,z)=0$ for all $z\in\mathbb{R}$, we conclude that there exists a smooth positive functions $S:\mathbb{R}^n\times\mathbb{R}\to(0,+\infty)$ such that the following inequality hold for all $(\xi,z)\in\mathbb{R}^n\times\mathbb{R}$:

$$\left|\frac{\partial V}{\partial \xi}(\xi,z)\right|+|k(\xi,z)|\leq S(\xi,z)|\xi| \tag{5.11}$$

Define for all $(\xi,\chi,z)\in\mathbb{R}^n\times\mathbb{R}\times\mathbb{R}$

$$\bar{\sigma}(\xi,\chi,z):=\left(1+2S^2(\xi,z)\right)\sigma(\xi,z)+4 \tag{5.12}$$

$$M(\xi,\chi,z):=\rho P(\xi,\chi,z)\left(1+\frac{P(\xi,\chi,z)}{4a}(\chi-k(\xi,z))^2\right) \tag{5.13}$$

where $\rho:=\kappa+\exp(z)$ and $P\in C^\infty\left(\mathbb{R}^n\times\mathbb{R}\times\mathbb{R};(0,+\infty)\right)$ is a function that is selected appropriately below.

In what follows, $\theta\in\Theta$, $d\in\mathbb{R}^m$, $b\in(0,+\infty)^q\cap B$, $\tilde{b}\in(0,+\infty)$, $(w,\xi,\chi,z)\in\mathbb{R}^l\times\mathbb{R}^n\times\mathbb{R}\times\mathbb{R}$ are arbitrary given vectors/scalars.

We next proceed to show the validity of (5.9) with $\bar{\sigma}(\xi,\chi,z)$ defined by (5.12). Using the triangle inequality $|\chi|\le|\chi-k(\xi,z)|+|k(\xi,z)|$ and the fact that $\left(|\chi-k(\xi,z)|+|k(\xi,z)|\right)^2\le 2|\chi-k(\xi,z)|^2+2|k(\xi,z)|^2$ we get:

$$|\xi|^2+|\chi|^2\le|\xi|^2+2|\chi-k(\xi,z)|^2+2|k(\xi,z)|^2$$

The above inequality in conjunction with definition (5.7) and inequality (5.11) gives:

$$|\xi|^2+|\chi|^2\le\left(1+2S^2(\xi,z)\right)|\xi|^2+4\bar{V}(\xi,\chi,z) \tag{5.14}$$

Using (5.7) and inequalities (5.2), (5.14) we obtain

$$|\xi|^2+|\chi|^2\le\left(\left(1+2S^2(\xi,z)\right)\sigma(\xi,z)+4\right)\bar{V}(\xi,\chi,z) \tag{5.15}$$

Combining (5.12) and (5.15) we obtain (5.9).

In what follows, we denote by $\frac{d\bar{V}}{dt}$ the quantity appearing in the left-hand side of (5.10). For brevity, we also use the following notation:

$$\zeta=\chi-k(\xi,z) \tag{5.16}$$

Before we proceed with the derivation of (5.10) we need the two following facts. Both facts are proved at the end of this proof.

FACT I: For every function $E\in C^\infty\left(\mathbb{R}^n\times\mathbb{R}\times\mathbb{R};(0,+\infty)\right)$ and for every $s\ge 0$ the following inequality holds:

$$\begin{aligned}
s|g(w,\xi,b,\theta)||\zeta| &\le \frac{E(\xi,\chi,z)}{2}s^2+\frac{a}{6\rho}\left(R(w)+\Lambda\right)\\
&+\frac{a}{4\rho}\left(\left(|\theta|+|b|+|b||\theta|-\rho\right)^+\right)^2+\frac{a}{4\rho}\left(\left(|\theta|+|b|+|b||\theta|-\rho\right)^+\right)^4\\
&+\frac{9\rho^3}{a^3}s^4\mu^4(\xi,\chi)\zeta^4+\left(\frac{\rho}{a}s^2+\frac{(1+\rho)^2}{2E(\xi,\chi,z)}+\frac{3\rho}{a}s^2(1+\rho)^2\right)\mu^2(\xi,\chi)\zeta^2
\end{aligned} \tag{5.17}$$

FACT II: For every function $E \in C^{\infty}\left(\mathbb{R}^n \times \mathbb{R} \times \mathbb{R};(0,+\infty)\right)$ and for every $s \geq 0$ the following inequality holds:

$$\begin{aligned}
& s(1+|\theta|)(1+|b|)\mu(\xi,\chi)\left(\left(1+\sqrt{R(w)+\Lambda}\right)|d|+|(\xi,\chi)|+\sqrt{R(w)+\Lambda}\right)|\zeta| \\
& \leq \frac{a}{3\rho}\left(R(w)+\Lambda\right)+\frac{a}{2\rho}\left(|d|^2+|d|^4\right)+\frac{E(\xi,\chi,z)}{2}|(\xi,\chi)|^2 \\
& +\frac{a}{4\rho}\left(\left(\left(|\theta|+|b|+|b|\left||\theta|-\rho\right)^+\right)^2+\left(\left(|\theta|+|b|+|b|\left||\theta|-\rho\right)^+\right)^4\right)+\frac{81\rho^7}{2a^7}s^8\mu^8(\xi,\chi)\zeta^8 \qquad (5.18) \\
& +\frac{\rho^3}{a^3}\left(9(1+\rho)^4+16\right)s^4\mu^4(\xi,\chi)\zeta^4+\frac{a}{2\rho}\left(\left(|\theta|+|b|+|b|\left||\theta|-\rho\right)^+\right)^8 \\
& +\left(\frac{4\rho}{a}(1+\rho)^2+\frac{(1+\rho)^2}{2E(\xi,\chi,z)}+\frac{\rho}{a}|(\xi,\chi)|^2\right)s^2\mu^2(\xi,\chi)\zeta^2
\end{aligned}$$

We are ready to proceed with the derivation of (5.10). Using (5.7), (5.8), (5.16) we get:

$$\begin{aligned}
\frac{d\bar{V}}{dt} &= \frac{\partial V}{\partial \xi}(\xi,z)\left(f(w,\xi,b,\theta,d)+g(w,\xi,b,\theta)\chi\right) \\
&+\frac{\partial V}{\partial z}(\xi,z)\Gamma\exp(-z)\left(\bar{V}(\xi,\chi,z)-\varepsilon\right)^+ \\
&+\zeta h(w,\xi,\chi,b,\theta,d)-\tilde{b}\tilde{g}(w,\xi,\chi,\theta)\frac{M(\xi,\chi,z)}{\eta(\xi,\chi)}\zeta^2 \qquad (5.19) \\
&-\zeta\frac{\partial k}{\partial \xi}(\xi,z)\left(f(w,\xi,b,\theta,d)+g(w,\xi,b,\theta)\chi\right) \\
&-\zeta\frac{\partial k}{\partial z}(\xi,z)\Gamma\exp(-z)\left(\bar{V}(\xi,\chi,z)-\varepsilon\right)^+
\end{aligned}$$

By virtue of (5.4) we obtain from (5.19):

$$\begin{aligned}
\frac{d\bar{V}}{dt} &\leq \frac{\partial V}{\partial \xi}(\xi,z)\left(f(w,\xi,b,\theta,d)+g(w,\xi,b,\theta)k(\xi,z)\right) \\
&+\frac{\partial V}{\partial z}(\xi,z)\Gamma\exp(-z)\left(V(\xi,z)-\varepsilon\right)^+-\zeta\frac{\partial k}{\partial z}(\xi,z)\Gamma\exp(-z)\left(\bar{V}(\xi,\chi,z)-\varepsilon\right)^+ \\
&+\frac{\partial V}{\partial \xi}(\xi,z)g(w,\xi,b,\theta)\zeta+\frac{\partial V}{\partial z}(\xi,z)\Gamma\exp(-z)\left(\left(\bar{V}(\xi,\chi,z)-\varepsilon\right)^+-\left(V(\xi,z)-\varepsilon\right)^+\right) \qquad (5.20) \\
&+\zeta h(w,\xi,\chi,b,\theta,d)-\tilde{b}M(\xi,\chi,z)\zeta^2-\zeta\frac{\partial k}{\partial \xi}(\xi,z)\left(f(w,\xi,b,\theta,d)+g(w,\xi,b,\theta)\chi\right)
\end{aligned}$$

Combining (5.3) and (5.20) we get:

$$
\begin{aligned}
\frac{d\bar{V}}{dt} \leq & -cV(\xi,z) - Q(\xi,z) + \frac{a}{\rho}\min_{i=1,\dots,q}(b_i)\left(\left(\frac{1}{\min_{i=1,\dots,q}(b_i)} - \rho\right)^{+}\right)^{2} - \tilde{b}M(\xi,\chi,z)\zeta^2 \\
& + \frac{a}{\rho}\left(\left((|\theta|+|b|+|b||\theta|-\rho)^{+}\right)^{4} + \left((|\theta|+|b|+|b||\theta|-\rho)^{+}\right)^{8}\right) + |h(w,\xi,\chi,b,\theta,d)||\zeta| \\
& + \frac{a}{\rho}\left(\left((|\theta|+|b|+|b||\theta|-\rho)^{+}\right)^{2} + |d|^2 + |d|^4 + R(w) + \Lambda\right) + |\zeta|\left|\frac{\partial k}{\partial z}(\xi,z)\right|\Gamma\exp(-z)\left(\bar{V}(\xi,\chi,z)-\varepsilon\right)^{+} \\
& + \left|\frac{\partial V}{\partial \xi}(\xi,z)\right||g(w,\xi,b,\theta)||\zeta| + \Gamma\left|\frac{\partial V}{\partial z}(\xi,z)\right|\exp(-z)\left|\left(\bar{V}(\xi,\chi,z)-\varepsilon\right)^{+} - \left(V(\xi,z)-\varepsilon\right)^{+}\right| \\
& + |\zeta|\left|\frac{\partial k}{\partial \xi}(\xi,z)\right||f(w,\xi,b,\theta,d)| + |\zeta|\left|\frac{\partial k}{\partial \xi}(\xi,z)\right||g(w,\xi,b,\theta)||\chi|
\end{aligned}
\tag{5.21}
$$

For the term $\left|\frac{\partial V}{\partial \xi}(\xi,z)\right||g(w,\xi,b,\theta)||\zeta|$ we use inequality (5.17) with $s = \left|\frac{\partial V}{\partial \xi}(\xi,z)\right|$, while for the term $|\zeta|\left|\frac{\partial k}{\partial \xi}(\xi,z)\right||g(w,\xi,b,\theta)||\chi|$ we use inequality (5.17) with $s = \left|\frac{\partial k}{\partial \xi}(\xi,z)\right||\chi|$, both for an appropriate function $E \in C^{\infty}\left(\mathbb{R}^n \times \mathbb{R} \times \mathbb{R};(0,+\infty)\right)$ (to be selected) and we obtain:

$$
\begin{aligned}
& \left|\frac{\partial V}{\partial \xi}(\xi,z)\right||g(w,\xi,b,\theta)||\zeta| + |\zeta|\left|\frac{\partial k}{\partial \xi}(\xi,z)\right||g(w,\xi,b,\theta)||\chi| \\
& \leq \frac{E(\xi,\chi,z)}{2}\left(\left|\frac{\partial V}{\partial \xi}(\xi,z)\right|^2 + \left|\frac{\partial k}{\partial \xi}(\xi,z)\right|^2|\chi|^2\right) + \frac{a}{3\rho}\left(R(w)+\Lambda\right) \\
& + \frac{a}{2\rho}\left(\left((|\theta|+|b|+|b||\theta|-\rho)^{+}\right)^{2} + \frac{a}{2\rho}\left((|\theta|+|b|+|b||\theta|-\rho)^{+}\right)^{4}\right. \\
& + \left(\frac{\rho}{a}\left(\left|\frac{\partial V}{\partial \xi}(\xi,z)\right|^2 + \left|\frac{\partial k}{\partial \xi}(\xi,z)\right|^2|\chi|^2\right)\left(1+3(1+\rho)^2\right) + \frac{(1+\rho)^2}{E(\xi,\chi,z)}\right)\mu^2(\xi,\chi)\zeta^2 \\
& + \frac{9\rho^3}{a^3}\left(\left|\frac{\partial V}{\partial \xi}(\xi,z)\right|^4 + \left|\frac{\partial k}{\partial \xi}(\xi,z)\right|^4|\chi|^4\right)\mu^4(\xi,\chi)\zeta^4
\end{aligned}
\tag{5.22}
$$

Exploiting (5.5), (5.18) with $s=1$ for the term $|h(w,\xi,\chi,b,\theta,d)||\zeta|$ and (5.5), (5.18) with $s = \left|\frac{\partial k}{\partial \xi}(\xi,z)\right|$ for the term $|\zeta|\left|\frac{\partial k}{\partial \xi}(\xi,z)\right||f(w,\xi,b,\theta,d)|$, both for an appropriate function $E \in C^{\infty}\left(\mathbb{R}^n \times \mathbb{R} \times \mathbb{R};(0,+\infty)\right)$ (to be selected), we obtain:

$$
\begin{aligned}
&\left|h(w,\xi,\chi,b,\theta,d)\right|\left|\zeta\right| + \left|\zeta\right|\left|\frac{\partial k}{\partial \xi}(\xi,z)\right|\left|f(w,\xi,b,\theta,d)\right| \\
&\leq \frac{2a}{3\rho}\left(R(w)+\Lambda\right) + \frac{a}{\rho}\left(|d|^2+|d|^4\right) + \frac{81\rho^7}{2a^7}\left(1+\left|\frac{\partial k}{\partial \xi}(\xi,z)\right|^8\right)\mu^8(\xi,\chi)\zeta^8 \\
&+\frac{a}{2\rho}\left(\left(\left(|\theta|+|b|+|b|\left|\theta\right|-\rho\right)^+\right)^2 + \left(\left(|\theta|+|b|+|b|\left|\theta\right|-\rho\right)^+\right)^4\right) + \frac{a}{\rho}\left(\left(|\theta|+|b|+|b|\left|\theta\right|-\rho\right)^+\right)^8 \\
&+\frac{\rho^3}{a^3}\left(9(1+\rho)^4+16\right)\left(1+\left|\frac{\partial k}{\partial \xi}(\xi,z)\right|^4\right)\mu^4(\xi,\chi)\zeta^4 + E(\xi,\chi,z)\left|(\xi,\chi)\right|^2 \\
&+\left(\frac{4\rho}{a}(1+\rho)^2 + \frac{(1+\rho)^2}{2E(\xi,\chi,z)} + \frac{\rho}{a}\left|(\xi,\chi)\right|^2\right)\left(1+\left|\frac{\partial k}{\partial \xi}(\xi,z)\right|^2\right)\mu^2(\xi,\chi)\zeta^2
\end{aligned}
\tag{5.23}
$$

Combining (5.21), (5.22) and (5.23) we get:

$$
\begin{aligned}
&\frac{d\bar{V}}{dt} \leq -cV(\xi,z) + \frac{a}{\rho}\min_{i=1,\ldots,q}(b_i)\left(\left(\frac{1}{\min\limits_{i=1,\ldots,q}(b_i)}-\rho\right)^+\right)^2 + |\zeta|\left|\frac{\partial k}{\partial z}(\xi,z)\right|\Gamma\exp(-z)\left(\bar{V}(\xi,\chi,z)-\varepsilon\right)^+ \\
&+\frac{2a}{\rho}\left(\left(\left(|\theta|+|b|+|b|\left|\theta\right|-\rho\right)^+\right)^8 + \left(\left(|\theta|+|b|+|b|\left|\theta\right|-\rho\right)^+\right)^4 + \left(\left(|\theta|+|b|+|b|\left|\theta\right|-\rho\right)^+\right)^2\right) \\
&-Q(\xi,z) + \Gamma\left|\frac{\partial V}{\partial z}(\xi,z)\right|\exp(-z)\left|\left(\bar{V}(\xi,\chi,z)-\varepsilon\right)^+ - \left(V(\xi,z)-\varepsilon\right)^+\right| - \tilde{b}M(\xi,\chi,z)\zeta^2 \\
&+\frac{E(\xi,\chi,z)}{2}\left(\left|\frac{\partial V}{\partial \xi}(\xi,z)\right|^2 + \left|\frac{\partial k}{\partial \xi}(\xi,z)\right|^2|\chi|^2 + 2\left|(\xi,\chi)\right|^2\right) + G_3(\xi,\chi,z)\mu^8(\xi,\chi)\zeta^8 \\
&+G_1(\xi,\chi,z)\mu^2(\xi,\chi)\zeta^2 + G_2(\xi,\chi,z)\mu^4(\xi,\chi)\zeta^4 + \frac{2a}{\rho}\left(|d|^2+|d|^4+R(w)+\Lambda\right)
\end{aligned}
\tag{5.24}
$$

where

$$
\begin{aligned}
G_1(\xi,\chi,z) &:= \frac{\rho}{a}\left(1+3(1+\rho)^2\right)\left(\left|\frac{\partial V}{\partial \xi}(\xi,z)\right|^2 + \left|\frac{\partial k}{\partial \xi}(\xi,z)\right|^2|\chi|^2\right) + \frac{(1+\rho)^2}{E(\xi,\chi,z)} \\
&+\left(\frac{4\rho}{a}(1+\rho)^2 + \frac{(1+\rho)^2}{2E(\xi,\chi,z)} + \frac{\rho}{a}\left|(\xi,\chi)\right|^2\right)\left(1+\left|\frac{\partial k}{\partial \xi}(\xi,z)\right|^2\right)
\end{aligned}
\tag{5.25}
$$

$$
G_2(\xi,\chi,z) := \frac{\rho^3}{a^3}\left(9\left(\left|\frac{\partial V}{\partial \xi}(\xi,z)\right|^4 + \left|\frac{\partial k}{\partial \xi}(\xi,z)\right|^4|\chi|^4\right) + \left(9(1+\rho)^4+16\right)\left(1+\left|\frac{\partial k}{\partial \xi}(\xi,z)\right|^4\right)\right) \tag{5.26}
$$

$$
G_3(\xi,\chi,z) := \frac{81\rho^7}{2a^7}\left(1+\left|\frac{\partial k}{\partial \xi}(\xi,z)\right|^8\right) \tag{5.27}
$$

Exploiting (5.11), the facts that $\left|v^+-\varphi^+\right|\le|v-\varphi|$ for all $v,\varphi\in\mathbb{R}$, $\left(\bar{V}(\xi,\chi,z)-\varepsilon\right)^+\le\bar{V}(\xi,\chi,z)$ and (5.7), (5.16), we obtain from (5.24):

$$\begin{aligned}
\frac{d\bar{V}}{dt}\le&-cV(\xi,z)+\frac{a}{\rho}\min_{i=1,\dots,q}\left(b_i\right)\left(\left(\frac{1}{\min\limits_{i=1,\dots,q}\left(b_i\right)}-\rho\right)^+\right)^2+|\zeta|\left|\frac{\partial k}{\partial z}(\xi,z)\right|\Gamma\exp(-z)V(\xi,z)\\
&+\frac{2a}{\rho}\left(\left(\left(|\theta|+|b|+|b||\theta|-\rho\right)^+\right)^8+\left(\left(|\theta|+|b|+|b||\theta|-\rho\right)^+\right)^4+\left(\left(|\theta|+|b|+|b||\theta|-\rho\right)^+\right)^2\right)\\
&-Q(\xi,z)+\frac{\Gamma}{2}\left|\frac{\partial V}{\partial z}(\xi,z)\right|\exp(-z)\zeta^2-\tilde{b}M(\xi,\chi,z)\zeta^2+\frac{\Gamma}{2}|\zeta|^3\left|\frac{\partial k}{\partial z}(\xi,z)\right|\exp(-z)\\
&+\frac{E(\xi,\chi,z)}{2}\left(S^2(\xi,z)|\xi|^2+\left|\frac{\partial k}{\partial \xi}(\xi,z)\right|^2|\chi|^2+2\left|(\xi,\chi)\right|^2\right)+G_3(\xi,\chi,z)\mu^8(\xi,\chi)\zeta^8\\
&+G_1(\xi,\chi,z)\mu^2(\xi,\chi)\zeta^2+G_2(\xi,\chi,z)\mu^4(\xi,\chi)\zeta^4+\frac{2a}{\rho}\left(|d|^2+|d|^4+R(w)+\Lambda\right)
\end{aligned}\qquad(5.28)$$

Using the inequalities

$$\Gamma|\zeta|\left|\frac{\partial k}{\partial z}(\xi,z)\right|\exp(-z)V(\xi,z)\le\frac{c}{10}V(\xi,z)+\frac{5\Gamma^2}{2c}\left|\frac{\partial k}{\partial z}(\xi,z)\right|^2\exp(-2z)V(\xi,z)\zeta^2$$

$$\frac{\Gamma}{2}|\zeta|^3\left|\frac{\partial k}{\partial z}(\xi,z)\right|\exp(-z)\le\frac{\Gamma}{4}\left(1+\left|\frac{\partial k}{\partial z}(\xi,z)\right|^2\right)\left(1+\zeta^2\right)\exp(-z)\zeta^2$$

$$\frac{\Gamma}{2}\left|\frac{\partial V}{\partial z}(\xi,z)\right|\exp(-z)\zeta^2\le\frac{\Gamma}{4}\left(1+\left|\frac{\partial V}{\partial z}(\xi,z)\right|^2\right)\exp(-z)\zeta^2$$

we obtain from (5.28) and the facts that $|\xi|^2\le\left|(\xi,\chi)\right|^2$, $|\chi|^2\le\left|(\xi,\chi)\right|^2$:

$$\begin{aligned}
\frac{d\bar{V}}{dt}\le&-\frac{9c}{10}V(\xi,z)+\frac{a}{\rho}\min_{i=1,\dots,q}\left(b_i\right)\left(\left(\frac{1}{\min\limits_{i=1,\dots,q}\left(b_i\right)}-\rho\right)^+\right)^2-Q(\xi,z)-\tilde{b}M(\xi,\chi,z)\zeta^2\\
&+\frac{2a}{\rho}\left(\left(\left(|\theta|+|b|+|b||\theta|-\rho\right)^+\right)^8+\left(\left(|\theta|+|b|+|b||\theta|-\rho\right)^+\right)^4+\left(\left(|\theta|+|b|+|b||\theta|-\rho\right)^+\right)^2\right)\\
&+G_4(\xi,\chi,z)\zeta^2+G_1(\xi,\chi,z)\mu^2(\xi,\chi)\zeta^2+G_2(\xi,\chi,z)\mu^4(\xi,\chi)\zeta^4\\
&+\frac{E(\xi,\chi,z)}{2}\left(S^2(\xi,z)+\left|\frac{\partial k}{\partial \xi}(\xi,z)\right|^2+2\right)\left|(\xi,\chi)\right|^2+G_3(\xi,\chi,z)\mu^8(\xi,\chi)\zeta^8\\
&+\frac{2a}{\rho}\left(|d|^2+|d|^4+R(w)+\Lambda\right)
\end{aligned}\qquad(5.29)$$

where

$$G_4(\xi,\chi,z)$$
$$:=\frac{\Gamma\exp(-z)}{4}\left(1+\left|\frac{\partial V}{\partial z}(\xi,z)\right|^2+\left(1+\left|\frac{\partial k}{\partial z}(\xi,z)\right|^2\right)\left(1+\zeta^2\right)+\frac{10\Gamma}{c}\left|\frac{\partial k}{\partial z}(\xi,z)\right|^2\exp(-z)V(\xi,z)\right) \tag{5.30}$$

Using (5.9), (5.7) and (5.16) we get from (5.29):

$$\begin{aligned}
\frac{d\bar{V}}{dt} &\le -\frac{c}{2}\bar{V}(\xi,\chi,z)+\frac{c}{4}\zeta^2-\frac{2c}{5}V(\xi,z)+\frac{a}{\rho}\min_{i=1,\dots,q}(b_i)\left(\left(\frac{1}{\min\limits_{i=1,\dots,q}(b_i)}-\rho\right)^+\right)^2-Q(\xi,z)\\
&+\frac{2a}{\rho}\left(\left((|\theta|+|b|+|b||\theta|-\rho)^+\right)^8+\left((|\theta|+|b|+|b||\theta|-\rho)^+\right)^4+\left((|\theta|+|b|+|b||\theta|-\rho)^+\right)^2\right)\\
&-\tilde{b}M(\xi,\chi,z)\zeta^2+G_4(\xi,\chi,z)\zeta^2+G_1(\xi,\chi,z)\mu^2(\xi,\chi)\zeta^2+G_2(\xi,\chi,z)\mu^4(\xi,\chi)\zeta^4\\
&+\frac{E(\xi,\chi,z)}{2}\left(S^2(\xi,z)+\left|\frac{\partial k}{\partial \xi}(\xi,z)\right|^2+2\right)\bar{\sigma}(\xi,\chi,z)V(\xi,z)+G_3(\xi,\chi,z)\mu^8(\xi,\chi)\zeta^8\\
&+\frac{E(\xi,\chi,z)}{4}\left(S^2(\xi,z)+\left|\frac{\partial k}{\partial \xi}(\xi,z)\right|^2+2\right)\bar{\sigma}(\xi,\chi,z)\zeta^2+\frac{2a}{\rho}\left(|d|^2+|d|^4+R(w)+\Lambda\right)
\end{aligned} \tag{5.31}$$

Selecting

$$E(\xi,\chi,z):=\frac{4c}{5\left(S^2(\xi,z)+\left|\frac{\partial k}{\partial \xi}(\xi,z)\right|^2+2\right)\bar{\sigma}(\xi,\chi,z)} \tag{5.32}$$

and setting $P(\xi,\chi,z)$ to be any smooth function that satisfies

$$\begin{aligned}
P(\xi,\chi,z)&\ge\frac{9c}{20}+G_4(\xi,\chi,z)+G_1(\xi,\chi,z)\mu^2(\xi,\chi)\\
&+G_2(\xi,\chi,z)\mu^4(\xi,\chi)\zeta^2+G_3(\xi,\chi,z)\mu^8(\xi,\chi)\zeta^6
\end{aligned} \tag{5.33}$$

we get from (5.31):

$$\begin{aligned}
\frac{d\bar{V}}{dt} &\le -\frac{c}{2}\bar{V}(\xi,\chi,z)-Q(\xi,z)+\frac{a}{\rho}\min_{i=1,\dots,q}(b_i)\left(\left(\frac{1}{\min\limits_{i=1,\dots,q}(b_i)}-\rho\right)^+\right)^2\\
&+\frac{2a}{\rho}\left(\left((|\theta|+|b|+|b||\theta|-\rho)^+\right)^8+\left((|\theta|+|b|+|b||\theta|-\rho)^+\right)^4+\left((|\theta|+|b|+|b||\theta|-\rho)^+\right)^2\right)\\
&-\tilde{b}M(\xi,\chi,z)\zeta^2+P(\xi,\chi,z)\zeta^2+\frac{2a}{\rho}\left(|d|^2+|d|^4+R(w)+\Lambda\right)
\end{aligned} \tag{5.34}$$

The fact that the map $(0,+\infty)\ni s\mapsto s\left(\left(\frac{1}{s}-\rho\right)^+\right)^2$ is non-increasing allows us to obtain from (5.34):

$$\begin{aligned}
&\frac{d\bar{V}}{dt}\leq\frac{a}{\rho}\min\left(\min_{i=1,...,q}(b_i),\tilde{b}\right)\left(\left(\frac{1}{\min\left(\min_{i=1,...,q}(b_i),\tilde{b}\right)}-\rho\right)^+\right)^2\\
&+\frac{2a}{\rho}\left(\left(\left(|\theta|+|b|+|b||\theta|-\rho\right)^+\right)^8+\left(\left(|\theta|+|b|+|b||\theta|-\rho\right)^+\right)^4+\left(\left(|\theta|+|b|+|b||\theta|-\rho\right)^+\right)^2\right)\\
&-\frac{c}{2}\bar{V}(\xi,\chi,z)-Q(\xi,z)-\min\left(\min_{i=1,...,q}(b_i),\tilde{b}\right)M(\xi,\chi,z)\zeta^2\\
&+\min\left(\min_{i=1,...,q}(b_i),\tilde{b}\right)\frac{1}{\min\left(\min_{i=1,...,q}(b_i),\tilde{b}\right)}P(\xi,\chi,z)\zeta^2+\frac{2a}{\rho}\left(|d|^2+|d|^4+R(w)+\Lambda\right)
\end{aligned}\tag{5.35}$$

Inequality (5.10) is a consequence of definition (5.13), inequality (5.35) and the following inequalities:

$$\frac{1}{\min\left(\min_{i=1,...,q}(b_i),\tilde{b}\right)}\leq\rho+\left(\frac{1}{\min\left(\min_{i=1,...,q}(b_i),\tilde{b}\right)}-\rho\right)^+$$

$$\left(\frac{1}{\min\left(\min_{i=1,...,q}(b_i),\tilde{b}\right)}-\rho\right)^+P(\xi,\chi,z)\zeta^2\leq\frac{a}{\rho}\left(\left(\frac{1}{\min\left(\min_{i=1,...,q}(b_i),\tilde{b}\right)}-\rho\right)^+\right)^2+\frac{\rho}{4a}P^2(\xi,\chi,z)\zeta^4$$

We finish the proof by showing the validity of (5.17) and (5.18). Using the inequality

$$\begin{aligned}
&s\left(1+|\theta|\right)\left(1+|b|\right)\mu(\xi,\chi)|\zeta|\sqrt{R(w)+\Lambda}\\
&\leq\frac{a}{6\rho}\left(R(w)+\Lambda\right)+\frac{3\rho}{2a}s^2\left(1+|\theta|\right)^2\left(1+|b|\right)^2\mu^2(\xi,\chi)\zeta^2
\end{aligned}$$

we obtain from (5.6) the inequality

$$\begin{aligned}
&s|g(w,\xi,b,\theta)||\zeta|\leq s\left(1+|\theta|\right)\left(1+|b|\right)\mu(\xi,\chi)|\zeta|\\
&+s\left(1+|\theta|\right)\left(1+|b|\right)\mu(\xi,\chi)|\zeta|\sqrt{R(w)+\Lambda}\\
&\leq s\left(1+|\theta|\right)\left(1+|b|\right)\mu(\xi,\chi)|\zeta|\\
&+\frac{a}{6\rho}\left(R(w)+\Lambda\right)+\frac{3\rho}{2a}s^2\left(1+|\theta|\right)^2\left(1+|b|\right)^2\mu^2(\xi,\chi)\zeta^2
\end{aligned}\tag{5.36}$$

Similarly, using the inequalities

$$s(1+|\theta|)(1+|b|)\mu(\xi,\chi)\sqrt{R(w)+\Lambda}\,|d||\zeta|$$
$$\leq \frac{a}{6\rho}(R(w)+\Lambda)+\frac{3\rho}{2a}s^2(1+|\theta|)^2(1+|b|)^2\mu^2(\xi,\chi)|d|^2\zeta^2$$

$$s(1+|\theta|)(1+|b|)\mu(\xi,\chi)\sqrt{R(w)+\Lambda}\,|\zeta|$$
$$\leq \frac{a}{6\rho}(R(w)+\Lambda)+\frac{3\rho}{2a}s^2(1+|\theta|)^2(1+|b|)^2\mu^2(\xi,\chi)\zeta^2$$

we get:

$$\begin{aligned}
&s(1+|\theta|)(1+|b|)\mu(\xi,\chi)\left(1+\sqrt{R(w)+\Lambda}\right)|d||\zeta| \\
&+s(1+|\theta|)(1+|b|)\mu(\xi,\chi)\left(|(\xi,\chi)|+\sqrt{R(w)+\Lambda}\right)|\zeta| \\
&\leq \frac{a}{3\rho}(R(w)+\Lambda)+\frac{3\rho}{2a}s^2(1+|\theta|)^2(1+|b|)^2\mu^2(\xi,\chi)|d|^2\zeta^2+s(1+|\theta|)(1+|b|)\mu(\xi,\chi)|d||\zeta| \\
&+\frac{3\rho}{2a}s^2(1+|\theta|)^2(1+|b|)^2\mu^2(\xi,\chi)\zeta^2+s|(\xi,\chi)|(1+|\theta|)(1+|b|)\mu(\xi,\chi)|\zeta|
\end{aligned} \tag{5.37}$$

The facts $|\theta|+|b|+|b||\theta|\leq \rho+(|\theta|+|b|+|b||\theta|-\rho)^+$ and $(1+|\theta|)^2(1+|b|)^2\leq 2(1+\rho)^2+2\left((|\theta|+|b|+|b||\theta|-\rho)^+\right)^2$ in conjunction with (5.36), (5.37) give:

$$\begin{aligned}
&s|g(w,\xi,b,\theta)||\zeta|\leq s(1+\rho)\mu(\xi,\chi)|\zeta|+s(|\theta|+|b|+|b||\theta|-\rho)^+\mu(\xi,\chi)|\zeta| \\
&+\frac{a}{6\rho}(R(w)+\Lambda)+\frac{3\rho}{a}s^2(1+\rho)^2\mu^2(\xi,\chi)\zeta^2 \\
&+\frac{3\rho}{a}s^2\left((|\theta|+|b|+|b||\theta|-\rho)^+\right)^2\mu^2(\xi,\chi)\zeta^2
\end{aligned} \tag{5.38}$$

$$\begin{aligned}
&s(1+|\theta|)(1+|b|)\mu(\xi,\chi)\left(1+\sqrt{R(w)+\Lambda}\right)|d||\zeta| \\
&+s(1+|\theta|)(1+|b|)\mu(\xi,\chi)\left(|(\xi,\chi)|+\sqrt{R(w)+\Lambda}\right)|\zeta| \\
&\leq \frac{a}{3\rho}(R(w)+\Lambda)+\frac{3\rho}{a}s^2(1+\rho)^2\mu^2(\xi,\chi)|d|^2\zeta^2+s|(\xi,\chi)|(1+\rho)\mu(\xi,\chi)|\zeta| \\
&+\frac{3\rho}{a}s^2\left((|\theta|+|b|+|b||\theta|-\rho)^+\right)^2\mu^2(\xi,\chi)|d|^2\zeta^2+s(1+\rho)\mu(\xi,\chi)|d||\zeta| \\
&+s(|\theta|+|b|+|b||\theta|-\rho)^+\mu(\xi,\chi)|d||\zeta|+s|(\xi,\chi)|(|\theta|+|b|+|b||\theta|-\rho)^+\mu(\xi,\chi)|\zeta| \\
&+\frac{3\rho}{a}s^2(1+\rho)^2\mu^2(\xi,\chi)\zeta^2+\frac{3\rho}{a}s^2\left((|\theta|+|b|+|b||\theta|-\rho)^+\right)^2\mu^2(\xi,\chi)\zeta^2
\end{aligned} \tag{5.39}$$

The inequalities

$$s(1+\rho)\mu(\xi,\chi)|\zeta| \le \frac{E(\xi,\chi,z)}{2}s^2 + \frac{1}{2E(\xi,\chi,z)}(1+\rho)^2\mu^2(\xi,\chi)\zeta^2$$

$$s\left(|\theta|+|b|+|b||\theta|-\rho\right)^+\mu(\xi,\chi)|\zeta| \le \frac{a}{4\rho}\left(\left(|\theta|+|b|+|b||\theta|-\rho\right)^+\right)^2 + \frac{\rho}{a}s^2\mu^2(\xi,\chi)\zeta^2$$

$$\begin{aligned}&\frac{3\rho}{a}s^2\left(\left(|\theta|+|b|+|b||\theta|-\rho\right)^+\right)^2\mu^2(\xi,\chi)\zeta^2\\ &\le \frac{a}{4\rho}\left(\left(|\theta|+|b|+|b||\theta|-\rho\right)^+\right)^4 + \frac{9\rho^3}{a^3}s^4\mu^4(\xi,\chi)\zeta^4\end{aligned}$$

in conjunction with (5.38) show that inequality (5.17) holds.

The inequalities

$$s(1+\rho)\mu(\xi,\chi)|d||\zeta| \le \frac{a}{4\rho}|d|^2 + \frac{\rho}{a}s^2(1+\rho)^2\mu^2(\xi,\chi)\zeta^2$$

$$s|(\xi,\chi)|(1+\rho)\mu(\xi,\chi)|\zeta| \le \frac{E(\xi,\chi,z)}{2}|(\xi,\chi)|^2 + \frac{(1+\rho)^2}{2E(\xi,\chi,z)}s^2\mu^2(\xi,\chi)\zeta^2$$

$$s\left(|\theta|+|b|+|b||\theta|-\rho\right)^+\mu(\xi,\chi)|d||\zeta| \le \frac{a}{4\rho}|d|^2 + \frac{\rho}{a}s^2\left(\left(|\theta|+|b|+|b||\theta|-\rho\right)^+\right)^2\mu^2(\xi,\chi)\zeta^2$$

in conjunction with (5.39) show that following inequality holds:

$$\begin{aligned}&s(1+|\theta|)(1+|b|)\mu(\xi,\chi)\left(1+\sqrt{R(w)+\Lambda}\right)|d||\zeta|\\ &+s(1+|\theta|)(1+|b|)\mu(\xi,\chi)\left(|(\xi,\chi)|+\sqrt{R(w)+\Lambda}\right)|\zeta|\\ &\le \frac{a}{3\rho}(R(w)+\Lambda)+\frac{3\rho}{a}s^2(1+\rho)^2\mu^2(\xi,\chi)|d|^2\zeta^2+\frac{E(\xi,\chi,z)}{2}|(\xi,\chi)|^2\\ &+\frac{3\rho}{a}s^2\left(\left(|\theta|+|b|+|b||\theta|-\rho\right)^+\right)^2\mu^2(\xi,\chi)|d|^2\zeta^2+\frac{a}{2\rho}|d|^2\\ &+s|(\xi,\chi)|\left(|\theta|+|b|+|b||\theta|-\rho\right)^+\mu(\xi,\chi)|\zeta|+\frac{(1+\rho)^2}{2E(\xi,\chi,z)}s^2\mu^2(\xi,\chi)\zeta^2\\ &+\frac{4\rho}{a}s^2(1+\rho)^2\mu^2(\xi,\chi)\zeta^2+\frac{4\rho}{a}s^2\left(\left(|\theta|+|b|+|b||\theta|-\rho\right)^+\right)^2\mu^2(\xi,\chi)\zeta^2\end{aligned} \tag{5.40}$$

It follows from (5.40) and the inequalities

$$\frac{3\rho}{a}s^2(1+\rho)^2\mu^2(\xi,\chi)|d|^2\zeta^2 \le \frac{a}{4\rho}|d|^4 + \frac{9\rho^3}{a^3}s^4(1+\rho)^4\mu^4(\xi,\chi)\zeta^4$$

$$\frac{3\rho}{a}s^2\left(\left(|\theta|+|b|+|b||\theta|-\rho\right)^+\right)^2\mu^2(\xi,\chi)|d|^2\zeta^2$$
$$\le\frac{a}{4\rho}|d|^4+\frac{9\rho^3}{a^3}s^4\left(\left(|\theta|+|b|+|b||\theta|-\rho\right)^+\right)^4\mu^4(\xi,\chi)\zeta^4$$

$$s|(\xi,\chi)|\left(|\theta|+|b|+|b||\theta|-\rho\right)^+\mu(\xi,\chi)|\zeta|$$
$$\le\frac{a}{4\rho}\left(\left(|\theta|+|b|+|b||\theta|-\rho\right)^+\right)^2+\frac{\rho}{a}s^2|(\xi,\chi)|^2\mu^2(\xi,\chi)\zeta^2$$

that the following estimate is valid:

$$\begin{aligned}
&s\left(1+|\theta|\right)\left(1+|b|\right)\mu(\xi,\chi)\left(1+\sqrt{R(w)+\Lambda}\right)|d||\zeta| \\
&+s\left(1+|\theta|\right)\left(1+|b|\right)\mu(\xi,\chi)\left(|(\xi,\chi)|+\sqrt{R(w)+\Lambda}\right)|\zeta| \\
&\le\frac{a}{3\rho}\left(R(w)+\Lambda\right)+\frac{a}{2\rho}\left(|d|^2+|d|^4\right)+\frac{E(\xi,\chi,z)}{2}|(\xi,\chi)|^2 \\
&+\frac{9\rho^3}{a^3}s^4(1+\rho)^4\mu^4(\xi,\chi)\zeta^4+\frac{9\rho^3}{a^3}s^4\left(\left(|\theta|+|b|+|b||\theta|-\rho\right)^+\right)^4\mu^4(\xi,\chi)\zeta^4 \\
&+\left(\frac{\rho}{a}|(\xi,\chi)|^2+\frac{(1+\rho)^2}{2E(\xi,\chi,z)}+\frac{4\rho}{a}(1+\rho)^2\right)s^2\mu^2(\xi,\chi)\zeta^2 \\
&+\frac{a}{4\rho}\left(\left(|\theta|+|b|+|b||\theta|-\rho\right)^+\right)^2+\frac{4\rho}{a}s^2\left(\left(|\theta|+|b|+|b||\theta|-\rho\right)^+\right)^2\mu^2(\xi,\chi)\zeta^2
\end{aligned}\tag{5.41}$$

Inequality (5.18) is a consequence of (5.41) and the inequalities

$$\frac{9\rho^3}{a^3}s^4\left(\left(|\theta|+|b|+|b||\theta|-\rho\right)^+\right)^4\mu^4(\xi,\chi)\zeta^4\le\frac{a}{2\rho}\left(\left(|\theta|+|b|+|b||\theta|-\rho\right)^+\right)^8+\frac{81\rho^7}{2a^7}s^8\mu^8(\xi,\chi)\zeta^8$$

$$\frac{4\rho}{a}s^2\left(\left(|\theta|+|b|+|b||\theta|-\rho\right)^+\right)^2\mu^2(\xi,\chi)\zeta^2\le\frac{a}{4\rho}\left(\left(|\theta|+|b|+|b||\theta|-\rho\right)^+\right)^4+\frac{16\rho^3}{a^3}s^4\mu^4(\xi,\chi)\zeta^4$$

The proof is complete. $\triangleleft$

We next provide the proof of Theorem 1.

**Proof of Theorem 1:** The proof is divided in four parts.

Part I: We show that for every $\kappa, a>0$ inequalities (2.9), (2.10) imply the following inequality holds for all $(\xi,w)\in\mathbb{R}^m\times\mathbb{R}^l$, $\delta\in\mathbb{R}^\alpha$, $\theta\in\Theta$, $z\in\mathbb{R}$:

$$\nabla\Phi(w)\, f(\xi_1, w, \delta, \theta) \le -\left(1-\frac{a}{\kappa}\right)R(w) + \gamma\left(|\delta|\right) + \bar{\gamma}\left(|\theta|\right)$$
$$+Q\left(\xi_1, z\right) + \frac{a}{\rho}\Lambda + \frac{a}{\rho}\left(\left(|\theta|-\rho\right)^+\right)^2 + \frac{a}{\rho}\left(\left(|\theta|-\rho\right)^+\right)^4 \tag{5.42}$$

where $\rho := \kappa + \exp(z)$ and

$$Q\left(\xi_1, z\right) := \xi_1^2\left(\psi\left(\xi_1\right) + \frac{\rho^3}{16a^3}\xi_1^2\mu_1^4(\xi_1) + \frac{\rho^3}{2a}\mu_1^2(\xi_1) + \frac{\rho}{4a}\xi_1^2\mu_1^2(\xi_1) + \rho\mu_1(\xi_1)\right) \tag{5.43}$$

Part II: We show that for every $a, c, \kappa > 0$ there exists a function $k_1 \in C^\infty\left(\mathbb{R}\times\mathbb{R}\right)$ with $k_1(0,z) = 0$ for $z \in \mathbb{R}$ such that the following inequality holds for all $\left(\xi, w\right) \in \mathbb{R}^m \times \mathbb{R}^l$, $\theta \in \Theta$, $d \in \mathbb{R}^q$, $b_1 \in \left(0, +\infty\right) \cap B_1$ and $z \in \mathbb{R}$:

$$\frac{\partial V_1}{\partial \xi_1}(\xi_1)\left(h_1(w, \xi_1, \theta, d) + b_1 g_1(w, \xi_1, \theta) k_1(\xi_1, z)\right)$$
$$\le -cV_1(\xi_1) - Q(\xi_1, z) + \frac{a}{\rho}\left(b_1\left(\left(\frac{1}{b_1} - \rho\right)^+\right)^2 + R\left(w\right) + \Lambda\right) \tag{5.44}$$
$$+\frac{a}{\rho}\left(\left(\left(|\theta|-\rho\right)^+\right)^2 + \left(\left(|\theta|-\rho\right)^+\right)^4 + \left(\left(|\theta|-\rho\right)^+\right)^8 + |d|^2 + |d|^4\right)$$

where

$$V_1\left(\xi_1\right) := \frac{1}{2}\xi_1^2 \tag{5.45}$$

Part III: Applying Lemma exactly $m-1$ times and exploiting Assumption (A), we construct for every tuple of constants $\kappa, \Gamma, \varepsilon, a, c > 0$ functions $k, V \in C^\infty\left(\mathbb{R}^m \times \mathbb{R}\right)$ with $V(0,z) = k(0,z) = 0$ for $z \in \mathbb{R}$, $V(\xi, z) > 0$ for $\xi \ne 0$, such that:

(i) for each $\tilde{M} > 0$ the set $\left\{ (\xi, z) \in \mathbb{R}^m \times \left[-\tilde{M}, \tilde{M}\right] : V(\xi, z) \le \tilde{M} \right\}$ is bounded, and

(ii) the following inequalities hold:

$$\frac{1}{2}\xi_1^2 \le V(\xi, z), \text{ for all } \left(\xi, z\right) \in \mathbb{R}^m \times \mathbb{R} \tag{5.46}$$

$$\sum_{i=1}^{m-1}\frac{\partial V}{\partial \xi_i}(\xi,z)\left(h_j(w,\xi_1,...,\xi_j,\theta,d)+b_j g_j(w,\xi_1,...,\xi_j,\theta)\xi_{j+1}\right)$$
$$+\frac{\partial V}{\partial \xi_m}(\xi,z)\left(h_m(w,\xi,\theta,d)+b_m g_m(w,\xi,\theta)k\left(\xi,z\right)\right)$$
$$+\frac{\partial V}{\partial z}(\xi,z)\Gamma\exp(-z)\left(V(\xi,z)-\varepsilon\right)^+$$
$$\leq -2^{1-m}cV(\xi,z)-Q(\xi_1,z)+\frac{2^{m-1}a}{\rho}\left(\min_{i=1,...,m}(b_i)\left(\left(\frac{1}{\min_{i=1,...,m}(b_i)}-\rho\right)^+\right)^2+|d|^2+|d|^4+R(w)+\Lambda\right)$$
$$+\frac{2^{m-1}a}{\rho}\left(\left(\left(|\theta|+|b|+|b||\theta|-\rho\right)^+\right)^2+\left(\left(|\theta|+|b|+|b||\theta|-\rho\right)^+\right)^4+\left(\left(|\theta|+|b|+|b||\theta|-\rho\right)^+\right)^8\right)$$
$$\text{for all } (\xi,w)\in\mathbb{R}^m\times\mathbb{R}^l,\ z\in\mathbb{R},\ d\in\mathbb{R}^q,\ \theta\in\Theta,\ b_j\in(0,+\infty)\cap B_j \text{ for } j=1,...,m \tag{5.47}$$

where $\rho:=\kappa+\exp(z)$.

More specifically, applying Lemma exactly $m-1$ times, we construct functions $k_j\in C^\infty\left(\mathbb{R}^j\times\mathbb{R}\right)$ for $j=2,...,m$ with $k_j(0,z)=0$ for $z\in\mathbb{R}$, such that

$$V(\xi,z)=\frac{1}{2}\xi_1^2+\frac{1}{2}\sum_{j=2}^{m}\left(\xi_j-k_{j-1}\left(\xi_1,...,\xi_{j-1},z\right)\right)^2$$

Notice that the above equation guarantees the validity of (5.46).

Part IV: We perform the stability analysis of the closed-loop system (2.1), (2.2) with (2.17), (2.18) for the functions $k,V\in C^\infty\left(\mathbb{R}^m\times\mathbb{R}\right)$ with $V(0,z)=k(0,z)=0$ for $z\in\mathbb{R}$, $V(\xi,z)>0$ for $\xi\neq 0$, that are designed in Part III for arbitrary constants $\kappa,\Gamma,\varepsilon>0$, $c=2^m$ and $a=\frac{\kappa}{2\left(1+2^{m-1}\right)}$.

Part I: Let arbitrary $\kappa,a>0$ be given. Using (2.9), (2.10) we obtain for all $(\xi,w)\in\mathbb{R}^m\times\mathbb{R}^l$, $\delta\in\mathbb{R}^\alpha$, $\theta\in\Theta$, $z\in\mathbb{R}$:

$$\begin{aligned}\nabla\Phi(w)f(\xi_1,w,\delta,\theta)&\leq -R(w)+\gamma\left(|\delta|\right)+\bar{\gamma}\left(|\theta|\right)+\xi_1^2\psi\left(\xi_1\right)+|\theta||\xi_1||\phi(\xi_1,w)|\\&\leq -R(w)+\gamma\left(|\delta|\right)+\bar{\gamma}\left(|\theta|\right)+\xi_1^2\psi\left(\xi_1\right)+|\theta||\xi_1|\mu_1(\xi_1)\sqrt{R(w)+\Lambda}+|\theta|\xi_1^2\mu_1(\xi_1)\end{aligned} \tag{5.48}$$

where $\rho:=\kappa+\exp(z)$. The inequality

$$|\theta||\xi_1|\mu_1(\xi_1)\sqrt{R(w)+\Lambda}\leq\frac{\rho}{4a}|\theta|^2\xi_1^2\mu_1^2(\xi_1)+\frac{a}{\rho}\left(R(w)+\Lambda\right)$$

in conjunction with (5.48) gives:

$$\nabla\Phi(w)f(\xi_1,w,\delta,\theta)\leq -R(w)+\gamma\left(|\delta|\right)+\bar{\gamma}\left(|\theta|\right)+\xi_1^2\psi\left(\xi_1\right)$$
$$+\frac{\rho}{4a}|\theta|^2\,\xi_1^2\mu_1^2(\xi_1)+\frac{a}{\rho}\left(R(w)+\Lambda\right)+|\theta|\,\xi_1^2\mu_1(\xi_1) \tag{5.49}$$

Exploiting (5.49) and using the facts that $\rho\geq\kappa$, $|\theta|\leq\left(|\theta|-\rho\right)^+ +\rho$, $\left(\left(|\theta|-\rho\right)^+ +\rho\right)^2\leq 2\left(\left(|\theta|-\rho\right)^+\right)^2+2\rho^2$, we get:

$$\nabla\Phi(w)f(\xi_1,w,\delta,\theta)\leq -\left(1-\frac{a}{\kappa}\right)R(w)+\gamma\left(|\delta|\right)+\bar{\gamma}\left(|\theta|\right)$$
$$+\xi_1^2\psi\left(\xi_1\right)+\left(\left(|\theta|-\rho\right)^+\right)^2\frac{\rho}{2a}\xi_1^2\mu_1^2(\xi_1) \tag{5.50}$$
$$+\frac{\rho^3}{2a}\xi_1^2\mu_1^2(\xi_1)+\frac{a}{\rho}\Lambda+\left(|\theta|-\rho\right)^+\xi_1^2\mu_1(\xi_1)+\rho\xi_1^2\mu_1(\xi_1)$$

Exploiting the inequalities

$$\left(\left(|\theta|-\rho\right)^+\right)^2\frac{\rho}{2a}\xi_1^2\mu_1^2(\xi_1)\leq\frac{a}{\rho}\left(\left(|\theta|-\rho\right)^+\right)^4+\frac{\rho^3}{16a^3}\xi_1^4\mu_1^4(\xi_1) \tag{5.51}$$

$$\left(|\theta|-\rho\right)^+\xi_1^2\mu_1(\xi_1)\leq\frac{a}{\rho}\left(\left(|\theta|-\rho\right)^+\right)^2+\frac{\rho}{4a}\xi_1^4\mu_1^2(\xi_1) \tag{5.52}$$

we obtain inequality (5.42) from (5.50) and definition (5.43).

Part II: Let arbitrary $a,c,\kappa>0$ be given. We define for all $\left(\xi_1,z\right)\in\mathbb{R}\times\mathbb{R}$:

$$k_1(\xi_1,z)=-\frac{\xi_1^2 S\left(\xi_1,z\right)+4a}{4a\eta_1(\xi_1)}\rho S\left(\xi_1,z\right)\xi_1 \tag{5.53}$$

where $\rho:=\kappa+\exp(z)$, $\eta_1(\xi_1)$ is the positive smooth function involved in (2.6) for $j=1$ and $S\left(\xi_1,z\right)$ is any smooth function that satisfies the following inequality for all $\left(\xi_1,z\right)\in\mathbb{R}\times\mathbb{R}$:

$$S\left(\xi_1,z\right)\geq\frac{c}{2}+\psi\left(\xi_1\right)+\left(1+\rho\right)^2\frac{3\rho}{2a}\mu_1^2\left(\xi_1\right)+\mu_1\left(\xi_1\right)\left(1+2\rho\right)$$
$$+\frac{5\rho^3}{8a^3}\mu_1^4\left(\xi_1\right)\xi_1^2+\frac{\rho^3}{2a}\mu_1^2(\xi_1)+\frac{\rho^3}{2a^3}\left(1+\rho\right)^4\mu_1^4\left(\xi_1\right)\xi_1^2 \tag{5.54}$$
$$+\frac{\rho}{2a}\mu_1^2\left(\xi_1\right)\xi_1^2+\frac{\rho^7}{16a^7}\mu_1^8\left(\xi_1\right)\xi_1^6$$

We next show that (5.44) holds for all $(\xi,w)\in\mathbb{R}^m\times\mathbb{R}^l$, $\theta\in\Theta$, $d\in\mathbb{R}^q$, $b_1\in(0,+\infty)\cap B_1$ and $z\in\mathbb{R}$. Using (5.45) and (5.53) we get:

$$\begin{aligned}&\frac{\partial V_1}{\partial \xi_1}(\xi_1)\left(h_1(w,\xi_1,\theta,d)+b_1 g_1(w,\xi_1,\theta)k_1(\xi_1,z)\right)\\&=\xi_1 h_1(w,\xi_1,\theta,d)+b_1 g_1(w,\xi_1,\theta)\xi_1 k_1(\xi_1,z)\\&=\xi_1 h_1(w,\xi_1,\theta,d)-b_1 g_1(w,\xi_1,\theta)\frac{\xi_1^2 S(\xi_1,z)+4a}{4a\eta_1(\xi_1)}\rho S(\xi_1,z)\xi_1^2\end{aligned}\tag{5.55}$$

Exploiting (2.6) for $j=1$, we obtain from (5.55):

$$\begin{aligned}&\frac{\partial V_1}{\partial \xi_1}(\xi_1)\left(h_1(w,\xi_1,\theta,d)+b_1 g_1(w,\xi_1,\theta)k_1(\xi_1,z)\right)\\&\le|\xi_1||h_1(w,\xi_1,\theta,d)|-b_1\frac{\xi_1^2 S(\xi_1,z)+4a}{4a}\rho S(\xi_1,z)\xi_1^2\end{aligned}\tag{5.56}$$

Using (2.8) for $j=1$ and the inequalities

$$|\xi_1|\mu_1(\xi_1)(1+|\theta|)\sqrt{R(w)+\Lambda}\le\frac{a}{2\rho}\left(R(w)+\Lambda\right)+\frac{\rho}{2a}\mu_1^2(\xi_1)(1+|\theta|)^2\xi_1^2$$

$$(1+|\theta|)\mu_1(\xi_1)|\xi_1||d|\le\frac{a}{\rho}|d|^2+\frac{\rho}{4a}\mu_1^2(\xi_1)(1+|\theta|)^2\xi_1^2$$

$$(1+|\theta|)\mu_1(\xi_1)|\xi_1||d|\sqrt{R(w)+\Lambda}\le\frac{a}{2\rho}\left(R(w)+\Lambda\right)+\frac{\rho}{2a}\mu_1^2(\xi_1)(1+|\theta|)^2|d|^2\xi_1^2$$

we obtain from (5.56):

$$\begin{aligned}&\frac{\partial V_1}{\partial \xi_1}(\xi_1)\left(h_1(w,\xi_1,\theta,d)+b_1 g_1(w,\xi_1,\theta)k_1(\xi_1,z)\right)\\&\le\mu_1(\xi_1)(1+|\theta|)\xi_1^2+\frac{a}{\rho}\left(|d|^2+R(w)+\Lambda\right)+\frac{3\rho}{4a}\mu_1^2(\xi_1)(1+|\theta|)^2\xi_1^2\\&+\frac{\rho}{2a}\mu_1^2(\xi_1)(1+|\theta|)^2|d|^2\xi_1^2-b_1\frac{\xi_1^2 S(\xi_1,z)+4a}{4a}\rho S(\xi_1,z)\xi_1^2\end{aligned}\tag{5.57}$$

The inequalities $|\theta|\le\rho+(|\theta|-\rho)^+$, $(1+|\theta|)^2\le 2(1+\rho)^2+2\left((|\theta|-\rho)^+\right)^2$ in conjunction with (5.57) give:

$$
\begin{aligned}
&\frac{\partial V_1}{\partial \xi_1}(\xi_1)\big(h_1(w,\xi_1,\theta,d)+b_1 g_1(w,\xi_1,\theta)k_1(\xi_1,z)\big)\\
&\le -Q(\xi_1,z)+\mu_1(\xi_1)(|\theta|-\rho)^+\xi_1^2+\frac{a}{\rho}\left(|d|^2+R(w)+\Lambda\right)\\
&+\xi_1^2\left(\psi(\xi_1)+(1+\rho)^2\frac{3\rho}{2a}\mu_1^2(\xi_1)+\mu_1(\xi_1)(1+\rho)+\frac{\rho^3}{16a^3}\xi_1^2\mu_1^4(\xi_1)\right)\\
&+\xi_1^2\left(\frac{\rho^3}{2a}\mu_1^2(\xi_1)+\frac{\rho}{4a}\xi_1^2\mu_1^2(\xi_1)+\rho\mu_1(\xi_1)\right)\\
&+\left((|\theta|-\rho)^+\right)^2\frac{3\rho}{2a}\xi_1^2\mu_1^2(\xi_1)+\frac{\rho}{a}|d|^2\xi_1^2\mu_1^2(\xi_1)(1+\rho)^2\\
&+\frac{\rho}{a}|d|^2\xi_1^2\mu_1^2(\xi_1)\left((|\theta|-\rho)^+\right)^2-b_1\frac{\xi_1^2 S(\xi_1,z)+4a}{4a}\rho S(\xi_1,z)\xi_1^2
\end{aligned}
\tag{5.58}
$$

It follows from (5.58) and the inequalities

$$
\begin{aligned}
&\mu_1(\xi_1)(|\theta|-\rho)^+\xi_1^2\le\frac{a}{\rho}\left((|\theta|-\rho)^+\right)^2+\frac{\rho}{4a}\mu_1^2(\xi_1)\xi_1^4\\
&\frac{3\rho}{2a}\left((|\theta|-\rho)^+\right)^2\mu_1^2(\xi_1)\xi_1^2\le\frac{a}{\rho}\left((|\theta|-\rho)^+\right)^4+\frac{9\rho^3}{16a^3}\mu_1^4(\xi_1)\xi_1^4\\
&\frac{\rho}{a}(1+\rho)^2\mu_1^2(\xi_1)|d|^2\xi_1^2\le\frac{a}{2\rho}|d|^4+\frac{\rho^3}{2a^3}(1+\rho)^4\mu_1^4(\xi_1)\xi_1^4\\
&\frac{\rho}{a}\mu_1^2(\xi_1)\left((|\theta|-\rho)^+\right)^2|d|^2\xi_1^2\le\frac{a}{2\rho}|d|^4+\frac{\rho^3}{2a^3}\mu_1^4(\xi_1)\left((|\theta|-\rho)^+\right)^4\xi_1^4
\end{aligned}
$$

that the following estimate holds:

$$
\begin{aligned}
&\frac{\partial V_1}{\partial \xi_1}(\xi_1)\big(h_1(w,\xi_1,\theta,d)+b_1 g_1(w,\xi_1,\theta)k_1(\xi_1,z)\big)\\
&\le -Q(\xi_1,z)+\frac{a}{\rho}\left(|d|^2+|d|^4+\left((|\theta|-\rho)^+\right)^2+\left((|\theta|-\rho)^+\right)^4+R(w)+\Lambda\right)\\
&+\xi_1^2\left(\psi(\xi_1)+\rho\mu_1(\xi_1)+(1+\rho)^2\frac{3\rho}{2a}\mu_1^2(\xi_1)+\mu_1(\xi_1)(1+\rho)+\frac{5\rho^3}{8a^3}\mu_1^4(\xi_1)\xi_1^2\right)\\
&+\xi_1^2\left(\frac{\rho^3}{2a}\mu_1^2(\xi_1)+\frac{\rho^3}{2a^3}(1+\rho)^4\mu_1^4(\xi_1)\xi_1^2+\frac{\rho}{4a}\xi_1^2\mu_1^2(\xi_1)+\frac{\rho}{4a}\mu_1^2(\xi_1)\xi_1^2\right)\\
&+\frac{\rho^3}{2a^3}\mu_1^4(\xi_1)\left((|\theta|-\rho)^+\right)^4\xi_1^4-b_1\frac{\xi_1^2 S(\xi_1,z)+4a}{4a}\rho S(\xi_1,z)\xi_1^2
\end{aligned}
$$

The above inequality and the inequality

$$
\frac{\rho^3}{2a^3}\mu_1^4(\xi_1)\left((|\theta|-\rho)^+\right)^4\xi_1^4\le\frac{a}{\rho}\left((|\theta|-\rho)^+\right)^8+\frac{\rho^7}{16a^7}\mu_1^8(\xi_1)\xi_1^8
$$

in conjunction with (5.45), (5.54) gives:

$$\begin{aligned}
&\frac{\partial V_1}{\partial \xi_1}(\xi_1)\left(h_1(w,\xi_1,\theta,d)+b_1 g_1(w,\xi_1,\theta)k_1(\xi_1,z)\right)\\
&\leq -cV_1(\xi_1)-Q\left(\xi_1,z\right)+\xi_1^2 S\left(\xi_1,z\right)-b_1\frac{\xi_1^2 S\left(\xi_1,z\right)+4a}{4a}\rho S\left(\xi_1,z\right)\xi_1^2\\
&+\frac{a}{\rho}\left(\left(\left(|\theta|-\rho\right)^+\right)^2+\left(\left(|\theta|-\rho\right)^+\right)^4+\left(\left(|\theta|-\rho\right)^+\right)^8+|d|^2+|d|^4+R\left(w\right)+\Lambda\right)
\end{aligned} \tag{5.59}$$

The inequality $1=b_1\dfrac{1}{b_1}\leq b_1\left(\dfrac{1}{b_1}-\rho\right)^+ +b_1\rho$ allows us to obtain from (5.59):

$$\begin{aligned}
&\frac{\partial V_1}{\partial \xi_1}(\xi_1)\left(h_1(w,\xi_1,\theta,d)+b_1 g_1(w,\xi_1,\theta)k_1(\xi_1,z)\right)\\
&\leq -cV_1(\xi_1)-Q\left(\xi_1,z\right)+b_1\left(\frac{1}{b_1}-\rho\right)^+ S\left(\xi_1,z\right)\xi_1^2-b_1\frac{\rho}{4a}S^2\left(\xi_1,z\right)\xi_1^4\\
&+\frac{a}{\rho}\left(\left(\left(|\theta|-\rho\right)^+\right)^2+\left(\left(|\theta|-\rho\right)^+\right)^4+\left(\left(|\theta|-\rho\right)^+\right)^8+|d|^2+|d|^4+R\left(w\right)+\Lambda\right)
\end{aligned} \tag{5.60}$$

Completing the squares in (5.60) we obtain (5.44).

<u>Part IV:</u> Let arbitrary $\kappa,\Gamma,\varepsilon,c>0$, $\beta\in(0,1)$ be given. Define $a=\dfrac{\beta\kappa}{1+2^{m-1}}$. By virtue of Part I and Part III there exist functions $k,V\in C^\infty\left(\mathbb{R}^m\times\mathbb{R}\right)$ with $V(0,z)=k(0,z)=0$ for $z\in\mathbb{R}$, $V(\xi,z)>0$ for $\xi\neq 0$, such that:

(a) for each $\tilde{M}>0$ the set $\left\{(\xi,z)\in\mathbb{R}^m\times\left[-\tilde{M},\tilde{M}\right]:V(\xi,z)\leq\tilde{M}\right\}$ is bounded, and

(b) Inequality (5.46) as well as the following inequalities hold:

$$\sum_{i=1}^{m-1}\frac{\partial V}{\partial \xi_i}(\xi,z)\left(h_j(w,\xi_1,...,\xi_j,\theta,d)+b_j g_j(w,\xi_1,...,\xi_j,\theta)\xi_{j+1}\right)$$
$$+\frac{\partial V}{\partial \xi_m}(w,z)\left(h_m(w,\xi,\theta,d)+b_m g_m(w,\xi,\theta)k\left(\xi,z\right)\right)$$
$$+\frac{\partial V}{\partial z}(\xi,z)\Gamma\exp(-z)\left(V(\xi,z)-\varepsilon\right)^{+}$$
$$\leq -2^{1-m}cV(\xi,z)-Q(\xi_1,z)+\frac{2^{m-1}\beta\kappa}{\left(1+2^{m-1}\right)\rho}\left(\min_{i=1,...,m}\left(b_i\right)\left(\left(\frac{1}{\min_{i=1,...,m}\left(b_i\right)}-\rho\right)^{+}\right)^2+|d|^2+|d|^4+R(w)+\Lambda\right)$$
$$+\frac{2^{m-1}\beta\kappa}{\left(1+2^{m-1}\right)\rho}\left(\left(\left(|\theta|+|b|+|b||\theta|-\rho\right)^{+}\right)^2+\left(\left(|\theta|+|b|+|b||\theta|-\rho\right)^{+}\right)^4+\left(\left(|\theta|+|b|+|b||\theta|-\rho\right)^{+}\right)^8\right)$$
$$\text{for all } \left(\xi,w\right)\in\mathbb{R}^m\times\mathbb{R}^l,\ z\in\mathbb{R},\ d\in\mathbb{R}^q,\ \theta\in\Theta,\ b_j\in\left(0,+\infty\right)\cap B_j \text{ for } j=1,...,m \tag{5.61}$$

$$\nabla\Phi(w)f(\xi_1,w,\delta,\theta)\leq-\left(1-\frac{\beta}{1+2^{m-1}}\right)R(w)+\gamma\left(|\delta|\right)+\bar{\gamma}\left(|\theta|\right)$$
$$+Q\left(\xi_1,z\right)+\frac{\beta\kappa}{\left(1+2^{m-1}\right)\rho}\Lambda+\frac{\beta\kappa}{\left(1+2^{m-1}\right)\rho}\left(\left(|\theta|-\rho\right)^{+}\right)^2+\frac{\beta\kappa}{\left(1+2^{m-1}\right)\rho}\left(\left(|\theta|-\rho\right)^{+}\right)^4$$
$$\text{for all } \left(\xi,w\right)\in\mathbb{R}^m\times\mathbb{R}^l,\ \delta\in\mathbb{R}^\alpha,\ \theta\in\Theta,\ z\in\mathbb{R} \tag{5.62}$$

where $\rho:=\kappa+\exp(z)$ and $Q\left(\xi,z\right)$ is given by (5.43).

We next perform the stability analysis of the closed-loop system (2.1), (2.2) with (2.17), (2.18). To this end, we need some facts.

First, the fact that the smooth functions $\Phi,R:\mathbb{R}^l\to\mathbb{R}_+$ are positive definite and radially unbounded, functions and Proposition 2.2 on page 107 in [7] implies the existence of functions $a_1,a_2,a_3,a_4\in K_\infty$ such that the following inequalities hold for all $w\in\mathbb{R}^l$:

$$a_1\left(\Phi\left(w\right)\right)\leq R\left(w\right)\leq a_2\left(\Phi\left(w\right)\right) \tag{5.63}$$

$$a_3\left(|w|\right)\leq\Phi\left(w\right)\leq a_4\left(|w|\right) \tag{5.64}$$

Secondly, we notice the existence of a function $a_5\in K_\infty$ such that the following inequality is valid for all $w\in\mathbb{R}^l$, $\xi\in\mathbb{R}^m$ and $z\in\mathbb{R}$:

$$a_5\left(\Phi\left(w\right)+V\left(\xi,z\right)\right)\leq\left(1-\beta\right)R\left(w\right)+2^{1-m}cV\left(\xi,z\right) \tag{5.65}$$

Indeed, (5.65) holds with $a_5\left(s\right)=\min\left(\left(1-\beta\right)a_1\left(\frac{s}{2}\right),2^{-m}cs\right)$ and notice that $a_5\in K_\infty$. To see why (5.65) is valid notice that since $a_5\left(s_1+s_2\right)\leq a_5\left(2s_1\right)+a_5\left(2s_2\right)$ for all $s_1,s_2\geq 0$ we get:

$$
\begin{aligned}
& a_5\left(\Phi(w)+V(\xi,z)\right) \leq a_5\left(2\Phi(w)\right)+a_5\left(2V(\xi,z)\right) \\
& = \min\left(\left(1-\beta\right)a_1\left(\Phi(w)\right), 2^{1-m}c\Phi(w)\right)+\min\left(\left(1-\beta\right)a_1\left(V(\xi,z)\right), 2^{1-m}cV(\xi,z)\right) \\
& \leq \left(1-\beta\right)a_1\left(\Phi(w)\right)+2^{1-m}cV(\xi,z)
\end{aligned}
$$

and inequality (5.65) is a consequence of the left inequality (5.63).

Thirdly, we notice that there exist functions $\omega \in KL$, $\bar{\zeta} \in K_\infty$ with the following property: for every $T \in (0,+\infty]$, $\tilde{\mu} \geq 0$ and for every absolutely continuous function $Y:[0,T) \to \mathbb{R}_+$ for which the following inequality holds for $t \in [0,T)$ a.e.:

$$
\dot{Y}(t) \leq -a_5\left(Y(t)\right)+\tilde{\mu} \tag{5.66}
$$

the following estimate holds for all $t \in [0,T)$:

$$
Y(t) \leq \omega\left(Y(0),t\right)+\bar{\zeta}\left(\tilde{\mu}\right) \tag{5.67}
$$

The existence of functions $\omega \in KL$, $\zeta \in K_\infty$ follows from Lemma 2.14 on page 82 in [7] and by noticing that (5.66) shows that the implication

$$
Y(t) \geq \tilde{\rho}^{-1}\left(2\tilde{\mu}\right) \Rightarrow \dot{Y}(t) \leq -\tilde{\rho}\left(Y(t)\right)/2
$$

holds for $t \in [0,T)$ a.e..

Fourthly, we notice the existence of a non-decreasing, continuous function $H:\mathbb{R}_+ \to \mathbb{R}_+$ such that the following inequality is valid for all $\xi \in \mathbb{R}^m$ and $z \in \mathbb{R}$:

$$
|\xi| \leq H\left(\max\left(V(\xi,z),|z|\right)\right) \tag{5.68}
$$

Indeed, we define for all $r \geq 0$:

$$
\tilde{H}(r) := \max\left\{|\xi| : V(\xi,z) \leq r \text{ and } |z| \leq r\right\} \tag{5.69}
$$

The fact that for each $r \geq 0$ the set $\left\{ (\xi,z) \in \mathbb{R}^m \times [-r,r] : V(\xi,z) \leq r \right\}$ is bounded guarantees that $\tilde{H}$ is well-defined. It follows from (5.69) that $\tilde{H}$ is non-decreasing. Moreover, it follows from (5.69) that the following inequality is valid for all $\xi \in \mathbb{R}^m$ and $z \in \mathbb{R}$:

$$
|\xi| \leq \tilde{H}\left(\max\left(V(\xi,z),|z|\right)\right) \tag{5.70}
$$

Define the continuous function $H(r) := \int_r^{r+1} \tilde{H}(l)dl$ for $r \geq 0$ (notice that since $\tilde{H}$ is non-decreasing, it is locally bounded and Riemann-integrable on every bounded interval of $\mathbb{R}_+$). We conclude from (5.70) that estimate (5.68) holds.

We also define the function:

$$
U(w,\xi,z) = \Phi(w)+V(\xi,z) \tag{5.71}
$$

Let arbitrary $(w_0,\xi_0,z_0)\in\mathbb{R}^l\times\mathbb{R}^m\times\mathbb{R}$, $d\in L^\infty\left(\mathbb{R}_+;\mathbb{R}^l\right)$, $\theta\in L^\infty\left(\mathbb{R}_+;\Theta\right)$, $\delta\in L^\infty\left(\mathbb{R}_+;\mathbb{R}^\alpha\right)$, $b_j\in L^\infty\left(\mathbb{R}_+;(0,+\infty)\cap B_j\right)$ with $\inf_{t\ge 0}\left(b_j(t)\right)>0$ for $j=1,...,m$ be given. Consider the unique solution $(w(t),\xi(t),z(t))$ of (2.1), (2.2) with (2.17), (2.18) and $(w(0),\xi(0),z(0))=(w_0,\xi_0,z_0)$ for which there exists $t_{\max}\in(0,+\infty]$ such that solution is defined on $[0,t_{\max})$. Moreover, if $t_{\max}<+\infty$ then $\limsup_{t\to t_{\max}^-}\left(|w(t)|+|\xi(t)|+|z(t)|\right)=+\infty$.

By virtue of (2.18) the map $[0,t_{\max})\ni t\mapsto z(t)$ is non-decreasing. Therefore, we get for all $t\in[0,t_{\max})$:

$$z_0\le z(t) \tag{5.72}$$

Define for $t\ge 0$ a.e.:

$$\begin{aligned}
v(t) &:= \gamma\left(|\delta(t)|\right)+\bar{\gamma}\left(|\theta(t)|\right)+\beta\left(|d(t)|^2+|d(t)|^4\right)\\
&+\beta\left(\left(|\theta(t)|+|b(t)|+|\theta(t)||b(t)|-\kappa-\exp\left(z(t)\right)\right)^+\right)^2\\
&+\beta\left(\left(|\theta(t)|+|b(t)|+|\theta(t)||b(t)|-\kappa-\exp\left(z(t)\right)\right)^+\right)^4\\
&+\beta\left(\left(|\theta(t)|+|b(t)|+|\theta(t)||b(t)|-\kappa-\exp\left(z(t)\right)\right)^+\right)^8\\
&+\beta\left(\min_{i=1,...,m}\left(b_i(t)\right)\left(\left(\frac{1}{\min_{i=1,...,m}\left(b_i(t)\right)}-\kappa-\exp\left(z(t)\right)\right)^+\right)^2+\Lambda\right)
\end{aligned} \tag{5.73}$$

We notice that since $\gamma,\bar{\gamma}\in C^0\left(\mathbb{R}_+;\mathbb{R}_+\right)$ are non-decreasing functions and since $d\in L^\infty\left(\mathbb{R}_+;\mathbb{R}^l\right)$, $\theta\in L^\infty\left(\mathbb{R}_+;\Theta\right)$, $\delta\in L^\infty\left(\mathbb{R}_+;\mathbb{R}^\alpha\right)$, $b_j\in L^\infty\left(\mathbb{R}_+;(0,+\infty)\cap B_j\right)$ with $\inf_{t\ge 0}\left(b_j(t)\right)>0$ for $j=1,...,m$, definition (5.73) implies that $v\in L^\infty\left(\mathbb{R}_+;\mathbb{R}\right)$ with

$$\begin{aligned}
\|v\|_\infty &\le G\left(\|d\|_\infty,\|\delta\|_\infty,\|\theta\|_\infty,\|b\|_\infty,\min_{i=1,...,m}\left(\inf_{s\ge 0}\left(b_i(s)\right)\right)\right)\\
&:= \gamma\left(\|\delta\|_\infty\right)+\bar{\gamma}\left(\|\theta\|_\infty\right)+\beta\left(\|\theta\|_\infty+\|b\|_\infty+\|\theta\|_\infty\|b\|_\infty\right)^2\\
&+\beta\left(\|\theta\|_\infty+\|b\|_\infty+\|\theta\|_\infty\|b\|_\infty\right)^4+\beta\left(\|\theta\|_\infty+\|b\|_\infty+\|\theta\|_\infty\|b\|_\infty\right)^8\\
&+\beta\left(\|d\|_\infty^2+\|d\|_\infty^4+\frac{1}{\min_{i=1,...,m}\left(\inf_{s\ge 0}\left(b_i(s)\right)\right)}+\Lambda\right)
\end{aligned} \tag{5.74}$$

Inequalities (5.61), (5.62) and definitions (5.71), (5.73) imply that the following differential inequality holds for $t\in[0,t_{\max})$ a.e.:

$$\frac{d}{dt}\Big(U\big(w(t),\xi(t),z(t)\big)\Big) \le -\big(1-\beta\big)R(w(t)) - 2^{1-m}cV(\xi(t),z(t)) + v(t) \tag{5.75}$$

Combining (5.65), (5.66), (5.67), (5.71) and (5.74), (5.75) we obtain for all $t \in \left[0, t_{\max}\right)$:

$$U\big(w(t),\xi(t),z(t)\big) \le \omega\Big(U\big(w_0,\xi_0,z_0\big),t\Big) + \bar{\zeta}\left(G\left(\|d\|_\infty, \|\delta\|_\infty, \|\theta\|_\infty, \|b\|_\infty, \min_{i=1,\ldots,m}\left(\inf_{s\ge 0}\big(b_i(s)\big)\right)\right)\right) \tag{5.76}$$

Inequalities (5.61), (5.76), the right inequality (5.63) and definition (5.71) imply that the following differential inequality holds for $t \in \left[0, t_{\max}\right)$ a.e.:

$$\begin{aligned} \frac{d}{dt}\Big(V\big(\xi(t),z(t)\big)\Big) &\le -2^{1-m}cV\big(\xi(t),z(t)\big) \\ &+ \frac{Z\left(w_0,\xi_0,z_0,\|d\|_\infty,\|\delta\|_\infty,\|\theta\|_\infty,\|b\|_\infty,\min_{i=1,\ldots,m}\left(\inf_{s\ge 0}\big(b_i(s)\big)\right)\right)}{1+\exp\big(z(t)\big)} \end{aligned} \tag{5.77}$$

where

$$\begin{aligned} &Z\left(w_0,\xi_0,z_0,\|d\|_\infty,\|\delta\|_\infty,\|\theta\|_\infty,\|b\|_\infty,\min_{i=1,\ldots,m}\left(\inf_{s\ge 0}\big(b_i(s)\big)\right)\right) \\ &:= \frac{2^{m-1}\beta\kappa}{\left(1+2^{m-1}\right)\min\left(1,\kappa\right)} a_2\left(\omega\Big(U\big(w_0,\xi_0,z_0\big),0\Big) + \bar{\zeta}\left(G\left(\|d\|_\infty,\|\delta\|_\infty,\|\theta\|_\infty,\|b\|_\infty,\min_{i=1,\ldots,m}\left(\inf_{s\ge 0}\big(b_i(s)\big)\right)\right)\right)\right) \\ &+ \frac{2^{m-1}\beta\kappa}{\left(1+2^{m-1}\right)\min\left(1,\kappa\right)}\left(\left(\|\theta\|_\infty + \|b\|_\infty + \|\theta\|_\infty\|b\|_\infty\right)^2 + \left(\|\theta\|_\infty + \|b\|_\infty + \|\theta\|_\infty\|b\|_\infty\right)^4\right) \\ &+ \frac{2^{m-1}\beta\kappa}{\left(1+2^{m-1}\right)\min\left(1,\kappa\right)}\left(\|d\|_\infty^2 + \|d\|_\infty^4 + \Lambda + \frac{1}{\min_{i=1,\ldots,m}\left(\inf_{s\ge 0}\big(b_i(s)\big)\right)} + \left(\|\theta\|_\infty + \|b\|_\infty + \|\theta\|_\infty\|b\|_\infty\right)^8\right) \end{aligned} \tag{5.78}$$

Lemma 4.1 on page 63 in [10], (2.1) and differential inequality (5.77) imply the following estimate for all $t \in \left[0, t_{\max}\right)$:

$$z_0 \le z(t) \le \ln\left(\tilde{Z}\left(w_0,\xi_0,z_0,\|d\|_\infty,\|\delta\|_\infty,\|\theta\|_\infty,\|b\|_\infty,\min_{i=1,\ldots,m}\left(\inf_{s\ge 0}\big(b_i(s)\big)\right)\right)\right) \tag{5.79}$$

where

$$\begin{aligned} &\tilde{Z}\left(w_0,\xi_0,z_0,\|d\|_\infty,\|\delta\|_\infty,\|\theta\|_\infty,\|b\|_\infty,\min_{i=1,\ldots,m}\left(\inf_{s\ge 0}\big(b_i(s)\big)\right)\right) \\ &:= \left(\frac{1}{\varepsilon} + \frac{2^{m-1}\Gamma}{c}\right)\frac{2^{m-1}\mu}{c} Z\left(w_0,\xi_0,z_0,\|d\|_\infty,\|\delta\|_\infty,\|\theta\|_\infty,\|b\|_\infty,\min_{i=1,\ldots,m}\left(\inf_{s\ge 0}\big(b_i(s)\big)\right)\right) \\ &+ \exp(z_0) + \frac{2^{m-1}\Gamma}{c}V\big(\xi_0,z_0\big) \end{aligned} \tag{5.80}$$

It follows from (5.68), the left inequality (5.64), definition (5.71) and estimates (5.79), (5.76) that the following estimates hold for all $t \in [0, t_{\max})$:

$$|w(t)| \le a_3^{-1}\left(\bar{Z}\right), \quad |\xi(t)| \le H\left(\bar{Z}+K\right), \quad |z(t)| \le K \tag{5.81}$$

where

$$\begin{aligned} \bar{Z} &:= \omega\big(U\left(w_0,\xi_0,z_0\right),0\big)+\bar{\zeta}\left(G\left(\|d\|_\infty,\|\delta\|_\infty,\|\theta\|_\infty,\|b\|_\infty,\min_{i=1,\dots,m}\left(\inf_{s\ge 0}\left(b_i(s)\right)\right)\right)\right) \\ K &:= |z_0|+\left|\ln\left(\tilde{Z}\left(w_0,\xi_0,z_0,\|d\|_\infty,\|\delta\|_\infty,\|\theta\|_\infty,\|b\|_\infty,\min_{i=1,\dots,m}\left(\inf_{s\ge 0}\left(b_i(s)\right)\right)\right)\right)\right| \end{aligned} \tag{5.82}$$

Estimates (5.81) imply that $t_{\max} = +\infty$. Moreover, estimates (5.81), (5.79) and the fact that $z(t)$ is non-decreasing, imply estimates (2.20), (2.21) for appropriate $\Xi \in C^0\left(\mathbb{R}_+^7 \times (0,+\infty)\right)$.

We next show estimate (2.19). Estimate (5.81) and the fact that $z(t)$ is non-decreasing guarantees that the function $z(t)$ has a finite limit as $t \to +\infty$. This implies that the function $\exp(z(t))$ has a finite limit as $t \to +\infty$. Moreover, the facts that $d \in L^\infty\left(\mathbb{R}_+;\mathbb{R}^l\right)$, $\theta \in L^\infty\left(\mathbb{R}_+;\Theta\right)$, $b_j \in L^\infty\left(\mathbb{R}_+;(0,+\infty)\cap B_j\right)$ with $\inf_{t\ge 0}\left(b_j(t)\right)>0$ for $j=1,\dots,m$ and (5.81), (2.1), (2.17), (2.18) imply that $\frac{d}{dt}\left(V(\xi(t),z(t))\right)$ is of class $L^\infty\left(\mathbb{R}_+\right)$. It follows that the function $\frac{d}{dt}\left(\exp(z(t))\right)=\Gamma\left(V(\xi(t),z(t))-\varepsilon\right)^+$ is uniformly continuous, in addition to $\int_0^{+\infty}\frac{d}{dt}\left(\exp(z(t))\right)dt=\lim_{t\to+\infty}\left(\exp(z(t))\right)-\exp(z_0)<+\infty$. From Barbălat's Lemma (see [8, 14]), we have:

$$\lim_{t\to+\infty}\left(\frac{d}{dt}\left(\exp(z(t))\right)\right)=\lim_{t\to+\infty}\left(\Gamma\left(V(\xi(t),z(t))-\varepsilon\right)^+\right)=0 \tag{5.83}$$

Therefore, the right estimate (2.19) holds. The left estimate (2.19) is a consequence of the right estimate (2.19) and inequality (5.46).

Finally, we consider the case $\gamma(0)=\Lambda=0$, $\bar{\gamma}(s)\equiv 0$, $\lim_{t\to+\infty}\left(d(t)\right)=0$, $\lim_{t\to+\infty}\left(\delta(t)\right)=0$, $\liminf_{t\to+\infty}\left(b_i(t)\right)\ge\frac{1}{\kappa}$ for $i=1,\dots,m$ and $\kappa+\exp\left(\lim_{s\to+\infty}\left(z(s)\right)\right)\ge\limsup_{t\to+\infty}\left(|\theta(t)|+|b(t)|+|\theta(t)||b(t)|\right)$. In this case, definition (5.73) shows that $\lim_{t\to+\infty}\left(v(t)\right)=0$. Combining (5.65), (5.66), (5.67), (5.71) and (5.75) we obtain:

$$\limsup_{t\to+\infty}\big(U\left(w(t),\xi(t),z(t)\right)\big)\le\bar{\zeta}\left(\limsup_{t\to+\infty}\left(v(t)\right)\right) \tag{5.84}$$

It follows from (5.84) and $\lim_{t\to+\infty}\left(v(t)\right)=0$ that $\lim_{t\to+\infty}\big(U\left(w(t),\xi(t),z(t)\right)\big)=0$. Clearly, the equation $\lim_{t\to+\infty}\big(U\left(w(t),\xi(t),z(t)\right)\big)=0$ in conjunction with definition (5.71) gives $\lim_{t\to+\infty}\left(\Phi\left(w(t)\right)\right)=\lim_{t\to+\infty}\left(V\left(\xi(t),z(t)\right)\right)=0$. Inequalities (5.64) and $\lim_{t\to+\infty}\left(\Phi\left(w(t)\right)\right)=0$ show that $\lim_{t\to+\infty}\left(|w(t)|\right)=0$. We next show by contradiction that $\lim_{t\to+\infty}\left(|\xi(t)|\right)=0$. Suppose that $\limsup_{t\to+\infty}\left(|\xi(t)|\right)>0$. In this case there exists $\bar{\varepsilon}>0$ and a sequence $\left\{t_k\in\mathbb{R}_+ : k=1,...\right\}$ with $\lim_{k\to+\infty}\left(t_k\right)=+\infty$ such that $\left|\xi\left(t_k\right)\right|\geq\bar{\varepsilon}$ for $k=1,2,...$. By virtue of (5.81) the sequence $\left\{\left(\xi\left(t_k\right),z\left(t_k\right)\right)\in\mathbb{R}^m\times\mathbb{R} : k=1,...\right\}$ is bounded. Consequently, without loss of generality (by restricting our attention to a subsequence of $\left\{t_k\in\mathbb{R}_+ : k=1,...\right\}$ if necessary), we may assume the existence of $\xi^*\in\mathbb{R}^m$ and $z^*\in\mathbb{R}$ such that $\lim_{k\to+\infty}\left(\xi\left(t_k\right)\right)=\xi^*$ and $\lim_{k\to+\infty}\left(z\left(t_k\right)\right)=z^*$. Clearly, $\left|\xi^*\right|\geq\bar{\varepsilon}$ and since $\lim_{t\to+\infty}\left(V\left(\xi(t),z(t)\right)\right)=0$ we get (by continuity of $V$) $V\left(\xi^*,z^*\right)=0$. However, $V(0,z)=0$ for $z\in\mathbb{R}$ and $V(\xi,z)>0$ for $\xi\neq 0$, we get $\xi^*=0$; a contradiction with $\left|\xi^*\right|\geq\bar{\varepsilon}>0$.

The proof is complete. ◁

We finally provide the proof of Theorem 2.

**Proof of Theorem 2:** The proof is divided in three parts.

Part I: We show that for every $a,\kappa>0$ there exist functions $k_1\in C^\infty\left(\mathbb{R}^n\times\mathbb{R}\times\mathbb{R}\right)$ with $k_1(0,0,z)=0$ for $z\in\mathbb{R}$, $V_1\in C^\infty\left(\mathbb{R}^n\times\mathbb{R}\times\mathbb{R}\right)$ being positive definite and radially unbounded and a constant $c>0$ such that the following inequalities hold for all $\left(\zeta,\xi,w\right)\in\mathbb{R}^n\times\mathbb{R}^m\times\mathbb{R}^l$, $\theta\in\Theta$, $d\in\mathbb{R}^q$, $b_1\in\left(0,+\infty\right)\cap B_1$ and $z\in\mathbb{R}$:

$$\begin{aligned}&\frac{\partial V_1}{\partial\zeta}(\zeta,\xi_1)\left(A\zeta+\bar{b}\xi_1\right)+\frac{\partial V_1}{\partial\xi_1}(\zeta,\xi_1)\left(h_1(w,\zeta,\xi_1,\theta,d)+b_1 g_1(w,\zeta,\xi_1,\theta)k_1(\zeta,\xi_1,z)\right)\\&\leq -cV_1(\zeta,\xi_1)-Q(\zeta,\xi_1)+\frac{a}{\rho}b_1\left(\left(\frac{1}{b_1}-\rho\right)^+\right)^2\\&+\frac{a}{\rho}\left(\left(\left(|\theta|-\rho\right)^+\right)^2+\left(\left(|\theta|-\rho\right)^+\right)^4+\left(\left(|\theta|-\rho\right)^+\right)^8+|d|^2+|d|^4+R(w)+\Lambda\right)\end{aligned}\tag{5.85}$$

$$|\zeta|^2+\xi_1^2\leq V_1(\zeta,\xi_1)\tag{5.86}$$

where $\rho=\kappa+\exp\left(z\right)$.

Part II: Applying Lemma exactly $m-1$ times and exploiting Assumption (B), we construct for every tuple of constants $\kappa,\Gamma,\varepsilon,a>0$ functions $k,V\in C^\infty\left(\mathbb{R}^n\times\mathbb{R}^m\times\mathbb{R}\right)$ with $V(0,0,z)=k(0,0,z)=0$ for $z\in\mathbb{R}$, $V(\zeta,\xi,z)>0$ for $\left(\zeta,\xi\right)\neq\left(0,0\right)$ and a constant $c>0$ such that:

(i) for each $\tilde{M} > 0$ the set $\left\{ (\zeta, \xi, z) \in \mathbb{R}^n \times \mathbb{R}^m \times \left[-\tilde{M}, \tilde{M}\right] : V(\zeta, \xi, z) \le \tilde{M} \right\}$ is bounded, and

(ii) the following inequalities hold:

$$|\zeta|^2 + \xi_1^2 \le V(\zeta, \xi, z), \text{ for all } (\zeta, \xi, z) \in \mathbb{R}^n \times \mathbb{R}^m \times \mathbb{R} \tag{5.87}$$

$$\begin{aligned}
&\frac{\partial V}{\partial \zeta}(\zeta, \xi, z)\left(A\zeta + \bar{b}\xi_1\right) + \sum_{i=1}^{m-1} \frac{\partial V}{\partial \xi_i}(\zeta, \xi, z)\left(h_j(w, \zeta, \xi_1, ..., \xi_j, \theta, d) + b_j g_j(w, \zeta, \xi_1, ..., \xi_j, \theta)\xi_{j+1}\right) \\
&+ \frac{\partial V}{\partial \xi_m}(\zeta, \xi, z)\left(h_m(w, \zeta, \xi, \theta, d) + b_m g_m(w, \zeta, \xi, \theta)k(\zeta, \xi, z)\right) \\
&+ \frac{\partial V}{\partial z}(\zeta, \xi, z)\Gamma \exp(-z)\left(V(\zeta, \xi, z) - \varepsilon\right)^+ \\
&\le -2^{1-m} c V(\zeta, \xi, z) - Q(\zeta, \xi_1) + \frac{2^{m-1}a}{\rho}\left( \min_{i=1,...,m}(b_i)\left(\left(\frac{1}{\min_{i=1,...,m}(b_i)} - \rho\right)^+\right)^2 + |d|^2 + |d|^4 + R(w) \right) \\
&+ \frac{2^{m-1}a}{\rho}\left( \left(\left(|\theta| + |b| + |\theta||b| - \rho\right)^+\right)^2 + \left(\left(|\theta| + |b| + |\theta||b| - \rho\right)^+\right)^4 + \left(\left(|\theta| + |b| + |\theta||b| - \rho\right)^+\right)^6 + \Lambda \right)
\end{aligned}$$

for all $(\zeta, \xi, w) \in \mathbb{R}^n \times \mathbb{R}^m \times \mathbb{R}^l$, $z \in \mathbb{R}$, $d \in \mathbb{R}^q$, $\theta \in \Theta$, $b_j \in (0, +\infty) \cap B_j$ for $j = 1, ..., m$ (5.88)

More specifically, applying Lemma exactly $m-1$ times, we construct functions $k_j \in C^\infty\left(\mathbb{R}^n \times \mathbb{R}^j \times \mathbb{R}\right)$ for $j = 2, ..., m$ with $k_j(0, 0, z) = 0$ for $z \in \mathbb{R}$, such that

$$V(\xi, z) = V_1(\zeta, \xi_1) + \frac{1}{2}\sum_{j=2}^{m}\left(\xi_j - k_{j-1}(\zeta, \xi_1, ..., \xi_{j-1}, z)\right)^2$$

Notice that the above equation and (5.86) guarantees the validity of (5.87).

<u>Part III:</u> We perform the stability analysis of the closed-loop system (2.3), (2.4), (2.5) with (2.24), (2.25) for the function $k \in C^\infty\left(\mathbb{R}^n \times \mathbb{R}^m \times \mathbb{R}\right)$ with $k(0, 0, z) = 0$ for $z \in \mathbb{R}$, that is designed in Part II for arbitrary constants $\kappa, \Gamma, \varepsilon > 0$ and $a = \frac{\kappa}{2\left(1 + 2^{m-1}\right)}$.

<u>Part I:</u> Let arbitrary $a, c, \kappa > 0$ be given.

Since the pair of matrices $\left(A, \bar{b}\right)$ is stabilizable, there exists a symmetric, positive definite matrix $P \in \mathbb{R}^{n \times n}$, a vector $\tilde{k} \in \mathbb{R}^n$ and a constant $r > 0$ such that the following inequality holds for all $\zeta \in \mathbb{R}^n$:

$$\zeta' P\left(A + \bar{b}\tilde{k}'\right)\zeta + \zeta'\left(A + \bar{b}\tilde{k}'\right)' P\zeta \le -r\zeta' P\zeta \tag{5.89}$$

We define the function:

$$\bar{V}(\zeta,\xi_1) := \zeta' P\zeta + \frac{1}{2}\left(\xi_1 - \tilde{k}'\zeta\right)^2 \tag{5.90}$$

Since $P \in \mathbb{R}^{n\times n}$ is a symmetric, positive definite matrix and since the function $\bar{V}(\zeta,\xi_1)$ defined by (5.90) is a positive definite quadratic form on $\mathbb{R}^n \times \mathbb{R}$, it follows that there exists a constant $\bar{c} > 0$ such that the following inequality holds for all $(\zeta,\xi_1) \in \mathbb{R}^n \times \mathbb{R}$:

$$\bar{V}(\zeta,\xi_1) \geq \bar{c}\left|(\zeta,\xi_1)\right|^2 \tag{5.91}$$

Since $Q:\mathbb{R}^n \times \mathbb{R} \to \mathbb{R}_+$ is a smooth function with $Q(0,0)=0$, there exists a smooth function $\omega:\mathbb{R}^n \times \mathbb{R} \to \mathbb{R}_+$ such that the following inequality holds for all $(\zeta,\xi_1) \in \mathbb{R}^n \times \mathbb{R}$:

$$Q(\zeta,\xi_1) \leq \omega(\zeta,\xi_1)\left|(\zeta,\xi_1)\right|^2 \tag{5.92}$$

Exploiting (5.91) and (5.92) we get:

$$Q(\zeta,\xi_1) \leq \frac{\omega(\zeta,\xi_1)}{\bar{c}}\bar{V}(\zeta,\xi_1) \tag{5.93}$$

Define for all $s \geq 0$:

$$\tilde{\omega}(s) := \max\left\{\omega(\zeta,\xi_1) : \bar{V}(\zeta,\xi_1) \leq s\right\} \tag{5.94}$$

Clearly, the function $\tilde{\omega}(s)$ is a non-negative, non-decreasing function that satisfies $\omega(\zeta,\xi_1) \leq \tilde{\omega}\left(\bar{V}(\zeta,\xi_1)\right)$ for all $(\zeta,\xi_1) \in \mathbb{R}^n \times \mathbb{R}$. Consequently, we obtain from (5.93) for all $(\zeta,\xi_1) \in \mathbb{R}^n \times \mathbb{R}$:

$$Q(\zeta,\xi_1) \leq \frac{1}{\bar{c}}\tilde{\omega}\left(\bar{V}(\zeta,\xi_1)\right)\bar{V}(\zeta,\xi_1) \tag{5.95}$$

We extend the domain of $\tilde{\omega}$ by defining $\tilde{\omega}(s)=0$ for $s<0$ and $\tilde{\omega}$ still remains a non-negative, non-decreasing function. Let $h:\mathbb{R} \to \mathbb{R}_+$ be a smooth function with $\int_0^1 h(y)dy \geq 1$, $h(y)=0$ for all $y \leq 0$ and $y \geq 1$. We define the following smooth function for all $s \geq 0$:

$$\bar{\omega}(s) := \int_{-\infty}^{+\infty} \tilde{\omega}(y)h(y-s)dy \tag{5.96}$$

Notice that definition (5.96) and the fact that $h(y)=0$ for all $y \leq 0$ and $y \geq 1$ implies that $\bar{\omega}(s) = \int_0^1 \tilde{\omega}(s+y)h(y)dy$ for all $s \geq 0$. The fact that $\tilde{\omega}$ is a non-negative, non-decreasing function and the fact that $h$ is non-negative implies that $\bar{\omega}:\mathbb{R}_+ \to \mathbb{R}$ is a non-decreasing function that satisfies $\bar{\omega}(s) \geq \tilde{\omega}(s)$ for all $s \geq 0$. Consequently, we obtain from (5.95) for all $(\zeta,\xi_1) \in \mathbb{R}^n \times \mathbb{R}$:

$$Q(\zeta,\xi_1) \le \frac{1}{\bar{c}}\bar{\omega}\left(\bar{V}(\zeta,\xi_1)\right)\bar{V}(\zeta,\xi_1) \tag{5.97}$$

We define for all $s \ge 0$:

$$\beta(s) := \frac{1}{\bar{c}}s + \frac{4}{\bar{c}r}\int_0^s \bar{\omega}(y)\,dy \tag{5.98}$$

Notice that $\beta \in K_\infty \cap C^\infty(\mathbb{R}_+)$. We define the function:

$$V_1(\zeta,\xi_1) := \beta\left(\bar{V}(\zeta,\xi_1)\right) \tag{5.99}$$

Definitions (5.98), (5.99) and inequality (5.91) imply that inequality (5.86) is valid for all $(\zeta,\xi_1) \in \mathbb{R}^n \times \mathbb{R}$. We also notice that since definition (5.98) gives $\frac{d\beta}{ds}(s) = \frac{1}{\bar{c}} + \frac{4}{\bar{c}r}\bar{\omega}(s)$ for all $s \ge 0$, the fact that $\bar{\omega}$ is non-decreasing implies for all $s \ge 0$:

$$\beta(s) \le \frac{1}{\bar{c}}s + \frac{4}{\bar{c}r}\bar{\omega}(s)s = \frac{d\beta}{ds}(s)s \tag{5.100}$$

Moreover, since $\frac{d\beta}{ds}(s) \ge \frac{4}{\bar{c}r}\bar{\omega}(s)$ we get from (5.97) for all $(\zeta,\xi_1) \in \mathbb{R}^n \times \mathbb{R}$:

$$Q(\zeta,\xi_1) \le \frac{r}{4}\frac{d\beta}{ds}\left(\bar{V}(\zeta,\xi_1)\right)\bar{V}(\zeta,\xi_1) \tag{5.101}$$

We define for all $(\zeta,\xi_1,z) \in \mathbb{R}^n \times \mathbb{R} \times \mathbb{R}$:

$$k_1(\zeta,\xi_1,z) := -\bar{c}\rho S(\zeta,\xi_1,z)\frac{S(\zeta,\xi_1,z)\left(\xi_1 - \tilde{k}'\zeta\right)^2 + 4a}{4a\eta_1(\zeta,\xi_1)}\left(\xi_1 - \tilde{k}'\zeta\right) \tag{5.102}$$

where $\rho := \kappa + \exp(z)$, $\eta_1(\zeta,\xi_1)$ is the positive smooth function involved in (2.12) for $j=1$ and $S(\zeta,\xi_1,z)$ is any smooth function that satisfies the following inequality for all $(\zeta,\xi_1,z) \in \mathbb{R}^n \times \mathbb{R} \times \mathbb{R}$:

$$
\begin{aligned}
S(\zeta,\xi_1,z) \geq & \frac{d\beta}{ds}\left(\bar{V}(\zeta,\xi_1)\right)\left(\mu_1(\zeta,\xi_1)(1+\rho)+\frac{3r}{8}+\left|\tilde{k}'\bar{b}\right|\right) \\
& +\frac{1}{r\bar{c}}\frac{d\beta}{ds}\left(\bar{V}(\zeta,\xi_1)\right)\left(2\left|P\bar{b}\right|+\left|\tilde{k}'\left(A+\bar{b}\tilde{k}'\right)\right|+\mu_1(\zeta,\xi_1)(1+\rho)\left(1+\left|\tilde{k}\right|\right)\right)^2 \\
& +\frac{\rho}{2a}\left(\frac{1}{2}\left|(\zeta,\xi_1)\right|^2+3(1+\rho)^2\right)\left(\frac{d\beta}{ds}\left(\bar{V}(\zeta,\xi_1)\right)\right)^2\mu_1^2(\zeta,\xi_1) \qquad (5.103) \\
& +\frac{\rho^3}{2a^3}\left(\frac{9}{2}+(1+\rho)^4\right)\left(\frac{d\beta}{ds}\left(\bar{V}(\zeta,\xi_1)\right)\right)^4\mu_1^4(\zeta,\xi_1)\left(\xi_1-\tilde{k}'\zeta\right)^2 \\
& +\frac{\rho^7}{16a^7}\left(\frac{d\beta}{ds}\left(\bar{V}(\zeta,\xi_1)\right)\right)^8\mu_1^8(\zeta,\xi_1)\left(\xi_1-\tilde{k}'\zeta\right)^6
\end{aligned}
$$

We next use the following notation:

$$
\begin{aligned}
\dot{V}_1\Big|_{\xi_2=k_1(\zeta,\xi_1,z)} = & \frac{\partial V_1}{\partial\zeta}(\zeta,\xi_1)\left(A\zeta+\bar{b}\xi_1\right) \\
& +\frac{\partial V_1}{\partial\xi_1}(\zeta,\xi_1)\left(h_1(w,\zeta,\xi_1,\theta,d)+b_1g_1(w,\zeta,\xi_1,\theta)k_1(\zeta,\xi_1,z)\right) \qquad (5.104)
\end{aligned}
$$

Combining (5.104), (5.90) and (5.99) we get for all $(\zeta,\xi,w)\in\mathbb{R}^n\times\mathbb{R}^m\times\mathbb{R}^l$, $\theta\in\Theta$, $d\in\mathbb{R}^q$, $b_1\in(0,+\infty)\cap B_1$ and $z\in\mathbb{R}$:

$$
\begin{aligned}
\dot{V}_1\Big|_{\xi_2=k_1(\zeta,\xi_1,z)} = & \frac{d\beta}{ds}\left(\bar{V}(\zeta,\xi_1)\right)\left(\zeta'P\left(A\zeta+\bar{b}\xi_1\right)+\left(A\zeta+\bar{b}\xi_1\right)'P\zeta\right) \\
& +\frac{d\beta}{ds}\left(\bar{V}(\zeta,\xi_1)\right)\left(\xi_1-\tilde{k}'\zeta\right)\left(h_1(w,\zeta,\xi_1,\theta,d)+b_1g_1(w,\zeta,\xi_1,\theta)k_1(\zeta,\xi_1,z)-\tilde{k}'A\zeta-\tilde{k}'\bar{b}\xi_1\right) \qquad (5.105)
\end{aligned}
$$

Since $\beta\in K_\infty\cap C^\infty(\mathbb{R}_+)$ it follows that $\frac{d\beta}{ds}(s)\geq 0$ for all $s\geq 0$. Exploiting the fact that $\frac{d\beta}{ds}(s)\geq 0$ for all $s\geq 0$ and (5.89), (5.102) we get from (5.105):

$$
\begin{aligned}
\dot{V}_1\Big|_{\xi_2=k_1(\zeta,\xi_1,z)} \leq & -r\frac{d\beta}{ds}\left(\bar{V}(\zeta,\xi_1)\right)\zeta'P\zeta+2\frac{d\beta}{ds}\left(\bar{V}(\zeta,\xi_1)\right)\left|\xi_1-\tilde{k}'\zeta\right|\left|P\bar{b}\right|\left|\zeta\right| \\
& +\frac{d\beta}{ds}\left(\bar{V}(\zeta,\xi_1)\right)\left|\xi_1-\tilde{k}'\zeta\right|\left|h_1(w,\zeta,\xi_1,\theta,d)-\tilde{k}'A\zeta-\tilde{k}'\bar{b}\xi_1\right| \qquad (5.106) \\
& -\frac{d\beta}{ds}\left(\bar{V}(\zeta,\xi_1)\right)b_1g_1(w,\zeta,\xi_1,\theta)\bar{c}\rho S(\zeta,\xi_1,z)\frac{S(\zeta,\xi_1,z)\left(\xi_1-\tilde{k}'\zeta\right)^2+4a}{4a\eta_1(\zeta,\xi_1)}\left(\xi_1-\tilde{k}'\zeta\right)^2
\end{aligned}
$$

Combining (2.12) for $j=1$, (5.101), (5.106) we get:

$$\begin{aligned}\dot{V}_1\Big|_{\xi_2=k_1(\zeta,\xi_1,z)} &\leq -r\frac{d\beta}{ds}\left(\bar{V}(\zeta,\xi_1)\right)\zeta' P\zeta - Q(\zeta,\xi_1)\\
&+\frac{r}{4}\frac{d\beta}{ds}\left(\bar{V}(\zeta,\xi_1)\right)\bar{V}(\zeta,\xi_1)+2\frac{d\beta}{ds}\left(\bar{V}(\zeta,\xi_1)\right)\left|\xi_1-\tilde{k}'\zeta\right|\left|P\bar{b}\right||\zeta|\\
&+\frac{d\beta}{ds}\left(\bar{V}(\zeta,\xi_1)\right)\left|\xi_1-\tilde{k}'\zeta\right|\left|h_1(w,\zeta,\xi_1,\theta,d)-\tilde{k}'\left(A+\bar{b}\tilde{k}'\right)\zeta-\tilde{k}'\bar{b}\left(\xi_1-\tilde{k}'\zeta\right)\right|\\
&-b_1\bar{c}\rho S(\zeta,\xi_1,z)\frac{d\beta}{ds}\left(\bar{V}(\zeta,\xi_1)\right)\frac{S(\zeta,\xi_1,z)\left(\xi_1-\tilde{k}'\zeta\right)^2+4a}{4a}\left(\xi_1-\tilde{k}'\zeta\right)^2\end{aligned} \tag{5.107}$$

Using the triangle inequality, the fact that $\frac{d\beta}{ds}(s)\geq\frac{1}{\bar{c}}$ for all $s\geq 0$ (recall (5.98)) and (5.90), we get from (5.107):

$$\begin{aligned}\dot{V}_1\Big|_{\xi_2=k_1(\zeta,\xi_1,z)} &\leq -\frac{3r}{4}\frac{d\beta}{ds}\left(\bar{V}(\zeta,\xi_1)\right)\zeta' P\zeta - Q(\zeta,\xi_1)\\
&+\frac{d\beta}{ds}\left(\bar{V}(\zeta,\xi_1)\right)\left(2\left|P\bar{b}\right|+\left|\tilde{k}'\left(A+\bar{b}\tilde{k}'\right)\right|\right)\left|\xi_1-\tilde{k}'\zeta\right||\zeta|\\
&+\frac{d\beta}{ds}\left(\bar{V}(\zeta,\xi_1)\right)\left|\xi_1-\tilde{k}'\zeta\right|\left|h_1(w,\zeta,\xi_1,\theta,d)\right|\\
&+\left(\frac{r}{8}+\left|\tilde{k}'\bar{b}\right|\right)\frac{d\beta}{ds}\left(\bar{V}(\zeta,\xi_1)\right)\left(\xi_1-\tilde{k}'\zeta\right)^2\\
&-b_1\rho S(\zeta,\xi_1,z)\frac{S(\zeta,\xi_1,z)\left(\xi_1-\tilde{k}'\zeta\right)^2+4a}{4a}\left(\xi_1-\tilde{k}'\zeta\right)^2\end{aligned} \tag{5.108}$$

We next notice that inequality (2.14) for $j=1$ in conjunction with the inequalities

$$\begin{aligned}&\frac{d\beta}{ds}\left(\bar{V}(\zeta,\xi_1)\right)\mu_1(\zeta,\xi_1)\left(1+|\theta|\right)\left|\xi_1-\tilde{k}'\zeta\right|\sqrt{R(w)+\Lambda}\\
&\leq\frac{a}{2\rho}\left(R(w)+\Lambda\right)+\frac{\rho}{2a}\left(\frac{d\beta}{ds}\left(\bar{V}(\zeta,\xi_1)\right)\right)^2\mu_1^2(\zeta,\xi_1)\left(1+|\theta|\right)^2\left(\xi_1-\tilde{k}'\zeta\right)^2\end{aligned}$$

$$\begin{aligned}&\frac{d\beta}{ds}\left(\bar{V}(\zeta,\xi_1)\right)\left(1+|\theta|\right)\mu_1(\zeta,\xi_1)\left|\xi_1-\tilde{k}'\zeta\right||d|\sqrt{R(w)+\Lambda}\\
&\leq\frac{a}{2\rho}\left(R(w)+\Lambda\right)+\frac{\rho}{2a}\left(\frac{d\beta}{ds}\left(\bar{V}(\zeta,\xi_1)\right)\right)^2\mu_1^2(\zeta,\xi_1)\left(1+|\theta|\right)^2|d|^2\left(\xi_1-\tilde{k}'\zeta\right)^2\end{aligned}$$

implies the following inequality:

$$\frac{d\beta}{ds}\left(\bar{V}(\zeta,\xi_1)\right)\left|\xi_1-\tilde{k}'\zeta\right|\left|h_1(w,\zeta,\xi_1,\theta,d)\right| \le \frac{d\beta}{ds}\left(\bar{V}(\zeta,\xi_1)\right)\mu_1(\zeta,\xi_1)\left(1+|\theta|\right)\left|\xi_1-\tilde{k}'\zeta\right|\left|(\zeta,\xi_1)\right|$$
$$+\frac{a}{\rho}\left(R(w)+\Lambda\right)+\frac{\rho}{2a}\left(\frac{d\beta}{ds}\left(\bar{V}(\zeta,\xi_1)\right)\right)^2\mu_1^2(\zeta,\xi_1)\left(1+|\theta|\right)^2\left(\xi_1-\tilde{k}'\zeta\right)^2$$
$$+\frac{d\beta}{ds}\left(\bar{V}(\zeta,\xi_1)\right)\left(1+|\theta|\right)\mu_1(\zeta,\xi_1)\left|\xi_1-\tilde{k}'\zeta\right||d| \tag{5.109}$$
$$+\frac{\rho}{2a}\left(\frac{d\beta}{ds}\left(\bar{V}(\zeta,\xi_1)\right)\right)^2\mu_1^2(\zeta,\xi_1)\left(1+|\theta|\right)^2|d|^2\left(\xi_1-\tilde{k}'\zeta\right)^2$$

Using the inequalities $|\theta| \le \rho+\left(|\theta|-\rho\right)^+$, $\left(1+|\theta|\right)^2 \le 2(1+\rho)^2+2\left(\left(|\theta|-\rho\right)^+\right)^2$, we get from (5.109):

$$\frac{d\beta}{ds}\left(\bar{V}(\zeta,\xi_1)\right)\left|\xi_1-\tilde{k}'\zeta\right|\left|h_1(w,\zeta,\xi_1,\theta,d)\right| \le \frac{d\beta}{ds}\left(\bar{V}(\zeta,\xi_1)\right)\mu_1(\zeta,\xi_1)(1+\rho)\left|\xi_1-\tilde{k}'\zeta\right|\left|(\zeta,\xi_1)\right|$$
$$+\frac{d\beta}{ds}\left(\bar{V}(\zeta,\xi_1)\right)\mu_1(\zeta,\xi_1)\left(|\theta|-\rho\right)^+\left|\xi_1-\tilde{k}'\zeta\right|\left|(\zeta,\xi_1)\right|$$
$$+\frac{a}{\rho}\left(R(w)+\Lambda\right)+\frac{\rho}{a}\left(\frac{d\beta}{ds}\left(\bar{V}(\zeta,\xi_1)\right)\right)^2\mu_1^2(\zeta,\xi_1)(1+\rho)^2\left(\xi_1-\tilde{k}'\zeta\right)^2$$
$$+\frac{\rho}{a}\left(\frac{d\beta}{ds}\left(\bar{V}(\zeta,\xi_1)\right)\right)^2\mu_1^2(\zeta,\xi_1)\left(\left(|\theta|-\rho\right)^+\right)^2\left(\xi_1-\tilde{k}'\zeta\right)^2$$
$$+\frac{d\beta}{ds}\left(\bar{V}(\zeta,\xi_1)\right)(1+\rho)\mu_1(\zeta,\xi_1)\left|\xi_1-\tilde{k}'\zeta\right||d|+\frac{d\beta}{ds}\left(\bar{V}(\zeta,\xi_1)\right)\left(|\theta|-\rho\right)^+\mu_1(\zeta,\xi_1)\left|\xi_1-\tilde{k}'\zeta\right||d|$$
$$+\frac{\rho}{a}\left(\frac{d\beta}{ds}\left(\bar{V}(\zeta,\xi_1)\right)\right)^2\mu_1^2(\zeta,\xi_1)(1+\rho)^2|d|^2\left(\xi_1-\tilde{k}'\zeta\right)^2$$
$$+\frac{\rho}{a}\left(\frac{d\beta}{ds}\left(\bar{V}(\zeta,\xi_1)\right)\right)^2\mu_1^2(\zeta,\xi_1)\left(\left(|\theta|-\rho\right)^+\right)^2|d|^2\left(\xi_1-\tilde{k}'\zeta\right)^2 \tag{5.110}$$

Estimate (5.110) and the inequalities

$$\frac{d\beta}{ds}\left(\bar{V}(\zeta,\xi_1)\right)(1+\rho)\mu_1(\zeta,\xi_1)\left|\xi_1-\tilde{k}'\zeta\right||d|$$
$$\le \frac{a}{2\rho}|d|^2+\frac{\rho}{2a}\left(\frac{d\beta}{ds}\left(\bar{V}(\zeta,\xi_1)\right)\right)^2(1+\rho)^2\mu_1^2(\zeta,\xi_1)\left(\xi_1-\tilde{k}'\zeta\right)^2$$

$$\frac{d\beta}{ds}\left(\bar{V}(\zeta,\xi_1)\right)\left(|\theta|-\rho\right)^+\mu_1(\zeta,\xi_1)\left|\xi_1-\tilde{k}'\zeta\right||d|$$
$$\le \frac{a}{2\rho}|d|^2+\frac{\rho}{2a}\left(\frac{d\beta}{ds}\left(\bar{V}(\zeta,\xi_1)\right)\right)^2\left(\left(|\theta|-\rho\right)^+\right)^2\mu_1^2(\zeta,\xi_1)\left(\xi_1-\tilde{k}'\zeta\right)^2$$

$$\frac{\rho}{a}\left(\frac{d\beta}{ds}(\bar{V}(\zeta,\xi_1))\right)^2 \mu_1^2(\zeta,\xi_1)(1+\rho)^2|d|^2\left(\xi_1-\tilde{k}'\zeta\right)^2$$
$$\leq \frac{a}{2\rho}|d|^4+\frac{\rho^3}{2a^3}\left(\frac{d\beta}{ds}(\bar{V}(\zeta,\xi_1))\right)^4(1+\rho)^4\mu_1^4(\zeta,\xi_1)\left(\xi_1-\tilde{k}'\zeta\right)^4$$

$$\frac{\rho}{a}\left(\frac{d\beta}{ds}(\bar{V}(\zeta,\xi_1))\right)^2 \mu_1^2(\zeta,\xi_1)\left(\left(|\theta|-\rho\right)^+\right)^2|d|^2\left(\xi_1-\tilde{k}'\zeta\right)^2$$
$$\leq \frac{a}{2\rho}|d|^4+\frac{\rho^3}{2a^3}\left(\frac{d\beta}{ds}(\bar{V}(\zeta,\xi_1))\right)^4\left(\left(|\theta|-\rho\right)^+\right)^4\mu_1^4(\zeta,\xi_1)\left(\xi_1-\tilde{k}'\zeta\right)^4$$

give:

$$\begin{aligned}
&\frac{d\beta}{ds}(\bar{V}(\zeta,\xi_1))\left|\xi_1-\tilde{k}'\zeta\right|\left|h_1(w,\zeta,\xi_1,\theta,d)\right| \leq \frac{d\beta}{ds}(\bar{V}(\zeta,\xi_1))\,\mu_1(\zeta,\xi_1)(1+\rho)\left|\xi_1-\tilde{k}'\zeta\right|\left|(\zeta,\xi_1)\right| \\
&+\frac{d\beta}{ds}(\bar{V}(\zeta,\xi_1))\,\mu_1(\zeta,\xi_1)\left(|\theta|-\rho\right)^+\left|\xi_1-\tilde{k}'\zeta\right|\left|(\zeta,\xi_1)\right|+\frac{a}{\rho}\left(|d|^2+|d|^4+R(w)+\Lambda\right) \\
&+\frac{3\rho}{a}\left(\frac{d\beta}{ds}(\bar{V}(\zeta,\xi_1))\right)^2\mu_1^2(\zeta,\xi_1)\left(\left(|\theta|-\rho\right)^+\right)^2\left(\xi_1-\tilde{k}'\zeta\right)^2 \\
&+\frac{3\rho}{2a}\left(\frac{d\beta}{ds}(\bar{V}(\zeta,\xi_1))\right)^2(1+\rho)^2\mu_1^2(\zeta,\xi_1)\left(\xi_1-\tilde{k}'\zeta\right)^2 \\
&+\frac{\rho^3}{2a^3}\left(\frac{d\beta}{ds}(\bar{V}(\zeta,\xi_1))\right)^4(1+\rho)^4\mu_1^4(\zeta,\xi_1)\left(\xi_1-\tilde{k}'\zeta\right)^4 \\
&+\frac{\rho^3}{2a^3}\left(\frac{d\beta}{ds}(\bar{V}(\zeta,\xi_1))\right)^4\left(\left(|\theta|-\rho\right)^+\right)^4\mu_1^4(\zeta,\xi_1)\left(\xi_1-\tilde{k}'\zeta\right)^4
\end{aligned} \tag{5.111}$$

Next, we employ the inequalities

$$\left|(\zeta,\xi_1)\right|\leq|\zeta|+|\xi_1|\leq|\zeta|+\left|\xi_1-\tilde{k}'\zeta\right|+\left|\tilde{k}'\zeta\right|\leq\left(1+\left|\tilde{k}\right|\right)|\zeta|+\left|\xi_1-\tilde{k}'\zeta\right|$$

$$\frac{d\beta}{ds}(\bar{V}(\zeta,\xi_1))\,\mu_1(\zeta,\xi_1)\left(|\theta|-\rho\right)^+\left|\xi_1-\tilde{k}'\zeta\right|\left|(\zeta,\xi_1)\right|$$
$$\leq\frac{a}{\rho}\left(\left(|\theta|-\rho\right)^+\right)^2+\frac{\rho}{4a}\left(\frac{d\beta}{ds}(\bar{V}(\zeta,\xi_1))\right)^2\mu_1^2(\zeta,\xi_1)\left|(\zeta,\xi_1)\right|^2\left(\xi_1-\tilde{k}'\zeta\right)^2$$

$$\frac{3\rho}{a}\left(\frac{d\beta}{ds}(\bar{V}(\zeta,\xi_1))\right)^2\mu_1^2(\zeta,\xi_1)\left(\left(|\theta|-\rho\right)^+\right)^2\left(\xi_1-\tilde{k}'\zeta\right)^2$$
$$\leq\frac{a}{\rho}\left(\left(|\theta|-\rho\right)^+\right)^4+\frac{9\rho^3}{4a^3}\left(\frac{d\beta}{ds}(\bar{V}(\zeta,\xi_1))\right)^4\mu_1^4(\zeta,\xi_1)\left(\xi_1-\tilde{k}'\zeta\right)^4$$

$$\frac{\rho^3}{2a^3}\left(\frac{d\beta}{ds}\left(\bar{V}(\zeta,\xi_1)\right)\right)^4\left(\left(|\theta|-\rho\right)^+\right)^4\mu_1^4(\zeta,\xi_1)\left(\xi_1-\tilde{k}'\zeta\right)^4$$
$$\leq\frac{a}{\rho}\left(\left(|\theta|-\rho\right)^+\right)^8+\frac{\rho^7}{16a^7}\left(\frac{d\beta}{ds}\left(\bar{V}(\zeta,\xi_1)\right)\right)^8\mu_1^8(\zeta,\xi_1)\left(\xi_1-\tilde{k}'\zeta\right)^8$$

and we obtain from (5.111) the following estimate:

$$\begin{aligned}
&\frac{d\beta}{ds}\left(\bar{V}(\zeta,\xi_1)\right)\left|\xi_1-\tilde{k}'\zeta\right|\left|h_1(w,\zeta,\xi_1,\theta,d)\right| \\
&\leq\frac{d\beta}{ds}\left(\bar{V}(\zeta,\xi_1)\right)\mu_1(\zeta,\xi_1)(1+\rho)\left(1+\left|\tilde{k}\right|\right)|\zeta|\left|\xi_1-\tilde{k}'\zeta\right| \\
&+\frac{d\beta}{ds}\left(\bar{V}(\zeta,\xi_1)\right)\mu_1(\zeta,\xi_1)(1+\rho)\left(\xi_1-\tilde{k}'\zeta\right)^2 \\
&+\frac{\rho}{2a}\left(\frac{1}{2}\left|(\zeta,\xi_1)\right|^2+3(1+\rho)^2\right)\left(\frac{d\beta}{ds}\left(\bar{V}(\zeta,\xi_1)\right)\right)^2\mu_1^2(\zeta,\xi_1)\left(\xi_1-\tilde{k}'\zeta\right)^2 \\
&+\frac{a}{\rho}\left(\left(\left(|\theta|-\rho\right)^+\right)^2+\left(\left(|\theta|-\rho\right)^+\right)^4+\left(\left(|\theta|-\rho\right)^+\right)^8+|d|^2+|d|^4+R(w)+\Lambda\right) \\
&+\frac{\rho^3}{2a^3}\left(\frac{9}{2}+(1+\rho)^4\right)\left(\frac{d\beta}{ds}\left(\bar{V}(\zeta,\xi_1)\right)\right)^4\mu_1^4(\zeta,\xi_1)\left(\xi_1-\tilde{k}'\zeta\right)^4 \\
&+\frac{\rho^7}{16a^7}\left(\frac{d\beta}{ds}\left(\bar{V}(\zeta,\xi_1)\right)\right)^8\mu_1^8(\zeta,\xi_1)\left(\xi_1-\tilde{k}'\zeta\right)^8
\end{aligned}\tag{5.112}$$

Combining (5.108) and (5.112) we get the estimate:

$$\dot{V}_1\Big|_{\xi_2=k_1(\zeta,\xi_1,z)} \leq -\frac{r}{2}V_1(\zeta,\xi_1)+\frac{r}{2}V_1(\zeta,\xi_1)-\frac{3r}{4}\frac{d\beta}{d\,s}\left(\bar{V}(\zeta,\xi_1)\right)\zeta' P\zeta - Q(\zeta,\xi_1)$$
$$+\frac{d\beta}{d\,s}\left(\bar{V}(\zeta,\xi_1)\right)\left(2\left|P\bar{b}\right|+\left|\tilde{k}'\left(A+\bar{b}\tilde{k}'\right)\right|+\mu_1(\zeta,\xi_1)(1+\rho)\left(1+\left|\tilde{k}\right|\right)\right)|\zeta|\left|\xi_1-\tilde{k}'\zeta\right|$$
$$+\frac{d\beta}{d\,s}\left(\bar{V}(\zeta,\xi_1)\right)\left(\mu_1(\zeta,\xi_1)(1+\rho)+\frac{r}{8}+\left|\tilde{k}'\bar{b}\right|\right)\left(\xi_1-\tilde{k}'\zeta\right)^2$$
$$+\frac{\rho}{2a}\left(\frac{1}{2}|(\zeta,\xi_1)|^2+3(1+\rho)^2\right)\left(\frac{d\beta}{d\,s}\left(\bar{V}(\zeta,\xi_1)\right)\right)^2\mu_1^2(\zeta,\xi_1)\left(\xi_1-\tilde{k}'\zeta\right)^2$$
$$+\frac{a}{\rho}\left(\left(\left(|\theta|-\rho\right)^+\right)^2+\left(\left(|\theta|-\rho\right)^+\right)^4+\left(\left(|\theta|-\rho\right)^+\right)^8+|d|^2+|d|^4+R(w)+\Lambda\right)$$
$$+\frac{\rho^3}{2a^3}\left(\frac{9}{2}+(1+\rho)^4\right)\left(\frac{d\beta}{d\,s}\left(\bar{V}(\zeta,\xi_1)\right)\right)^4\mu_1^4(\zeta,\xi_1)\left(\xi_1-\tilde{k}'\zeta\right)^4$$
$$+\frac{\rho^7}{16a^7}\left(\frac{d\beta}{d\,s}\left(\bar{V}(\zeta,\xi_1)\right)\right)^8\mu_1^8(\zeta,\xi_1)\left(\xi_1-\tilde{k}'\zeta\right)^8$$
$$-b_1\rho S(\zeta,\xi_1,z)\frac{S(\zeta,\xi_1,z)\left(\xi_1-\tilde{k}'\zeta\right)^2+4a}{4a}\left(\xi_1-\tilde{k}'\zeta\right)^2$$

(5.113)

Definitions (5.90), (5.99) and estimate (5.100) in conjunction with (5.113) imply the estimate:

$$\dot{V}_1\Big|_{\xi_2=k_1(\zeta,\xi_1,z)} \leq -\frac{r}{2}V_1(\zeta,\xi_1)-\frac{r}{4}\frac{d\beta}{d\,s}\left(\bar{V}(\zeta,\xi_1)\right)\zeta' P\zeta - Q(\zeta,\xi_1)$$
$$+\frac{d\beta}{d\,s}\left(\bar{V}(\zeta,\xi_1)\right)\left(2\left|P\bar{b}\right|+\left|\tilde{k}'\left(A+\bar{b}\tilde{k}'\right)\right|+\mu_1(\zeta,\xi_1)(1+\rho)\left(1+\left|\tilde{k}\right|\right)\right)|\zeta|\left|\xi_1-\tilde{k}'\zeta\right|$$
$$+\frac{d\beta}{d\,s}\left(\bar{V}(\zeta,\xi_1)\right)\left(\mu_1(\zeta,\xi_1)(1+\rho)+\frac{3r}{8}+\left|\tilde{k}'\bar{b}\right|\right)\left(\xi_1-\tilde{k}'\zeta\right)^2$$
$$+\frac{\rho}{2a}\left(\frac{1}{2}|(\zeta,\xi_1)|^2+3(1+\rho)^2\right)\left(\frac{d\beta}{d\,s}\left(\bar{V}(\zeta,\xi_1)\right)\right)^2\mu_1^2(\zeta,\xi_1)\left(\xi_1-\tilde{k}'\zeta\right)^2$$
$$+\frac{a}{\rho}\left(\left(\left(|\theta|-\rho\right)^+\right)^2+\left(\left(|\theta|-\rho\right)^+\right)^4+\left(\left(|\theta|-\rho\right)^+\right)^8+|d|^2+|d|^4+R(w)+\Lambda\right)$$
$$+\frac{\rho^3}{2a^3}\left(\frac{9}{2}+(1+\rho)^4\right)\left(\frac{d\beta}{d\,s}\left(\bar{V}(\zeta,\xi_1)\right)\right)^4\mu_1^4(\zeta,\xi_1)\left(\xi_1-\tilde{k}'\zeta\right)^4$$
$$+\frac{\rho^7}{16a^7}\left(\frac{d\beta}{d\,s}\left(\bar{V}(\zeta,\xi_1)\right)\right)^8\mu_1^8(\zeta,\xi_1)\left(\xi_1-\tilde{k}'\zeta\right)^8$$
$$-b_1\rho S(\zeta,\xi_1,z)\frac{S(\zeta,\xi_1,z)\left(\xi_1-\tilde{k}'\zeta\right)^2+4a}{4a}\left(\xi_1-\tilde{k}'\zeta\right)^2$$

(5.114)

Combining (5.114), (5.103), (5.91) (which implies that $\zeta' P\zeta \geq \bar{c}|\zeta|^2$) and the inequality

$$\left(2\left|P\bar{b}\right|+\left|\tilde{k}'\left(A+\bar{b}\tilde{k}'\right)\right|+\mu_1(\zeta,\xi_1)(1+\rho)\left(1+\left|\tilde{k}\right|\right)\right)|\zeta|\left|\xi_1-\tilde{k}'\zeta\right|$$
$$\leq \frac{r}{4}\bar{c}|\zeta|^2+\frac{1}{r\bar{c}}\left(2\left|P\bar{b}\right|+\left|\tilde{k}'\left(A+\bar{b}\tilde{k}'\right)\right|+\mu_1(\zeta,\xi_1)(1+\rho)\left(1+\left|\tilde{k}\right|\right)\right)^2\left(\xi_1-\tilde{k}'\zeta\right)^2$$

we obtain the estimate:

$$\begin{aligned}&\dot{V}_1\Big|_{\xi_2=k_1(\zeta,\xi_1,z)} \leq -\frac{r}{2}V_1(\zeta,\xi_1)-Q(\zeta,\xi_1)+S(\zeta,\xi_1,z)\left(\xi_1-\tilde{k}'\zeta\right)^2\\&+\frac{a}{\rho}\left(\left(\left(|\theta|-\rho\right)^+\right)^2+\left(\left(|\theta|-\rho\right)^+\right)^4+\left(\left(|\theta|-\rho\right)^+\right)^8+|d|^2+|d|^4+R(w)+\Lambda\right)\\&-b_1\rho S(\zeta,\xi_1,z)\frac{S(\zeta,\xi_1,z)\left(\xi_1-\tilde{k}'\zeta\right)^2+4a}{4a}\left(\xi_1-\tilde{k}'\zeta\right)^2\end{aligned} \tag{5.115}$$

Inequality (5.85) with $c=r/2$ is a consequence of (5.115) and the inequalities

$$\begin{aligned}&S(\zeta,\xi_1,z)\left(\xi_1-\tilde{k}'\zeta\right)^2=b_1\frac{1}{b_1}S(\zeta,\xi_1,z)\left(\xi_1-\tilde{k}'\zeta\right)^2\\&=b_1\rho S(\zeta,\xi_1,z)\left(\xi_1-\tilde{k}'\zeta\right)^2+b_1\left(\frac{1}{b_1}-\rho\right)S(\zeta,\xi_1,z)\left(\xi_1-\tilde{k}'\zeta\right)^2\\&\leq b_1\rho S(\zeta,\xi_1,z)\left(\xi_1-\tilde{k}'\zeta\right)^2+b_1\left(\frac{1}{b_1}-\rho\right)^+S(\zeta,\xi_1,z)\left(\xi_1-\tilde{k}'\zeta\right)^2\\&\leq b_1\rho S(\zeta,\xi_1,z)\left(\xi_1-\tilde{k}'\zeta\right)^2+\frac{a}{\rho}b_1\left(\left(\frac{1}{b_1}-\rho\right)^+\right)^2+b_1\frac{\rho}{4a}S^2(\zeta,\xi_1,z)\left(\xi_1-\tilde{k}'\zeta\right)^4\end{aligned}$$

<u>Part III:</u> Let arbitrary $\kappa,\Gamma,\varepsilon>0$ be given. Define $a=\frac{\kappa}{2\left(1+2^{m-1}\right)}$. By virtue of Part II there exist functions $k,V\in C^\infty\left(\mathbb{R}^n\times\mathbb{R}^m\times\mathbb{R}\right)$ with $V(0,0,z)=k(0,0,z)=0$ for $z\in\mathbb{R}$, $V(\zeta,\xi,z)>0$ for $(\zeta,\xi)\neq(0,0)$ and a constant $c>0$ such that:

(a) for each $\tilde{M}>0$ the set $\left\{(\zeta,\xi,z)\in\mathbb{R}^n\times\mathbb{R}^m\times\left[-\tilde{M},\tilde{M}\right]:V(\zeta,\xi,z)\leq\tilde{M}\right\}$ is bounded, and

(b) inequality (5.87) as well as the following inequality holds:

$$\frac{\partial V}{\partial \zeta}(\zeta,\xi,z)\left(A\zeta+b\xi_1\right)+\sum_{i=1}^{m-1}\frac{\partial V}{\partial \xi_i}(\zeta,\xi,z)\left(h_j(w,\zeta,\xi_1,...,\xi_j,\theta,d)+b_j g_j(w,\zeta,\xi_1,...,\xi_j,\theta)\xi_{j+1}\right)$$
$$+\frac{\partial V}{\partial \xi_m}(\zeta,\xi,z)\left(h_m(w,\zeta,\xi,\theta,d)+b_m g_m(w,\zeta,\xi,\theta)k\left(\zeta,\xi,z\right)\right)$$
$$+\frac{\partial V}{\partial z}(\zeta,\xi,z)\Gamma\exp(-z)\left(V(\zeta,\xi,z)-\varepsilon\right)^+$$
$$\leq -2^{1-m}cV(\zeta,\xi,z)-Q(\zeta,\xi_1)+\frac{2^{m-1}\kappa}{2\rho\left(1+2^{m-1}\right)}\left(\min_{i=1,...,m}\left(b_i\right)\left(\left(\frac{1}{\min_{i=1,...,m}\left(b_i\right)}-\rho\right)^+\right)^2+\left|d\right|^2+\left|d\right|^4+R(w)\right)$$
$$+\frac{2^{m-1}\kappa}{2\rho\left(1+2^{m-1}\right)}\left(\left(\left(|\theta|+|b|+|\theta||b|-\rho\right)^+\right)^2+\left(\left(|\theta|+|b|+|\theta||b|-\rho\right)^+\right)^4+\left(\left(|\theta|+|b|+|\theta||b|-\rho\right)^+\right)^8+\Lambda\right)$$

for all $\left(\zeta,\xi,w\right)\in\mathbb{R}^n\times\mathbb{R}^m\times\mathbb{R}^l$, $z\in\mathbb{R}$, $d\in\mathbb{R}^q$, $\theta\in\Theta$, $b_j\in\left(0,+\infty\right)\cap B_j$ for $j=1,...,m$ (5.116)

where $\rho:=\kappa+\exp(z)$.

We next perform the stability analysis of the closed-loop system (2.3), (2.4), (2.5) with (2.24), (2.25). To this end, we need some facts.

First, the fact that the smooth functions $\Phi,R:\mathbb{R}^l\to\mathbb{R}_+$ are positive definite and radially unbounded, functions and Proposition 2.2 on page 107 in [7] implies the existence of functions $a_1,a_2,a_3,a_4\in K_\infty$ such that inequalities (5.63), (5.64) hold for all $w\in\mathbb{R}^l$.

Secondly, we notice the existence of a function $a_5\in K_\infty$ such that the following inequality is valid for all $w\in\mathbb{R}^l$, $\xi\in\mathbb{R}^m$, $\zeta\in\mathbb{R}^n$ and $z\in\mathbb{R}$:

$$a_5\left(\Phi\left(w\right)+V\left(\zeta,\xi,z\right)\right)\leq\frac{1}{2}R\left(w\right)+2^{1-m}cV\left(\zeta,\xi,z\right) \tag{5.117}$$

Indeed, (5.117) holds with $a_5\left(s\right)=\min\left(\frac{1}{2}a_1\left(\frac{s}{2}\right),2^{-m}cs\right)$ and notice that $a_5\in K_\infty$. To see why (5.107) is valid notice that since $a_5\left(s_1+s_2\right)\leq a_5\left(2s_1\right)+a_5\left(2s_2\right)$ for all $s_1,s_2\geq 0$ we get:

$$\begin{aligned}&a_5\left(\Phi\left(w\right)+V\left(\zeta,\xi,z\right)\right)\leq a_5\left(2\Phi\left(w\right)\right)+a_5\left(2V\left(\zeta,\xi,z\right)\right)\\&=\min\left(\frac{1}{2}a_1\left(\Phi\left(w\right)\right),2^{1-m}c\Phi\left(w\right)\right)+\min\left(\frac{1}{2}a_1\left(V\left(\zeta,\xi,z\right)\right),2^{1-m}cV\left(\zeta,\xi,z\right)\right)\\&\leq\frac{1}{2}a_1\left(\Phi\left(w\right)\right)+2^{1-m}cV\left(\zeta,\xi,z\right)\end{aligned}$$

and inequality (5.117) is a consequence of the left inequality (5.63).

Thirdly, we notice that there exist functions $\omega \in KL$, $\bar{\zeta} \in K_\infty$ with the following property: estimate (5.67) holds for all $t \in [0,T)$, for every $T \in (0,+\infty]$, $\tilde{\mu} \geq 0$ and for every absolutely continuous function $Y:[0,T) \to \mathbb{R}_+$ for which inequality (5.66) holds for $t \in [0,T)$ a.e..

Fourthly, we notice the existence of a non-decreasing, continuous function $H:\mathbb{R}_+ \to \mathbb{R}_+$ such that the following inequality is valid for all $\xi \in \mathbb{R}^m$, $\zeta \in \mathbb{R}^n$ and $z \in \mathbb{R}$:

$$|(\zeta,\xi)| \leq H\left(\max\left(V\left(\zeta,\xi,z\right),|z|\right)\right) \tag{5.118}$$

The proof of (5.118) is exactly the same as the proof of (5.68) in the proof of Theorem 1.

We define the function:

$$U\left(w,\zeta,\xi,z\right) = \Phi\left(w\right) + V\left(\zeta,\xi,z\right) \tag{5.119}$$

Let arbitrary $\left(w_0,\zeta_0,\xi_0,z_0\right) \in \mathbb{R}^l \times \mathbb{R}^n \times \mathbb{R}^m \times \mathbb{R}$, $d \in L^\infty\left(\mathbb{R}_+;\mathbb{R}^l\right)$, $\theta \in L^\infty\left(\mathbb{R}_+;\Theta\right)$, $\delta \in L^\infty\left(\mathbb{R}_+;\mathbb{R}^\alpha\right)$, $b_j \in L^\infty\left(\mathbb{R}_+;(0,+\infty) \cap B_j\right)$ with $\inf_{t\geq 0}\left(b_j(t)\right) > 0$ for $j=1,...,m$ be given. Consider the unique solution $\left(w(t),\zeta(t),\xi(t),z(t)\right)$ of (2.3), (2.4), (2.5) with (2.24), (2.25) and $\left(w(0),\zeta(0),\xi(0),z(0)\right) = \left(w_0,\zeta_0,\xi_0,z_0\right)$ for which there exists $t_{\max} \in (0,+\infty]$ such that solution is defined on $[0,t_{\max})$. Moreover, if $t_{\max} < +\infty$ then $\limsup_{t\to t_{\max}^-}\left(|w(t)|+|\zeta(t)|+|\xi(t)|+|z(t)|\right) = +\infty$.

By virtue of (2.25) the map $[0,t_{\max}) \ni t \mapsto z(t)$ is non-decreasing. Therefore, we get for all $t \in [0,t_{\max})$:

$$z_0 \leq z(t) \tag{5.120}$$

Define for $t \geq 0$ a.e.:

$$\begin{aligned} v(t) &:= \gamma\left(|\delta(t)|\right) + \bar{\gamma}\left(|\theta(t)|\right) + \frac{1}{2}\left(\left(|\theta(t)|+|b(t)|+|\theta(t)||b(t)|-\kappa-\exp\left(z(t)\right)\right)^+\right)^2 \\ &+\frac{1}{2}\left(|d(t)|^2+|d(t)|^4+\min_{i=1,...,m}\left(b_i(t)\right)\left(\left(\left(\frac{1}{\min_{i=1,...,m}\left(b_i(t)\right)}-\kappa-\exp\left(z(t)\right)\right)^+\right)^2+\Lambda\right) \\ &+\frac{1}{2}\left(\left(|\theta(t)|+|b(t)|+|\theta(t)||b(t)|-\kappa-\exp\left(z(t)\right)\right)^+\right)^4 \\ &+\frac{1}{2}\left(\left(|\theta(t)|+|b(t)|+|\theta(t)||b(t)|-\kappa-\exp\left(z(t)\right)\right)^+\right)^8 \end{aligned} \tag{5.121}$$

We notice that since $\gamma,\bar{\gamma} \in C^0\left(\mathbb{R}_+;\mathbb{R}_+\right)$ are non-decreasing functions and since $d \in L^\infty\left(\mathbb{R}_+;\mathbb{R}^l\right)$, $\theta \in L^\infty\left(\mathbb{R}_+;\Theta\right)$, $\delta \in L^\infty\left(\mathbb{R}_+;\mathbb{R}^\alpha\right)$, $b_j \in L^\infty\left(\mathbb{R}_+;(0,+\infty) \cap B_j\right)$ with $\inf_{t\geq 0}\left(b_j(t)\right) > 0$ for $j=1,...,m$, definition (5.121) implies that $v \in L^\infty\left(\mathbb{R}_+;\mathbb{R}\right)$ with

$$\|v\|_\infty \le G\left(\|d\|_\infty, \|\delta\|_\infty, \|\theta\|_\infty, \|b\|_\infty, \min_{i=1,\dots,m}\left(\inf_{s\ge 0}\left(b_i(s)\right)\right)\right)$$
$$:= \gamma\left(\|\delta\|_\infty\right) + \bar{\gamma}\left(\|\theta\|_\infty\right) + \frac{1}{2}\left(\|d\|_\infty^2 + \|d\|_\infty^4 + \left(\|\theta\|_\infty + \|b\|_\infty + \|\theta\|_\infty \|b\|_\infty\right)^2 + \Lambda\right) \tag{5.122}$$
$$+\frac{1}{2}\left(\left(\|\theta\|_\infty + \|b\|_\infty + \|\theta\|_\infty \|b\|_\infty\right)^4 + \left(\|\theta\|_\infty + \|b\|_\infty + \|\theta\|_\infty \|b\|_\infty\right)^8 + \frac{1}{\min_{i=1,\dots,m}\left(\inf_{s\ge 0}\left(b_i(s)\right)\right)}\right)$$

Inequalities (2.15), (5.116) and definitions (5.119), (5.121) imply that the following differential inequality holds for $t \in [0, t_{\max})$ a.e.:

$$\frac{d}{dt}\left(U\left(w(t), \zeta(t), \xi(t), z(t)\right)\right) \le -\frac{1}{2} R(w(t)) - 2^{1-m} c V(\zeta(t), \xi(t), z(t)) + v(t) \tag{5.123}$$

Combining (5.117), (5.119), (5.66), (5.67) and (5.122), (5.123) we obtain for all $t \in [0, t_{\max})$:

$$U\left(w(t), \zeta(t), \xi(t), z(t)\right)$$
$$\le \omega\left(U\left(w_0, \zeta_0, \xi_0, z_0\right), t\right) + \bar{\zeta}\left(G\left(\|d\|_\infty, \|\delta\|_\infty, \|\theta\|_\infty, \|b\|_\infty, \min_{i=1,\dots,m}\left(\inf_{s\ge 0}\left(b_i(s)\right)\right)\right)\right) \tag{5.124}$$

Inequalities (5.116), (5.124), the right inequality (5.63) and definition (5.119) imply that the following differential inequality holds for $t \in [0, t_{\max})$ a.e.:

$$\frac{d}{dt}\left(V\left(\zeta(t), \xi(t), z(t)\right)\right) \le -2^{1-m} c V\left(\zeta(t), \xi(t), z(t)\right)$$
$$+\frac{Z\left(w_0, \zeta_0, \xi_0, z_0, \|d\|_\infty, \|\delta\|_\infty, \|\theta\|_\infty, \|b\|_\infty, \min_{i=1,\dots,m}\left(\inf_{s\ge 0}\left(b_i(s)\right)\right)\right)}{1 + \exp\left(z(t)\right)} \tag{5.125}$$

where

$$Z\left(w_0, \zeta_0, \xi_0, z_0, \|d\|_\infty, \|\delta\|_\infty, \|\theta\|_\infty, \|b\|_\infty, \min_{i=1,\dots,m}\left(\inf_{s\ge 0}\left(b_i(s)\right)\right)\right)$$
$$:= \frac{\kappa}{2\min(1,\kappa)} a_2\left(\omega\left(U\left(w_0, \zeta_0, \xi_0, z_0\right), 0\right) + \bar{\zeta}\left(G\left(\|d\|_\infty, \|\delta\|_\infty, \|\theta\|_\infty, \|b\|_\infty, \min_{i=1,\dots,m}\left(\inf_{s\ge 0}\left(b_i(s)\right)\right)\right)\right)\right)$$
$$+\frac{\kappa}{2\min(1,\kappa)}\left(\left(\|\theta\|_\infty + \|b\|_\infty + \|\theta\|_\infty \|b\|_\infty\right)^2 + \|d\|_\infty^2 + \|d\|_\infty^4 + \frac{1}{\min_{i=1,\dots,m}\left(\inf_{s\ge 0}\left(b_i(s)\right)\right)}\right)$$
$$+\frac{\kappa}{2\min(1,\kappa)}\left(\left(\|\theta\|_\infty + \|b\|_\infty + \|\theta\|_\infty \|b\|_\infty\right)^4 + \left(\|\theta\|_\infty + \|b\|_\infty + \|\theta\|_\infty \|b\|_\infty\right)^8 + \Lambda\right)$$
$$\tag{5.126}$$

Lemma 4.1 on page 63 in [10], (2.5) and differential inequality (5.125) imply the following estimate for all $t \in [0, t_{\max})$:

$$z_0 \le z(t) \le \ln\left(\tilde{Z}\left(w_0,\zeta_0,\xi_0,z_0,\|d\|_\infty,\|\delta\|_\infty,\|\theta\|_\infty,\|b\|_\infty,\min_{i=1,\ldots,m}\left(\inf_{s\ge 0}\left(b_i(s)\right)\right)\right)\right) \tag{5.127}$$

where

$$\begin{aligned}
&\tilde{Z}\left(w_0,\zeta_0,\xi_0,z_0,\|d\|_\infty,\|\delta\|_\infty,\|\theta\|_\infty,\|b\|_\infty,\min_{i=1,\ldots,m}\left(\inf_{s\ge 0}\left(b_i(s)\right)\right)\right)\\
&:=\left(\frac{1}{\varepsilon}+\frac{2^{m-1}\Gamma}{c}\right)\frac{2^{m-1}}{c}Z\left(w_0,\zeta_0,\xi_0,z_0,\|d\|_\infty,\|\delta\|_\infty,\|\theta\|_\infty,\|b\|_\infty,\min_{i=1,\ldots,m}\left(\inf_{s\ge 0}\left(b_i(s)\right)\right)\right)\\
&+\exp(z_0)+\frac{2^{m-1}\Gamma}{c}V\left(\zeta_0,\xi_0,z_0\right)
\end{aligned} \tag{5.128}$$

It follows from (5.118), the left inequality (5.64), definition (5.119) and estimates (5.127), (5.124) that the following estimates hold for all $t \in [0,t_{\max})$:

$$|w(t)| \le a_3^{-1}\left(\bar{Z}\right),\quad |\xi(t)| \le H\left(\bar{Z}+K\right),\quad |z(t)| \le K \tag{5.129}$$

where

$$\begin{aligned}
&\bar{Z} := \omega\left(U\left(w_0,\zeta_0,\xi_0,z_0\right),0\right)+\bar{\zeta}\left(G\left(\|d\|_\infty,\|\delta\|_\infty,\|\theta\|_\infty,\|b\|_\infty,\min_{i=1,\ldots,m}\left(\inf_{s\ge 0}\left(b_i(s)\right)\right)\right)\right)\\
&K := |z_0|+\left|\ln\left(\tilde{Z}\left(w_0,\zeta_0,\xi_0,z_0,\|d\|_\infty,\|\delta\|_\infty,\|\theta\|_\infty,\|b\|_\infty,\min_{i=1,\ldots,m}\left(\inf_{s\ge 0}\left(b_i(s)\right)\right)\right)\right)\right|
\end{aligned} \tag{5.130}$$

Estimates (5.129) imply that $t_{\max}=+\infty$. Moreover, estimates (5.129), (5.127) and the fact that $z(t)$ is non-decreasing, imply estimates (2.27), (2.28) for appropriate $\Xi \in C^0\left(\mathbb{R}_+^8 \times (0,+\infty)\right)$.

We next show estimate (2.26). Estimate (5.127) and the fact that $z(t)$ is non-decreasing guarantees that the function $z(t)$ has a finite limit as $t \to +\infty$. This implies that the function $\exp(z(t))$ has a finite limit as $t \to +\infty$. Moreover, the facts that $d \in L^\infty\left(\mathbb{R}_+;\mathbb{R}^l\right)$, $\theta \in L^\infty\left(\mathbb{R}_+;\Theta\right)$, $b_j \in L^\infty\left(\mathbb{R}_+;(0,+\infty)\cap B_j\right)$ with $\inf_{t\ge 0}\left(b_j(t)\right)>0$ for $j=1,\ldots,m$ and (5.129), (2.3), (2.4), (2.24), (2.25) imply that $\frac{d}{dt}\left(V(\zeta(t),\xi(t),z(t))\right)$ is of class $L^\infty\left(\mathbb{R}_+\right)$. It follows that the function $\frac{d}{dt}\left(\exp(z(t))\right)=\Gamma\left(V(\zeta(t),\xi(t),z(t))-\varepsilon\right)^+$ is uniformly continuous, in addition to $\int_0^{+\infty}\frac{d}{dt}\left(\exp(z(t))\right)dt=\lim_{t\to+\infty}\left(\exp(z(t))\right)-\exp(z_0)<+\infty$. From Barbălat's Lemma (see [8, 14]), we have:

$$\lim_{t\to+\infty}\left(\frac{d}{dt}\left(\exp(z(t))\right)\right)=\lim_{t\to+\infty}\left(\Gamma\left(V(\zeta(t),\xi(t),z(t))-\varepsilon\right)^+\right)=0 \tag{5.131}$$

Therefore, the right estimate (2.26) holds. The left estimate (2.26) is a consequence of the right estimate (2.26) and inequality (5.87).

Finally, we consider the case $\gamma(0)=\Lambda=0$, $\bar{\gamma}(s)\equiv 0$, $\lim_{t\to+\infty}\left(d(t)\right)=0$, $\lim_{t\to+\infty}\left(\delta(t)\right)=0$, $\liminf_{t\to+\infty}\left(b_i(t)\right)\geq\frac{1}{\kappa}$ for $i=1,...,m$ and $\kappa+\exp\left(\lim_{s\to+\infty}\left(z(s)\right)\right)\geq\limsup_{t\to+\infty}\left(|\theta(t)|+|b(t)|+|\theta(t)||b(t)|\right)$. In this case, definition (5.121) shows that $\lim_{t\to+\infty}\left(v(t)\right)=0$. Combining (5.117), (5.66), (5.67), (5.119) and (5.123) we obtain:

$$\limsup_{t\to+\infty}\left(U\left(w(t),\zeta(t),\xi(t),z(t)\right)\right)\leq\bar{\zeta}\left(\limsup_{t\to+\infty}\left(v(t)\right)\right) \tag{5.132}$$

It follows from (5.132) and $\lim_{t\to+\infty}\left(v(t)\right)=0$ that $\lim_{t\to+\infty}\left(U\left(w(t),\zeta(t),\xi(t),z(t)\right)\right)=0$. Clearly, the equation $\lim_{t\to+\infty}\left(U\left(w(t),\zeta(t),\xi(t),z(t)\right)\right)=0$ in conjunction with definition (5.119) gives $\lim_{t\to+\infty}\left(\Phi\left(w(t)\right)\right)=\lim_{t\to+\infty}\left(V\left(\zeta(t),\xi(t),z(t)\right)\right)=0$. Inequalities (5.64) and $\lim_{t\to+\infty}\left(\Phi\left(w(t)\right)\right)=0$ show that $\lim_{t\to+\infty}\left(|w(t)|\right)=0$. We next show by contradiction that $\lim_{t\to+\infty}\left(\left|\left(\zeta(t),\xi(t)\right)\right|\right)=0$. Suppose that $\limsup_{t\to+\infty}\left(\left|\left(\zeta(t),\xi(t)\right)\right|\right)>0$. In this case there exists $\bar{\varepsilon}>0$ and a sequence $\left\{t_k\in\mathbb{R}_+ : k=1,...\right\}$ with $\lim_{k\to+\infty}\left(t_k\right)=+\infty$ such that $\left|\left(\zeta\left(t_k\right),\xi\left(t_k\right)\right)\right|\geq\bar{\varepsilon}$ for $k=1,2,...$. By virtue of (5.129) the sequence $\left\{\left(\zeta\left(t_k\right),\xi\left(t_k\right),z\left(t_k\right)\right)\in\mathbb{R}^n\times\mathbb{R}^m\times\mathbb{R} : k=1,...\right\}$ is bounded. Consequently, without loss of generality (by restricting our attention to a subsequence of $\left\{t_k\in\mathbb{R}_+ : k=1,...\right\}$ if necessary), we may assume the existence of $\xi^*\in\mathbb{R}^m$, $\zeta^*\in\mathbb{R}^n$ and $z^*\in\mathbb{R}$ such that $\lim_{k\to+\infty}\left(\xi\left(t_k\right)\right)=\xi^*$, $\lim_{k\to+\infty}\left(\zeta\left(t_k\right)\right)=\zeta^*$ and $\lim_{k\to+\infty}\left(z\left(t_k\right)\right)=z^*$. Clearly, $\left|\left(\zeta^*,\xi^*\right)\right|\geq\bar{\varepsilon}$ and since $\lim_{t\to+\infty}\left(V\left(\zeta(t),\xi(t),z(t)\right)\right)=0$ we get (by continuity of $V$) $V\left(\zeta^*,\xi^*,z^*\right)=0$. However, $V(0,0,z)=0$ for $z\in\mathbb{R}$ and $V(\zeta,\xi,z)>0$ for $\left(\zeta,\xi\right)\neq\left(0,0\right)$, we get $\left(\zeta^*,\xi^*\right)=\left(0,0\right)$; a contradiction with $\left|\left(\zeta^*,\xi^*\right)\right|\geq\bar{\varepsilon}>0$.

The proof is complete. $\triangleleft$

## 6. Concluding Remarks

It is clear that DADS is an adaptive control scheme that can provide robust, partial-state, adaptive feedback laws. However, the extension of DADS to the case where only one output is measured is not trivial: an observer-based control structure is needed. Moreover, despite the fact that comparisons of DADS with leakage methodologies have been performed, a detailed comparison of DADS with other existing adaptive feedback design methodologies is also needed.